# Poc Sets, Median Algebras and Group Actions

An extended Study
of Dunwoody's Construction and Sageev's Theorem

# Poc Sets, Median Algebras and Group Actions
## An extended Study of Dunwoody's Construction and Sageev's Theorem

Martin Roller, Holzstr. 1, 80469 München, Germany
Martin.Roller@dynaware.de

BACKGROUND AND MOTIVATION

A Leitmotiv of geometric group theory is the study of group actions on geometric structures, from which one can draw conclusions about both the group and the geometry. A typical and particularly successful example is the theory of groups acting on trees. Here we can distinguish three different approaches.

- TITS [1970] describes the dynamic behaviour of automorphisms of trees. An automorphism either fixes a vertex or an edge of a tree, or it shifts along an axis. In analogy to classical hyperbolic geometry one can divide groups of automorphisms in three classes.
    (i) Elliptic groups, which leave an edge or a vertex fixed.
    (ii) Hyperbolic groups, which contain a shift.
    (iii) Parabolic groups, which stabilize an end of the tree.

- Bass and SERRE [1977] analyze the algebraic structure of a group acting on a tree in terms of the stabilizers of vertices and edges. For example, if $G$ acts hyperbolically and transitively on the edges of a tree, then $G$ splits algebraically as an amalgamated product or an HNN-extension, with the stabilizer of an edge as the splitting subgroup. In the general case the notion of a graph of groups describes the splitting adequately. This theory makes it desirable find a criterion to decide whether a given group can act hyperbolically on a tree.

- STALLINGS [1969,1971] provides just such a criterion: A finitely generated group splits over a finite subgroup if and only if it has more than one end. Ends of groups where introduced by FREUDENTHAL [1931] and HOPF [1943] as a means of compactification. Later SPECKER [1949] (see also EPSTEIN [1962]) exhibited the cohomological nature of this invariant by showing that an infinite group has more than one end iff $H^1(G, \mathbb{Z}G) \neq 0$. Nowadays we take the latter as a definition.

This successful theory invited many generalizations of the concept of trees. We just briefly indicate a few.

- ℝ-trees (and, more generally, $R$-trees) where introduced by MORGAN-SHALEN [1985]. They turned out to be instrumental for the study of free groups and 3-manifold groups. A deep theorem of Rips classifies the groups acting freely on ℝ-trees.
- DUNWOODY [1996] introduced pro-trees to produce examples of finitely generated, inaccessible groups.
- RIBES-ZALESSKII [1993] used profinite trees to prove hard results about the profinite topology of free groups.
- The hyperbolic spaces of GROMOV [1987] share many qualitative properties with trees, and the theory of word hyperbolic groups is one of the corner stones of geometric group theory.
- The forthcoming book NEUMANN-ADELEKE [1998] studies various tree-like structures and the permutation groups acting on them. BOWDITCH [1995] used tree-like structures to understand the boundary of hyperbolic groups.

Stallings' Theorem itself started a whole series of generalizations and applications (e.g., COHEN [1972], SWAN [1969]). A first break through was made by Dunwoody [1979]. Stallings originally found the algebraic splitting by means of a bipolar structure, which models the normal forms of amalgamated products. Dunwoody reinterpreted Stallings' proof and connected it with the theory of Bass and Serre. Dunwoody's proof consists of two steps.

CAREFUL CHOICE

The fact that $G$ has more than one end is equivalent the existence of a proper almost invariant subset $A \subset G$, i.e., $A + Ag$ is finite for all $g \in G$, but neither $A$ nor $G \smallsetminus A$ is finite. In the first step (which is also present in Stallings' proof), one modifies the set $A$ carefully, until the set $\widetilde{E} := \{gA, gA^* \mid g \in G\}$ of $G$-translates of $A$ and its complement $A^* := G \smallsetminus A$ satisfies the following two conditions.

(a) $\widetilde{E}$ is interval finite, i.e., for all $a, b \in \widetilde{E}$ there exist only finitely many $c \in \widetilde{E}$ with
$$a \subset b \subset c.$$

(b) $\widetilde{E}$ is nested, i.e., for all $a, b \in \widetilde{E}$ one of the following holds
$$a \subseteq b, \qquad a \subseteq b^*, \qquad a^* \subset b, \qquad a^* \subset b^*.$$

DUNWOODY'S CONSTRUCTION

In the second step Dunwoody interprets $\widetilde{E}$ as the set of oriented edges of a $G$-tree $T$. To do this he has to construct a set $V$ of vertices together with incidence relations. He considers the relation
$$a \approx b \quad \Leftrightarrow \quad a = b \quad \text{or} \quad a^* \subset b \text{ and there exists no } c \in \widetilde{E} \text{ with } a^* \subset c \subset b.$$

On a nested set $\widetilde{E}$ this defines an equivalence relation. Now he puts $V := \widetilde{E}/\approx$ and declares the class $[a]_\approx$ to be the terminal vertex of the edge $a$ to define a graph, which is then shown to be a tree, thanks to conditions (a) and (b).

From there Bass-Serre theory takes over to provide an algebraic splitting of $G$.

DICKS-DUNWOODY [1989] is based on an important modification of this construction. For an element $a \in \widetilde{E}$ consider the set

$$U_a := \{b \in \widetilde{E} \mid b \supseteq a \text{ or } b \supset a^*\}.$$

The class $[a]_\approx$ can be recovered as the set of minimal elements of $U_a$. Moreover, $U_a$ has the following two properties.

- For any $b \in \widetilde{E}$ either $b \in U_a$ or $b^* \in U_a$.

- If $b \in U_a$ and $c \supset b$ then $c \in U_a$.

A subset of $\widetilde{E}$ with these two properties is called an ultra filter. A vertex of a tree gives rise to the ultra filter of all oriented edges that point towards it, and these ultra filters form one almost equality class. All other ultra filters are related to ends of the tree.

We record a few results that were possible after this conceptual insight.

- DUNWOODY [1979] used it for a homological characterization of accessible groups, i.e., groups for which the process of splitting repeatedly over finite subgroups stops after finitely many steps.

- MÜLLER [1981] made a careful analysis of relative splittings, which together with the work of Bieri, Eckmann, Linnell and Strebel culminated in the classification of two dimensional Poincaré duality groups.

- The conclusive work about the classical splitting theory is DICKS-DUNWOODY [1989]. Their main effort lies in a refinement of the first step of Stallings' proof. The Almost Stability Theorem they obtain is the most general and powerful result one could wish for.

In contrast, the theory of splittings over infinite subgroups still contains large uncharted territories. Ends of pairs of groups — in a vague sense the analogue to homology of pairs — were already anticipated by SPECKER [1949] and EPSTEIN [1962]. SCOTT [1977] studied the necessary modifications of Dunwoody's proof. One now works with $K$-almost invariant groups, where $K$ is some subgroup of $G$, i.e., $A + Ag$ should be a finite union of right $K$-cosets. One half of Stallings' Theorem carries over easily: If $G$ splits over $K$ then $e(G, K) > 1$. But Scott could prove the converse only under additional assumptions ($G$ and $K$ have to be finitely generated, $K$ should be closed in the profinite topology of $G$), and he gave geometrical examples to show that the converse is not true in general.

KROPHOLLER-ROLLER [1989] examined a very special case, when $G$ and $K$ are Poincaré duality groups of dimension $n$ and $n - 1$, respectively. Here one has $e(G, K) = 2$, essentially. They identified an obstruction (a certain cohomology class) which vanishes precisely if $G$ splits over some subgroup commensurable with $K$.

Stallings explained that he found his proof by contemplating the Sphere Theorem of three dimensional topology. For infinite splitting theory the analogy with low dimensional topology is even stronger and very fruitful. While relative ends detect the algebraic analogue of an immersed codimension one submanifold, splitting theorems find the analogue of an embedded submanifold. KROPHOLLER [1993] gave a completely algebraic proof of the Torus Theorem for 3-manifolds, based on the Kropholler-Roller Splitting Theorem. In a series of very exciting papers DUNWOODY-SWENSON [1997] and DUNWOODY-SAGEEV [1997] managed to translate the Jaco-Shalen-Johannsen Theorem, one of the corner stones of three dimensional topology, into a completely algebraic setting.

Let us consider the two conditions in Dunwoody's Construction once more. Interval finiteness is usually considered "harmless", it is an easy consequence of the almost invariance of $A$, and vice versa. In his theory of finitely generated inaccessible groups Dunwoody initiated the study of nested sets without finite intervals. Applying Dunwoody's Construction one obtains a corresponding set of vertices and incidence relations. The whole object is no longer a tree, but the projective limit of a family of finite trees, hence Dunwoody calls it a pro-tree. Pro-trees include the category of $\mathbb{R}$-trees, but not the pro-finite trees in the sense of Ribes and Zalesskii.

In splitting theory one often has to struggle with the nesting condition. The careful choice procedures require not only to change the set $A$ but also the commensurability class of the subgroup $K$. And the presence of obstructions indicates that sometimes there simply does not exist a good choice. The theorems of splitting theory therefore typically offer an alternative. For example the Torus Theorem says: If a 3-manifold group $G$ contains a $\mathbb{Z} \times \mathbb{Z}$-subgroup, then either $G$ splits (over $\mathbb{Z}$ or $\mathbb{Z} \times \mathbb{Z}$) or there exists a normal $\mathbb{Z}$ in $G$.

If we have neither interval finiteness nor nesting, what is left? A close analysis of splitting theory shows that the arena for most arguments is a very simple structure, which we will call a poc set. It is a partially ordered set $(P, \leq)$ together with a complement operation $^*$, which is an involution that reverses the ordering. Furthermore, for any $a, b \in P$ at most one of the following relations should hold.

$$a \leq b, \qquad a^* \leq b, \qquad a < b^*, \qquad a^* < b^*.$$

Some parts of the theory run more smoothly if we further require a smallest element $0 \in P$ — though in other parts this may seem a nuisance. The above condition can now be expressed by saying

$$\forall a \in P : \qquad a \leq a^* \qquad \Rightarrow \qquad a = 0.$$

Ultra filters still make sense for poc sets, but the full picture doesn't emerge until we find the right structure for them.

A surprising new application of Dunwoody's Construction was found by SAGEEV [1995].

Starting with a group $G$ and a subgroup $K$ satisfying $e(G, K) > 1$ we get a poc set $P$ that is still interval finite but not necessarily nested. Sageev now views the proper elements of $P$ as oriented hyper planes; furthermore, if $a, b \in P$ are nested, then the corresponding hyper planes should be disjoint, otherwise they should meet transversely. This begs the question, in which ambient space these hyper planes should live. Applying Dunwoody's Construction to $P$, Sageev comes up with a set $M$ of "ultra filters", then groups the elements of $M$ to cubes to form a cubical cell complex $\mathscr{K}M$, where the hyper planes live as the extensions of the midplanes of the individual cubes. The subgroup $K$ stabilizes one of the hyper planes, and the full stabilizer of that hyperplane contains $K$ as a subgroup of finite index.

Sageev then shows that $\mathscr{K}M$ enjoys the following two properties.

- Its underlying topological space is simply connected.
- The complex is locally flat.

The second condition can be expressed combinatorially by saying that the link of any vertex is a flag complex. If we endow $\mathscr{K}M$ with a metric which turns every cube into a Euclidean unit cube, then it says that $\mathscr{K}M$ satisfies locally the CAT(0) inequality (see BRIDSON-HAEFLIGER [1995]). Both conditions can be put together neatly by saying that $\mathscr{K}M$ satisfies the global CAT(0) condition. Sageev calls such a CAT(0) cube complex a cubing.

The major part of Sageevs work is the converse of this construction. Starting with a cubing $X$ on which $G$ acts, one can define hyper planes by developing the mid planes of the individual cubes. Sageev then proves that every hyper plane embeds in $X$ and cuts $X$ into two connected components. Sageev shows that for the (oriented) stabilizer $K$ of every essential hyper plane one has $e(G, K) > 1$. Note the adjective "essential" here, we will say more about that later.

The present paper is an extended study of Dunwoody's Construction and Sageev's Theorem. It is based on the observation that Sageevs geometric characterization of cubings has an algebraic counter part. It is a well known fact that in a tree for any three vertices $x$, $y$, $z$ there exists a unique vertex $m(x, y, z)$, their median, that lies on the geodesics connecting any two of them. Cubings look locally like products of trees. Between two points $x$ and $y$ there may exist many shortest paths, let $[x, y]$ denote the set of vertices lying on all those paths. Then it is still true that there exists a vertex $m(x, y, z)$ with

$$[x, y] \cap [y, z] \cap [z, x] = \{m(x, y, z)\}.$$

The median operation $m(x, y, z)$ has been known and studied for a long time (see BIRKHOFF-KISS [1947], GRAU [1947]), mainly in the context of universal algebra. The best surveys I am aware of are by ISBELL [1980] and BANDELT-HEDLÍKOVÁ [1983]. The archetypal median algebra is a power set $\mathscr{P}X$ with the Boolean median operation

$$m(x, y, z) = (x \cap y) \cup (y \cap z) \cup (z \cup x) \qquad \text{for } x, y, z \in \mathscr{P}X.$$

A general median algebra is a subalgebra of a power set. Together with a notion of intervals one also has convex sets; half spaces are convex sets whose complements are convex, or prime ideals in the sense of universal algebra. Complementary half spaces may be viewed as the two orientations of a hyper plane, though we usually take the reverse point of view and regard half spaces as the main objects and a hyper plane as an unordered pair $\{H, H^*\}$ of a proper half space and its complement.

The key observation here is that the half spaces of median algebra naturally form a poc set, and a simple argument shows that the family of ultra filters of a poc set is closed under the median operation. Thus we have the combinatorial basis of a duality. This duality is an anti isomorphism if we restrict it to the categories of finite median algebras and poc sets; in general a median algebra need not be the dual of a poc set and vice versa. Dunwoody's construction, however, yields the anti dual to the poc set of edges of a tree. In order to completely understand this duality we need to topologize the dual objects. Taking this topology into account, one can then decide whether a poc set or a median algebra has an anti dual. Here we discover a major difference between the two categories: For a median algebra there exists at most one topology that turns it into a dual, but for a poc set there may be many different topologies, and there appears to be no canonical choice. This is why we finally use discrete representations to find the anti dual for a poc set, which carry a little more information than just the bare poc structure.

One of our results, Theorem 10.3, says that $M$ is an interval finite median algebra precisely if $M$ is the 0-skeleton of a cubing, with the median operation defined as above. (Examples of non interval finite median algebras are the vertex sets of Dunwoody's pro-trees.) This result encapsulates all topological arguments of Sageev's work. The theory of groups acting on interval finite median algebras is identical with the theory of groups acting on cubings, but in the former one can concentrate mainly on combinatorial arguments. That philosophy will be exemplified in §11.

To explain the essentiality condition for hyper planes in Sageev's Theorem we have to go a bit more into details. For any hyper plane $H$ one can find a $K$-almost invariant subset $A$. But $A$ only detects the presence of ends if it is proper, i.e., neither $K$-finite nor -cofinite. Sageev calls $H$ essential if the corresponding $A$ is proper, and a cubing is called essential if it contains an essential hyper plane. His characterization of essential hyper planes is just a rewriting of the algebraic condition, but it is not easy to understand geometrically. For example, he asks whether a cubing is essential if $G$ acts on it with unbounded orbits. He can show this for finite dimensional cubings by finding a hyper plane that is shifted by some element of the group. (For a more general condition see Proposition 11.3). Theorem 11.5 says that Sageev's conjecture is true for simple median algebras, which contain a single orbit of hyper planes. Similarly, it is true for finitely generated groups, see Theorem 11.6.

Our main tool to deal with bounded actions is Theorem 11.9, which says that a group acts on a discrete median algebra with bounded orbits precisely if it leaves a finite cube fixed, which then means that it fixes the center of this cube in the corresponding cubing. In general, unbounded orbits are not enough to find essential hyper planes, as the examples of parabolic group actions demonstrate. Our geometric characterization goes as follows: a hyper plane $H$ is essential iff the orbit $GH$ is infinite and $H$ divides some vertex orbit $Gx$ into two unbounded sets.

I hope that I can convince the reader that discrete median algebras are the ideal language to explain and work with Dunwoody's Construction. I also hope that median algebras will prove useful to deal with the other half of splitting theory, the careful choice procedures. Median algebras may help to understand the intermediate steps in the procedure and, since they are much more flexible than just trees, may provide new choice procedures.

I was surprised to learn that median algebras are very little known amongst group theorists and only recently attracted their interest. Thus it was at a very late stage of this work that two papers were brought to my attention which have overlaps with the present paper. The first is by Basarab, who has a number of papers on median algebras, which he calls generalized trees. BASARAB [1992] contains the same material as the first subsection of § 5, and although his language may be different, his methods are more or less the same as mine. GERASIMOV [1997] contains the same results as § 10 and a proof of Theorem 11.6, but his methods appear to be essentially different from mine. Part of § 11, in particular the equivalence $a \Leftrightarrow b$ of Theorem 11.5 has been published as joint work in NIBLO-ROLLER [1998]. In that paper we formulated our results carefully so as to expurgate any reference to median algebras and poc sets.

CONTENTS

### § 1. Axioms

We briefly introduce the axiomatic definitions of median algebras and poc sets, together with some standard examples that we shall meet repeatedly. The main example to keep in mind is that of a tree, whose vertices carry a natural median structure and whose oriented edges carry a poc structure.

### § 2. Geometry of Median Algebras

This is a survey of the elementary theory of median algebras that forms the basis of the remainder of this work. The most important concept is that of a half space; the most important result, Theorem 2.7, says that half spaces separate points.

### § 3. The Structure of Poc Sets

Poc sets are much simpler, so there is little to say about their general structure. The main concept here are ultra filters. The Extension Algorithm 3.3 is the counter part to Theorem 2.7. The Ramsey argument Corollary 3.7 will play a part in § 5 and § 3. The question to which extend poc sets are just products of trees is discussed in Proposition 3.12.

## §4. Duality

The half spaces of a median algebra form a poc set, the ultra filters of a poc set carry a median structure. This duality lies at the heart of Dunwoody's Construction. For the categories of finite median algebras and poc sets we get a a complete anti isomorphism. In general, the canonical double dual map is not surjective.

## §5. Duality for Poc Sets

We show that the category of poc sets is anti dual to the category of median algebras which carry a Stone topology such that the median operation is continuous.

We introduce the convex topology, defined on any median algebra, which is a Stone topology iff the double dual map is an isomorphism. This is the case iff the dual is a poc set of type $\omega$.

## §6. Duality for Median Algebras

In analogy to the previous section we show that the category of median algebras is anti dual to the category of poc sets which carry a compatible Stone topology. We give a similar characterization of poc sets which are canonically isomorphic to their double dual. Congruence relations on median algebras turn out to be related to the Stone topology on the dual. They will be used in §11 to construct simple median algebras.

## §7. Duality for Maps

Dual maps have more or less the properties that we expect for categorical reasons. As applications we study free median algebras and give a simple characterization of generating sets for median algebras.

## §8. Discrete Median Algebras

Discrete (i.e., interval finite) median algebras are the "geometries" we need to do geometric group theory. Here we collect a few useful characterizations and examples.

## §9. Discrete Poc Sets

Dunwoody uses a discrete (i.e., interval finite), nested poc set $P$ to construct a tree. The whole dual $P^\circ$ consists, however, of both vertices and ends of the tree. What we need is an anti dual, a discrete median algebra whose dual is the given poc set $P$. We push Dunwoody's Construction as far as possible, but, surprisingly, not every discrete poc set admits an anti dual.

In the group theoretic applications we get a little more information than just a discrete poc set. It comes as a family of subsets of the ambient group. We formalize this with the concept of discrete representations.

## § 10. Cubings

We use Sageev's topological methods to show that cubings and discrete median algebras are one and the same.

## § 11. Groups Acting on Discrete Median Algebras

We identify three classes of group actions on discrete median algebras. Elliptic actions are characterized by finite orbits; we prove a fixed point theorem for these actions. Parabolic actions have non principal fixed points in the double dual and only arise for large groups. Essential actions are characterized by Sageev's Theorem.

Notation

Words that are defined in the text are emphasized in sans serif font. For a set $X$ we denote the power set by $\mathscr{P}X$ and the family of finite subsets by $\mathscr{F}X$. We write $a \smallsetminus b := \{x \in a \mid x \notin b\}$ for the set theoretic difference and $a + b := (a \smallsetminus b) \cup (b \smallsetminus a)$ for the symmetric difference. Two subsets $a, b \subseteq X$ meet if $a \cap b \neq \varnothing$, and $a$ cuts $b$ if both $a \cap b \neq \varnothing$ and $b \smallsetminus a \neq \varnothing$. A family $\mathscr{S} \subseteq \mathscr{P}X$ is centered if for any finite subfamily $\mathscr{F} \subseteq \mathscr{S}$ the intersection $\bigcap \mathscr{F}$ is non empty. We shall freely use the abbreviation iff to mean if and only if, and xor to denote the exclusive or.

As is usual, the letters $\mathbb{N}$, $\mathbb{Z}$ and $\mathbb{R}$ denote the positive integers, the integers and the real numbers, respectively. The set $\mathbf{2} := \{0, 1\}$ plays a major part in this work and appears in many different rôles as a median algebra, a poc set, a Boolean algebra, a Stone space, a cyclic group and perhaps more — the context should make it clear which structure it has put on.

If $H$ is a subset of a median algebra $M$ (usually a half space), then we write $H^*$ to denote its complement. If $A$ is a subset of a poc set, then we write $A^*$ for the set $\{a^* \mid a \in A\}$. For an ultra filter $U$ the two interpretations of $U^*$ have the same meaning.

The symbol $\square$ denotes the absence of a proof, or the end of proof, example or remark like this. $\square$

Acknowledgements

This work has taken a long time in completing, and I have talked to a lot of people about it. I would like to thank Marcus Tressl, who prepared a special lecture series on Stone's duality theory just for my education. Claus Scheiderer and Manfred Knebusch invited me to report on my work in their Oberseminar. Theodor Bröcker helped me with the counter example in §9.

It was fun to work together with Martin Lustig and Jens Harlander who run the Frankfurt Seminar on Geometric Group Theory, and they suffered several talks of mine gladly. Şerban Basarab kindly made his unpublished papers available to me. Sarah Rees used her visit to Germany to work with me on median groups.

Most heartily I would like to thank Peter Kropholler, Martin Dunwoody and Graham Niblo, with whom I had the great honour and pleasure to do joint work. The Southampton Mathematics Department invited me several times, which gave me the opportunity to work with Graham and Martin, and meet Michah Sageev to discuss his theorem. Graham taught me the use of fixed point theorems that initiated the results of the last chapter. That chapter was vastly improved after my attempts to explain it to Peter. Last but not least, Martin's insights and ideas are a great stimulus for this entire work, and for a large part of geometric group theory.

# §1. Axioms

## 1.1. Median Axioms.

A median algebra is a set $M$ together with a function $m : M \times M \times M \to M$, the median, satisfying the following axioms:

(**Med 1**) $\quad \forall x, y, z \in M : \quad m(x, y, z) = m(x, z, y) = m(y, z, x),$

(**Med 2**) $\quad \forall x, y \in M : \quad m(x, x, y) = x,$

(**Med 3**) $\quad \forall x, y, z, u, v \in M : \quad m(m(x, y, z), u, v)) = m(x, m(y, u, v), m(z, u, v)).$

A morphism in the category of median algebras is a median preserving function $M \to M'$. We denote by $\mathrm{Hom}_{\mathrm{Med}}(M, M')$ the set of all median morphisms between $M$ and $M'$.

## 1.2. Examples.

(i) The two element set $\mathbf{2} := \{0, 1\}$ has a canonical median. If $x$, $y$ and $z$ are elements of $\mathbf{2}$, then at least two must agree, and by (**Med 1**) and (**Med 2**) this is the median of the three. To check that (**Med 3**) holds, observe that either $u = v$, in which case both sides of (**Med 3**) agree with $u$ and $v$, or $u \neq v$, and then both sides agree with $m(x, y, z)$.

(ii) If $(M_i \mid i \in I)$ is a family of median algebras, then the cartesian product $\prod_{i \in I} M_i$ with the median taken component wise, i.e. $m((x_i), (y_i), (z_i)) := (m(x_i, y_i, z_i))$, is again a median algebra.

(iii) In particular, if $I$ is any set, then the power set $\mathscr{P}I$ is a median algebra with the Boolean median operation

$$m(x, y, z) := (x \cap y) \cup (y \cap z) \cup (z \cap x) = (x \cup y) \cap (y \cup z) \cap (z \cup x).$$

(iv) A subset $M \subseteq \mathscr{P}I$ is a median algebra iff it is closed under $m$.

(v) Let $O$ be a totally ordered set, then for any three elements $x, y, z \in O$ there is one element contained in the interval bounded by the other two, which we define to be the median $m(x, y, z)$. The median algebra $\mathrm{Med}(O)$ induced from a total ordering in this way is called linear. It is easy to see that a linear median algebra is characterised by the property that the median of three elements actually equals one of these elements. □

## 1.3. Poc Set Axioms.

A poc set is a partially ordered set $P$ with a smallest element $0$ and an involution $^*$, satisfying the following axioms:

(**Poc 1**) $\quad \forall a, b \in P : \quad a \leq b \quad \Rightarrow \quad b^* \leq a^*,$

(**Poc 2**) $\quad \forall a \in P : \quad a \leq a^* \quad \Rightarrow \quad a = 0.$

A morphism in the category of poc sets is a monotone function $P \to P'$ which fixes 0 and respects $*$. By $\mathrm{Hom}_{\mathrm{Poc}}(P, P')$ we denote the set of all poc morphisms between $P$ and $P'$.

**1.4. Examples.**
(i) The set $\mathbf{2} = \{0, 1\}$ has a canonical poc structure, where $0 < 1$ and $0^* = 1$.

(ii) If $(P_i \mid i \in I)$ is a family of poc sets, then the cartesian product $\prod_{i \in I} P_i$ with the ordering $(a_i) \leq (b_i)$ iff $a_i \leq b_i$ for all $i \in I$, and $(a_i)^* := (a_i^*)$ is again a poc set.

(iii) In particular, if $I$ is any set, then the power set $\mathscr{P}I$ has a poc structure, where the order is set inclusion, and $*$ is the setwise complement in $I$.

(iv) A subset $P \subseteq \mathscr{P}I$ is a poc set iff $\varnothing \in P$ and $P$ is closed under $*$.

(v) Let $X$ be any set. We define $\mathrm{Poc}(X) := \{0, 0^*\} \sqcup X \sqcup X^*$, where $X^*$ is a copy of $X$ and $*: X \to X^*$ is a bijection. The only order relations we require are that 0 is the smallest and $0^*$ the largest element of $\mathrm{Poc}(X)$, all other elements are incomparable. Any poc set arising in this way is called orthogonal.

(vi) Let $O$ be a partially ordered set. We define $\mathrm{Poc}(O) := \{0, 0^*\} \sqcup O \sqcup O^*$ with a poc set structure, where 0 and $0^*$ are the smallest and the largest element of $\mathrm{Poc}(O)$, respectively, while $O^*$ is the same as $O$, but with the reverse order, and no element of $O$ is comparable with any element of $O^*$. A poc set arising in this way is called binary. If $O$ is totally ordered, then we call $\mathrm{Poc}(O)$ a linear poc set. One can show that a poc set $P$ is linear iff for any three elements in $P$ two are comparable with each other. $\square$

We call an element $a \in P$ proper if $a \notin \{0, 0^*\}$. The poc set axioms imply that for any two proper elements $a, b \in P$ at most one of the following relations holds.

(**Nest**)    $a \geq b, \quad a < b, \quad a \geq b^* \quad \text{or} \quad a < b^*.$

We say that two elements $a, b \in P$ are nested, in symbols: $a \parallel b$, if one of the relations (**Nest**) holds, otherwise they are transverse, in symbols: $a \pitchfork b$. A subset $S \subseteq P$ is called nested or transverse, if any two distinct elements of $S$ are nested or transverse, respectively.

**1.5. Convention.**
Observe that our axioms do not exclude the trivial median algebra $\varnothing$, nor the trivial poc set with a single element $0 = 0^*$. However, we will follow the usual practice and assume tacitly that all median algebras and all poc sets are non trivial.

**1.6. A Model Example.**
For a graph $\Gamma$ let $V\Gamma$ denote its vertex set and $\widetilde{E}\Gamma$ the oriented edge set of $\Gamma$. An oriented edge $e$ has an initial vertex $\iota(e)$, a terminal vertex $\tau(e)$, and an opposite edge $e^*$, satisfying $e^{**} = e$, $\iota(e^*) = \tau(e)$ and $\tau(e^*) = \iota(e)$. We allow that $\tau(e) = \iota(e)$, but forbid $e^* = e$. The unoriented edge set can be interpreted as $E\Gamma := \{\{e, e^*\} \mid e \in \widetilde{E}\Gamma\}$.

A tree is graph $T$ with the property, that cutting through any unoriented edge from $T$ yields a graph with exactly two connected components. By orienting the edge, we can distinguish the two components; say $T_e^\iota$ is the component containing the vertex $\iota(e)$, and $T_e^\tau$ the other component. There is a fundamental relation between the vertex set and the oriented edge set which we shall describe now. For $v \in VT$ and $e \in \widetilde{E}T$ we say $e$ **points to** $v$, and write $e \to v$, if $v$ is contained $T_e^\tau$. Otherwise, if $v$ lies in the component of $\iota(e)$, we say $v$ **points to** $e$ and write $v \to e$.

We define a partial ordering on $\widetilde{E}T$ by setting

$$e < f \quad :\Leftrightarrow \quad T_e^\iota \subset T_f^\iota,$$

i.e., if the geodesic from $\iota(e)$ to $\tau(f)$ contains $\tau(e)$ and $\iota(f)$ as interior points. Informally, $e$ and $f$ are comparable iff they point in the same direction, and the arrowheads indicate the ordering.

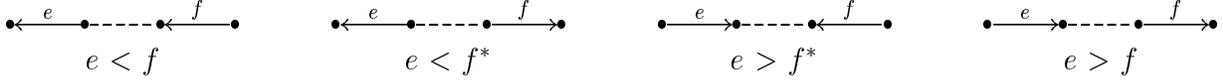

It is clear that $*$ reverses this order. Since there are no loops in $T$, no edge can be comparable to its opposite. If we adjoin a smallest element $0$ and a largest element $0^*$ to $\widetilde{E}T$ we obtain a poc set $\mathrm{Poc}(T)$. We point out two immediately obvious and characteristic properties of this poc set.

(**Nesting**)      Any pair of edges is nested.

(**Interval Finiteness**)      For any pair $e_1, e_2 \in \widetilde{E}T$ there exist only finitely many $f \in \widetilde{E}T$ with $e_1 < f < e_2$.

On the other hand we have a median algebra $\mathrm{Med}(T)$, whose underlying set is the vertex set $VT$ with the following median structure: The median of three vertices is the intersection of the geodesics connecting any two of them.

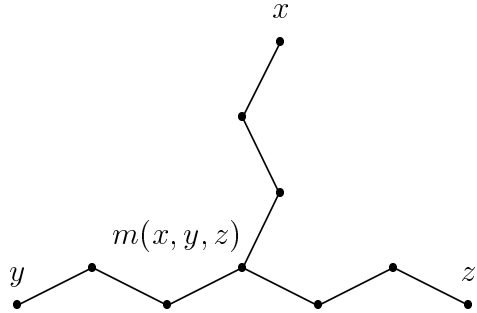

The first two axioms (**Med 1**) and (**Med 2**) obviously hold, but it seems tricky to verify (**Med 3**). To do this, we embed $VT$ into the power set of $\widetilde{E}T$. For a vertex $v \in V$ we define the set $\widehat{v} := \{e \in \widetilde{E}T \mid e \to v\}$. This defines an injection $VT \to \mathscr{P}\widetilde{E}T$, $v \mapsto \widehat{v}$, because if $v_1 \neq v_2$, then $\widehat{v}_1 + \widehat{v}_2$ contains all edges that lie on the geodesic from $v_1$ to $v_2$. A moment of thought convinces that

$e$ points to $m(x, y, z)$      $\Leftrightarrow$      $e$ points to at least two out of $x$, $y$, $z$.

This means that
$$\widehat{m(x,y,z)} = (\widehat{x} \cap \widehat{y}) \cup (\widehat{y} \cap \widehat{z}) \cup (\widehat{z} \cap \widehat{x}).$$

Thus the median on $VT$ agrees with the Boolean median operation on $\mathscr{P}\widetilde{E}T$, which satisfies (**Med 3**), as we know. □

**Notes.**

The definition of median algebras originated in the 1940's with BIRKHOFF-KISS [1947]. Since then they have been studied as examples of universal algebra, from a graph theoretic point of view and in the context of generalized convexity. The best survey of the history and classical theory of median algebras is BANDELT-HEDLÍKOVÁ [1983], for a detailed study of more general structures see ISBELL [1980].

Poc sets appeared briefly in ISBELL [1980] and WERNER [1981]. As our model example suggests, they are present in the background of the theory of group actions on trees, see e.g. STALLINGS [1971] and DICKS-DUNWOODY [1989].

# §2. The Geometry of Median Algebras

In this section we collect some basic facts about median algebras that will be needed throughout this paper. Most of the results are known, but for the benefit of the reader we provide full proofs.

### INTERVALS

Throughout this section, let $M$ denote a median algebra. For $x, y \in M$ we define the interval $[x, y] := \{z \in M \mid z = m(x, y, z)\}$. We will deduce some elementary properties of intervals in general median algebras. In the following arguments we will use the axioms (**Med 1**) and (**Med 2**) without further comment, but indicate the use of (**Med 3**).

(**Int 1**) $\quad [x, x] = \{x\} \quad$ and $\quad \{x, y\} \subseteq [x, y]$.

*Proof.* This follows immediately from (**Med 1**). $\square$

(**Int 2**) $\quad [x, y] = [y, x]$.

*Proof.* Again, this is obvious from (**Med 2**). $\square$

(**Int 3**) $\quad y \in [x, z] \quad \Rightarrow \quad [x, y] \subseteq [x, z]$.

*Proof.* Assume that $u \in [x, y]$, then

$$\begin{aligned}
m(u, x, z) &= m(m(u, x, y), x, z), & u &\in [x, y] \\
&= m(u, m(x, x, z), m(y, x, z)), & &(\textbf{Med 3}) \\
&= m(u, x, y), & y &\in [x, z] \\
&= u, & u &\in [x, y]
\end{aligned}$$

so $u \in [x, z]$. $\square$

(**Int 4**) $\quad m(x, y, z) \in [x, y]$.

*Proof.*
$$\begin{aligned}
m(m(x, y, z), x, y) &= m(z, m(x, x, y), m(x, y, y)), & (\textbf{Med 3}) \\
&= m(x, y, z).
\end{aligned}$$
$\square$

(**Int 5**) $\quad [x, y] \cap [x, z] = [x, m(x, y, z)]$.

*Proof.* We know from (**Int 3**) and (**Int 4**) that $[x, m(x, y, z)] \subseteq [x, y] \cap [x, z]$. To prove the converse assume that $u \in [x, y] \cap [x, z]$, then

$$\begin{aligned}
m(x, u, m(x, y, z)) &= m(x, m(x, u, y), m(x, u, z)), & (\textbf{Med 3}) \\
&= m(x, u, u), & \text{by assumption} \\
&= u,
\end{aligned}$$

so $u \in [x, m(x, y, z)]$. □

**(Int 6)**   $y \in [x, z]$ ⇔ $[x, y] \cap [y, z] = \{y\}$.

*Proof.* Since $[x, y] \cap [y, z] = [y, m(x, y, z)]$, this set is a singleton iff $y = m(x, y, z)$. □

**(Int 7)**   $[x, y] \cap [y, z] \cap [z, x] = \{m(x, y, z)\}$.

*Proof.*
$$[x, y] \cap [y, z] \cap [z, x] = [y, m(x, y, z)] \cap [z, m(x, y, z)] \quad \text{by (\textbf{Int 5})}$$
$$= \{m(x, y, z)\} \quad \text{by (\textbf{Int 4}) and (\textbf{Int 6})} \quad □$$

**(Int 8)**   $y \in [x, z]$ ⇒ $[x, w] \cap [z, w] \subseteq [y, w]$.

*Proof.* It is enough to show that if $y \in [x, z]$ then $m(x, z, w) \in [y, w]$.
$$m(y, w, m(x, z, w)) = m(w, m(x, y, z), m(x, z, w)), \qquad y \in [x, z]$$
$$= m(m(w, y, w), x, z), \qquad (\textbf{Med 3})$$
$$= m(x, z, w). \qquad □$$

**(Int 9)**   $y \in [x, z]$ and $x \in [w, z]$ ⇒ $x \in [w, y]$.

*Proof.*
$$m(w, x, y) = m(w, x, m(y, x, z)) \qquad y \in [x, z]$$
$$= m(y, m(w, x, x), m(w, x, z)) \qquad \text{by (\textbf{Med 3})}$$
$$= m(y, x, x) \qquad x \in [w, z]$$
$$= x. \qquad □$$

CONVEX SETS

We say that a subset $C \subseteq M$ is convex, if $[x, y] \subseteq C$ for all $x, y \in C$. In other words, $C = \{m(x, y, z) \mid x, y \in C, z \in M\}$, so the convex sets are precisely the ideals of the algebra $M$. The preimage of a convex set under a median morphism is also convex.
The characteristic properties of convex sets listed in the next proposition all follow immediately from the definitions.

**2.1. Proposition.**
  (o) *The empty set and $M$ are convex, singletons are convex.*
  (i) *For any family $\mathscr{C}$ of convex subsets of $M$ the intersection $\bigcap \mathscr{C}$ is convex.*
  (ii) *If $\mathscr{D}$ is an updirected family of convex subsets of $M$, i.e. for any $A, B \in \mathscr{D}$ there exists a $C \in \mathscr{D}$ with $A \cup B \subseteq C$, then the union $\bigcup \mathscr{D}$ is convex.* □

The next property is well known for intervals of the real line, and the proof carries over to our situation.

**2.2. Helly's Theorem.**
If $C_0, \ldots, C_n \subseteq M$ are convex and $C_i \cap C_j \neq \varnothing$ for all $i$ and $j$, then $\bigcap_{i=0}^{n} C_i \neq \varnothing$.

*Proof.* Consider three convex sets $C_0$, $C_1$ and $C_2$ and choose elements $x_i \in C_{i-1} \cap C_{i+1}$ (indices mod 3). Then both $x_{i-1}$ and $x_{i+1}$ lie in $C_i$, thus $m(x_0, x_1, x_2)$ lies in $C_i$ for $i = 0, 1, 2$. The general case follows by induction and applying the previous case to the sets $C_0$, $C_1$ and $C_2 \cap C_3 \cap \ldots \cap C_n$. □

The convex hull $\mathrm{conv}(X)$ of a subset $X \subseteq M$ is defined as the smallest convex set which contains $X$. By (i) above, $\mathrm{conv}(X)$ is the intersection of all convex sets containing $X$. For subsets $X, Y \subseteq M$ we define their join $[X, Y] := \bigcup \{[x, y] \mid x \in X,\ y \in Y\}$.

### 2.3. Proposition.
*If $X, Y \subseteq M$ are convex subsets, then $\mathrm{conv}(X \cup Y) = [X, Y]$.*

*Proof.* Obviously $[X, Y]$ is contained in the convex hull of $X \cup Y$ and contains both $X$ and $Y$, thus it is enough to show that $[X, Y]$ is convex. So if $x_1, x_2 \in X$, $y_1, y_2 \in Y$, further $z_1 \in [x_1, y_1]$ and $z_2 \in [x_2, y_2]$, we have to show that any $z_3 \in [z_1, z_2]$ lies in $[X, Y]$, i.e., is contained in an interval $[x_3, y_3]$. This is the content of the following lemma, where we show that one can choose $x_3$ as the median of $x_1$, $x_2$ and $z_3$, which lies in the convex set $X$, and similarly $y_3 = m(y_1, y_2, z_3)$.

### 2.4. Lemma.
*Let $x_i, y_i, z_i \in M$. Suppose that*

(i) $z_1 \in [x_1, y_1]$ and $z_2 \in [x_2, y_2]$,
(ii) $z_3 \in [z_1, z_2]$.

*Let*

(iii) $x_3 := m(x_1, x_2, z_3)$ and $y_3 := m(y_1, y_2, z_3)$.

*Then $z_3 \in [x_3, y_3]$.*

*Proof.*
$$
\begin{aligned}
[z_3, x_3] \cap [z_3, y_3] &= [z_3, x_1] \cap [z_3, x_2] \cap [z_3, y_1] \cap [z_3, y_2] && \text{by (\textbf{Int 5})} \\
&= [z_3, x_1] \cap [z_3, y_1] \cap [z_3, x_2] \cap [z_3, y_2] \\
&\subseteq [z_3, z_1] \cap [z_3, z_2] && \text{by (\textbf{Int 8}) and (\textbf{Int 3})} \\
&= \{z_3\} && \text{since } z_3 \in [z_1, z_2]
\end{aligned}
$$
so $z_3 \in [x_3, y_3]$, by (**Int 6**). □

### 2.5. Corollary.
*For any set $X \subseteq M$ we have $\mathrm{conv}(X) = \bigcup \{\mathrm{conv}(F) \mid F \subseteq X,\ F \text{ finite}\}$.*

*Proof.* Let $C := \bigcup \{\mathrm{conv}(F) \mid F \subseteq X,\ F \text{ finite}\}$, then obviously any convex set containing $X$ must also contain $C$, so it is enough to show that $C$ is convex. If $x_1, x_2$ lie in $C$, then there exist finite sets $F_1, F_2 \subseteq X$ such that $x_i \in \mathrm{conv}(F_i)$. The previous proposition now says that $[x_1, x_2] \subseteq \mathrm{conv}(F_1 \cup F_2) \subseteq C$. □

Suppose that $N$ is a median subalgebra of $M$ and let $X$ be a subset of $N$. In the next proposition we denote the convex hull of $X$ in $N$ by $\mathrm{conv}_N(X)$, and similarly we write $[x, y]_N$ for an interval in $N$.

**2.6. Proposition.**
If $N \leq M$ is a median subalgebra and $X \subseteq N$ then
$$\operatorname{conv}_N(X) = N \cap \operatorname{conv}_M(X).$$

*Proof.* If $X$ has at most one element, nothing is to show. Now $N$ is a median subalgebra, so for all $x, y \in N$ we have
$$\begin{aligned}[x, y]_N &= \{z \in N \mid z = m(x, y, z)\} \\ &= \{z \in M \mid z \in N \text{ and } z = m(x, y, z)\} \\ &= N \cap [x, y]_M,\end{aligned}$$
which proves the proposition in case that $X$ has two elements.

For finite $X$ we can prove the proposition inductively, using Proposition 2.3.
$$\begin{aligned}\operatorname{conv}_N(X \cup \{z\}) = [X, \{z\}]_N &= \bigcup_{x \in X} [x, z]_N \\ &= \bigcup_{x \in X} N \cap [x, z]_M = N \cap \bigcup_{x \in X} [x, z]_M \\ &= N \cap [X, \{z\}]_M = N \cap \operatorname{conv}_M(X \cup \{z\}).\end{aligned}$$

For an arbitrary subset $X \subseteq N$ we can use a similar idea together with Corollary 2.5.
$$\begin{aligned}\operatorname{conv}_N(X) &= \bigcup \{\operatorname{conv}_N(F) \mid F \in \mathscr{F}X\} \\ &= \bigcup \{N \cap \operatorname{conv}_M(F) \mid F \in \mathscr{F}X\} \\ &= N \cap \bigcup \{\operatorname{conv}_M(F) \mid F \in \mathscr{F}X\} \\ &= N \cap \operatorname{conv}_M(X).\end{aligned}$$

This completes our proof. □

Half Spaces

A subset $H \subseteq M$ is a **half space** if both $H$ and $M \smallsetminus H$ are convex. If $h\colon M \to \mathbf{2}$ is a median morphism, then both its support $h^{-1}(1)$ and its kernel $h^{-1}(0)$ are convex, thus they are both half spaces. This includes the degenerate case of a constant map, where the half spaces are empty or all of $M$. Conversely, for any half space $H \subseteq M$ the characteristic function $\chi_H\colon M \to \mathbf{2}$ is a median morphism.

The next theorem says that half spaces separate convex sets from disjoint points, in particular every convex set is an intersection of half spaces.

**2.7. Theorem.**
Let $A \subseteq M$ be a convex subset and $b \in M \smallsetminus A$. Then there exists a maximal convex set $H \subset M$ with $A \subseteq H$ and $b \notin H$, and every such $H$ is a half space.

*Proof.* Consider $\mathscr{H} := \{C \subseteq M \mid C \text{ convex}, A \subseteq C, b \notin C\}$. If $\mathscr{C} \subseteq \mathscr{H}$ is an ascending chain with respect to inclusion, then $\bigcup \mathscr{C}$ also lies in $\mathscr{H}$, thus by Zorn's Lemma $\mathscr{H}$ contains a maximal element.

Let $H$ be any maximal element of $\mathscr{H}$. To show that $H$ is a half space it remains to prove that $M \smallsetminus H$ is convex. Suppose that $x_1, x_2 \in M \smallsetminus H$ and there exists $x_3 \in [x_1, x_2]$ with $x_3 \in H$. By maximality, both $[H, \{x_1\}]$ and $[H, \{x_2\}]$ contain $b$, thus there exist $h_1, h_2 \in H$ with $b \in [h_1, x_1]$ and $b \in [h_2, x_2]$. Set $h_3 := m(h_1, h_2, b)$, then $h_3$ lies in both $[b, h_1]$ and $[b, h_2]$. Thus by (**Int 9**) and (**Int 8**) above we have $b \in [x_1, h_3] \cap [x_2, h_3] \subseteq [x_3, h_3]$. But both $h_3$ and $x_3$ lie in $H$, and $H$ is convex, thus $b \in H$, which is a contradiction. □

We can jazz up this proof to show that half spaces separate disjoint convex sets.

**2.8. Theorem.**
Let $A$ and $B$ be disjoint convex subsets of $M$. Then there exists a half space $H \subseteq M$ with $A \subseteq H$ and $B \cap H = \varnothing$.

*Proof.* Let again $H$ be a maximal convex set with $A \subseteq H$ and $B \cap H = \varnothing$. Suppose that $x_1, x_2 \in M \smallsetminus H$ but there exists $x_3 \in [x_1, x_2]$ with $x_3 \in H$. The maximality of $H$ implies that there exist $b_1, b_2 \in B$ and $h_1, h_2 \in H$ with $b_i \in [x_i, h_i]$. Using the lemma below we find elements $b_3 \in B$ and $h_3 \in H$, such that $b_3 \in [x_3, h_3]$. But the latter interval is contained in $H$, thus $B$ meets $H$, which is a contradiction. □

**2.9. Lemma.**
Let $x_i, b_i, h_i \in M$. Suppose that

   (i) $x_3 \in [x_1, x_2]$,
   (ii) $b_1 \in [x_1, h_1]$ and $b_2 \in [x_2, h_2]$.

Let

   (iii) $b_3 := m(b_1, b_2, x_3)$,
   (iv) $h_3 := m(h_1, h_2, b_3)$.

Then $b_3 \in [x_3, h_3]$.

*Proof.* Suppose that $b_3 \notin [x_3, h_3]$. Then by Theorem 2.7 there exists a half space $H \subseteq M$ such that $b_3 \in H$ and $x_3, h_3 \notin H$. We will show that this is not possible.

If $h_1$ or $h_2$ lies in $H$, then by (iv) we also have $h_3 \in H$. Thus neither $h_1$ nor $h_2$ lies in $H$. If $b_1$ or $b_2$ lies in $M \smallsetminus H$, then by (iii) we also have $b_3 \in M \smallsetminus H$. Thus both $b_1$ and $b_2$ lie in $H$. Now (ii) implies that both $x_1$ and $x_2$ lie in $H$. From (i) we infer that $x_3 \in H$, which is a contradiction. □

**Remark.**
Theorem 2.7 and Theorem 2.8 are of fundamental importance to the theory of median algebras. They imply that any median algebra arises as a subalgebra of a power set with the Boolean median operation. This is a much more useful characterization than any axiomatic definition. In the language of Universal Algebra this means that the median algebras form a variety generated by **2**.

A consequence of this observation is that any identity

(∗) $\quad\quad\quad \forall x_1, \ldots, x_n \in M : \quad \varphi_1(x_1, \ldots, x_n) = \varphi_2(x_1, \ldots, x_n),$

where $\varphi_1$ and $\varphi_2$ are expressions obtained by repeated taking of medians, is deducible from the median axioms if and only if it holds true for $M = \mathbf{2}$.

Thus an identity (e.g., Lemma 2.4 and Lemma 2.9) could in principle be checked with a straightforward computer program. Rather than referring the reader to a median version of MAPLE we use another device. When trying to prove a formula like (∗) that is supposed to be universally true, we assume that it is false for some $x_1, \ldots, x_n \in M$. By Theorem 2.7 there exists a half space $H \subset M$ which separates $\varphi_1(x_1, \ldots, x_n)$ from $\varphi_2(x_1, \ldots, x_n)$. We then start checking which of the $x_i$ lie inside or outside $H$, until we reach a contradiction (see e.g. the proof of Lemma 2.9). This method of proof conveys the same feeling of inevitability as the diagram chases in Homological Algebra.

Yet another consequence is a concrete construction of free median algebras, which we will see in §6. □

**Notes.**

Many systems of axioms for median algebras are known, see BANDELT-HEDLÍKOVÁ [1983]. By (**Int 4**) and (**Int 7**) one can recover the median operation from the interval structure. A particularly elegant system of axioms was found by SHOLANDER [1954].

$$\forall x \in M : \quad [x, x] = \{x\}.$$

$$\forall x, y, z \in M : \quad y \in [x, z] \quad \Rightarrow \quad [x, y] \subseteq [z, y].$$

$$\forall x, y, z \in M \; \exists w \in M : \quad [x, y] \cap [y, z] \cap [z, x] = \{w\} \quad \text{(where } w = m(x, y, z)\text{)}.$$

Theorem 2.8 appeared first in NIEMINEN [1978], Thm. 1. The proof given there depends on Thm. 2.2 of BALBES [1969], but it contains a small gap, so Nieminen only proves our Theorem 2.7 (see also Lemma 1.4 in Bandelt-Hedlíková [1983]). Since Theorem 2.8 is of fundamental importance for the theory of median algebras I have included a proof from first principles. Thm. 2.5 of VAN DE VEL [1984] contains our Theorem 2.8 as a special case; another proof can be found in §5.2 of BASARAB [1992].

SEPARATORS, RETRACTS AND GATES

Let $H \subseteq M$ be a half space. We write $H^* := M \smallsetminus H$ for the complement of $H$ in $M$. For subsets $X, Y \subseteq M$ we define their **separator** $\Delta(X, Y)$ as the set of all half spaces $H \subseteq M$ such that $X \subseteq H$ and $Y \subseteq H^*$. For elements $x, y \in M$ we will also write $\Delta(x, Y) := \Delta(\{x\}, Y)$ and $\Delta(x, y) := \Delta(\{x\}, \{y\})$.

It is clear that $\Delta(X, Y) = \Delta(\mathrm{conv}(X), \mathrm{conv}(Y))$, and Theorem 2.8 says that this set is non empty iff $\mathrm{conv}(X) \cap \mathrm{conv}(Y) = \varnothing$.

According to the final remark of the previous section it is sometimes useful to reformulate identities in median algebras in terms of separators. The next lemma is an example of this philosophy.

**2.10. Lemma.**
Let $x_1, x_2, x_3 \in M$. The following are equivalent.

(i) $x_3 \in [x_1, x_2]$.

(ii) $\forall y \in M : \quad \Delta(y, x_3) \subseteq \Delta(y, x_1) \cup \Delta(y, x_2)$.

(iii) $\Delta(x_1, x_3) \subseteq \Delta(x_1, x_2)$.

(iv) $\Delta(x_3, x_2) \subseteq \Delta(x_1, x_2)$.

(v) $\Delta(x_1, x_3) \cup \Delta(x_3, x_2) = \Delta(x_1, x_2)$.

*Proof.* Suppose $H \in \Delta(y, x_3) \smallsetminus (\Delta(y, x_1) \cup \Delta(y, x_2))$, then $H$ separates $x_3$ from $\{x_1, x_2\}$, thus (i) implies (ii). Both (iii) and (iv) are special cases of (ii) when $y = x_1$ or $y = x_2$. If either (iii) or (iv) holds, then no half space can separate $x_3$ from $\{x_1, x_2\}$, thus $x_3 \in \mathrm{conv}(\{x_1, x_2\}) = [x_1, x_2]$, i.e., (i) holds. Observe that $\Delta(x_1, x_2) \subseteq \Delta(x_1, x_3) \cup \Delta(x_3, x_2)$, regardless whether $x_3$ is contained in $[x_1, x_2]$ or not. Thus (v) holds iff (iii) and (iv) hold. □

Now suppose $X \subseteq M$ is nonempty. For any $y \in M$ we say that the element $x \in X$ is closest to $y$, and write $x = \mathrm{retr}_X(y)$, if $x \in [y, z]$ for all $z \in X$. Note that a closest point, if it exists, is necessarily unique: If $x_1, x_2 \in X$ are closest to $y$, then $x_1 \in [y, x_2]$ and $x_2 \in [y, x_1]$, thus
$$x_1 = m(x_1, x_2, y) = x_2.$$
Again, this property can be expressed in terms of separators.

**2.11. Lemma.**
For $x \in X \subseteq M$ and $y \in M$ the following are equivalent.

(i) $x \in X$ is closest to $y$.

(ii) $\Delta(x, y) = \Delta(X, y)$.

*Proof.* Obviously we have $\Delta(X, y) \subseteq \Delta(x, y)$. There exists a half space $H \in \Delta(x, y) \smallsetminus \Delta(X, y)$ iff there is an element $z \in X$ with $H \in \Delta(x, \{y, z\})$. This is the case iff $x \notin [y, z]$, i.e., iff $x \in X$ is not closest to $y$. □

We say that $X \subseteq M$ is a retract, if for every $y \in M$ there exists a point $x \in X$ closest to $y$. In this case $X$ must be convex: Suppose that $z_1, z_2 \in X$ and $z_3 \in [z_1, z_2]$, then there exists an $x_3 \in X$ closest to $z_3$; in particular
$$x_3 \in [z_1, z_3] \cap [z_2, z_3] = \{z_3\} \qquad \text{by } (\mathbf{Int\ 6}),$$
hence $z_3 = x_3 \in X$.

If $X$ is a retract then the map $M \to M$, $y \mapsto \mathrm{retr}_X(y)$ is called the **retraction** onto $X$. Our next proposition collects useful properties of retracts and retractions. They can all be proved by half space chasing arguments. We provide one example of such a proof and leave the rest as an amusement for the reader. Axiomatic proofs can be found in ISBELL [1980] and BANDELT-HEDLÍKOVÁ [1983].

**2.12. Proposition.**

(i) *A median morphism $r\colon M \to M$ is a retraction iff the image $r(M)$ is convex and $r^2 = r$.*

(ii) *Any map $r\colon M \to M$ is a retraction iff*
$$\forall x, y, z \in M: \quad r(m(x,y,z)) = m(x, r(y), r(z)).$$

(iii) *Let $X, Y \subseteq M$ be retracts. Then $X \cap Y \neq \varnothing$ iff $\mathrm{retr}_X$ and $\mathrm{retr}_Y$ commute, and then $\mathrm{retr}_{X \cap Y} = \mathrm{retr}_X \circ \mathrm{retr}_Y$.*

(iv) *If $X, Y \subseteq M$ are retracts, then the join $[X, Y]$ is again a retract with retraction*
$$\mathrm{retr}_{[X,Y]}(z) = m(\mathrm{retr}_X(z), \mathrm{retr}_Y(z), z).$$

*Proof of* (iv). Let $z \in M$ and choose $x \in X$ and $y \in Y$ closest to $z$. We claim that $m := m(x, y, z) \in [X, Y]$ is closest to $z$. If not, then there exists a $w \in [x', y']$ with $x' \in X$, $y' \in Y$, such that $m \notin [w, z]$. This means that there exists a half space $H \in \Delta(\{w, z\}, m)$. Now suppose that $x'$ lies in $H$, then by the choice of $x$ we have $x \in [x', z] \subseteq H$, therefore $m \in [x, z] \subseteq H$, contradicting the choice of $H$. Thus we have $x' \in H^*$, and by the same argument we have $y' \in H^*$. Hence $w \in [x', y'] \subseteq H^*$, which is our final contradiction. □

**2.13. Examples.**

(i) Singletons are retracts.

(ii) By (iv) above, any interval $[x, y]$ is a retract with $\mathrm{retr}_{[x,y]}(z) = m(x, y, z)$.

This gives a new interpretation of the median operation; in particular the axiom (**Med 3**) can be understood as a special case of (ii) above.

(iii) More generally, if $X := \{x_1, \ldots, x_n\}$ is a finite subset of $M$, then $\mathrm{conv}(X)$ is a retract. The retraction onto $\mathrm{conv}(X)$ has a particularly neat expression in the case when $M$ is a subalgebra of a power set.
$$\mathrm{retr}_{\mathrm{conv}(X)}(z) = \bigcap X \cup (z \cap \bigcup X).$$
The proof is a straightforward induction on $n$.

(iv) The simplest examples of convex sets that are not retracts are found amongst linear median algebras. Let $O$ be a totally ordered set. A subset of $O$ is convex in the median algebra $\mathrm{med}(O)$ iff it is convex as a subset of $O$, and it is a retract iff it is a closed interval in the order topology on $O$, i.e., if it is of the form
$$\{z \in O \mid z \leq y\} \quad \text{or} \quad \{z \in O \mid x \leq z \leq y\} \quad \text{or} \quad \{z \in O \mid x \leq z\}. \quad \square$$

The next lemma is a symmetrical version of Lemma 2.11. The proof is analogous.

**2.14. Lemma.**
Let $x \in X \subseteq M$ and $y \in Y \subseteq M$. The following are equivalent.

(i) $x \in X$ is closest to $y$ and $y \in Y$ is closest to $x$.

(ii) $\Delta(x, y) = \Delta(X, Y)$. □

If such a pair $(x, y) \in X \times Y$ exists, we call it a gate for $X$ and $Y$.

**2.15. Lemma.**
If $X$ and $Y$ are retracts of $M$, then they have a gate.

*Proof.* If $X$ and $Y$ meet, then choose any $x = y \in X \cap Y$. Otherwise choose $x_0 \in X$ arbitrarily, let $y := \mathrm{retr}_Y(x_0)$ and $x := \mathrm{retr}_X(y)$. It remains to show that $y = \mathrm{retr}_Y(x)$. Take any $z \in Y$. By definition of $y$ we have $y \in [z, x_0]$. By definition of $x$ we have $x \in [x_0, y]$. In view of (**Int 9**) it follows that $y \in [z, x]$. Thus $y \in Y$ is closest to $x$. □

DISCRETE MEDIAN ALGEBRAS AND MEDIAN GRAPHS

To a half space $H \subseteq M$ we associate a hyper plane $\overline{H} := \{H, H^*\}$, and we may think of $H$ and $H^*$ as the two orientations of $\overline{H}$. We also define the set of hyper planes that separate $x$ and $y$ as $\overline{\Delta}(x, y) := \{\overline{H} \mid H \in \Delta(x, y)\}$. The following simple lemma will be very useful.

**2.16. Lemma.**
Let $x, y, z \in M$. Then

(i) $\overline{\Delta}(x, y) + \overline{\Delta}(y, z) = \overline{\Delta}(x, z)$.   *(Triangle Equality)*

(ii) $y \in [x, z] \iff \overline{\Delta}(x, y) \sqcup \overline{\Delta}(y, z) = \overline{\Delta}(x, z)$.

*Proof.* (i) For all hyper planes $\overline{H}$ in $M$ we have

$\overline{H}$ separates $x$ and $z$ $\iff$ ($\overline{H}$ separates $x$ and $y$)   xor   ($\overline{H}$ separates $y$ and $z$)

(ii) In view of (i) we have
$$\overline{\Delta}(x, y) \sqcup \overline{\Delta}(y, z) = \overline{\Delta}(x, z)$$
$$\iff \overline{\Delta}(x, y) \cap \overline{\Delta}(y, z) = \varnothing$$
$$\iff \Delta(\{x, z\}, y) = \varnothing$$
$$\iff y \in \mathrm{conv}(x, z) = [x, z]. \qquad \square$$

We say that $M$ is discrete if $\Delta(x, y)$ is finite for all $x, y \in M$. In this case the function $d(x, y) := |\Delta(x, y)| = |\overline{\Delta}(x, y)|$ is an integer valued metric on $M$: It is clearly symmetric and Theorem 2.7 says that $d(x, y) = 0$ iff $x = y$; the triangle inequality for $d$ is a consequence of our triangle equality.

It is instructive to have a more geometric interpretation of median algebras. The **median graph** $\Gamma M$ has $M$ as vertex set, and two distinct vertices $x$ and $y$ are adjacent if $\Delta(x,y)$ is a singleton, or, equivalently, $[x,y] = \{x,y\}$. A **path** of length $n$ in $\Gamma M$ from $x$ to $y$ is a sequence $x_0 = x$, $x_1$, ..., $x_n = y \in M$, where $x_{i-1}$ is adjacent to $x_i$ for all $i = 1, \ldots, n$. The distance between $x$ and $y$ in the **path metric** is the length of the shortest path connecting them.

### 2.17. Proposition.
*For a discrete median algebra the path metric on $\Gamma M$ agrees with $d$. Moreover, $M$ is discrete iff $\Gamma M$ is connected.*

*Proof.* Let $\overline{H} \in \overline{\Delta}(x,y)$, then any path $x_0 = x$, $x_1$, ..., $x_n = y$ must cross $\overline{H}$ in the sense that $\overline{\Delta}(x_{i-1}, x_i) = \{\overline{H}\}$ for some $i$. Therefore $n \geq d(x,y)$.
Conversely, we can prove that there exists a path of length $d(x,y)$ between $x$ and $y$. For $d(x,y) = 0$ or $1$ nothing is to prove. If $d(x,y) \geq 2$ then there exist two distinct half spaces $H_1, H_2 \in \Delta(x,y)$. At least one of $H_1 \smallsetminus H_2$ and $H_2 \smallsetminus H_1$ must be non empty, say $z \in H_1 \smallsetminus H_2$. Replacing $z$ by $m(x,y,z)$, we may assume that $z \in (H_1 \smallsetminus H_2) \cap [x,y]$. Then $d(x,y) = d(x,z) + d(z,y)$ and $d(x,z), d(z,y) \geq 1$, hence the proposition follows by induction. $\square$

We have just proved that the interval $[x,y]$ contains at least $d(x,y)+1$ distinct points. On the other hand, any point $z \in [x,y]$ is uniquely determined by $\Delta(\{x,z\}, y)$, hence we have
$$d(x,y) + 1 \leq |[x,y]| \leq 2^{d(x,y)}.$$
It follows that $[x,y]$ is finite iff $\Delta(x,y)$ is finite.
Observe that we can recover the median operation of $M$ from $\Gamma M$: The interval between $x$ and $y$ is given by $[x,y] = \{z \in M \mid d(x,y) = d(x,z) + d(z,y)\}$, and then $m(x,y,z)$ is determined by (**Int 7**). It follows that all graph automorphisms of $\Gamma M$ are also median automorphisms of $M$.

### 2.18. Examples.
Here are the median graphs with up to five vertices.

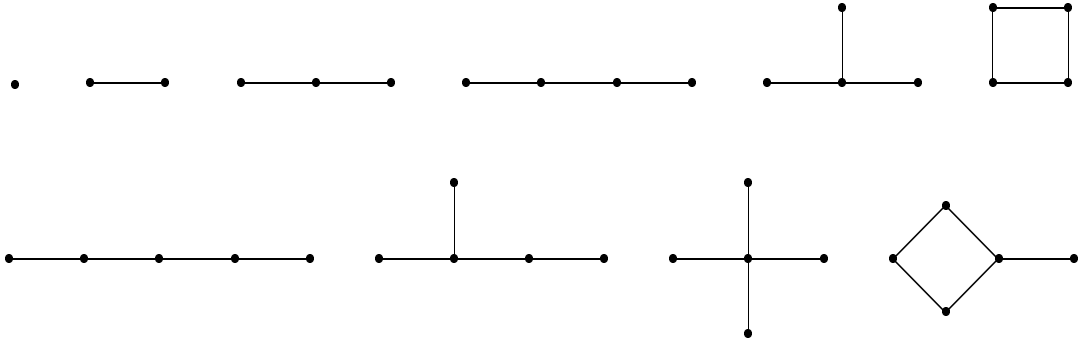

It is easy to see that median graphs must be bipartite. Choose any $x \in M$ and set $M_0 := \{y \in M \mid d(x,y) \text{ even}\}$ and $M_1 := M \smallsetminus M_0$. The triangle equality contradicts the assumption that any two elements of $M_i$ are adjacent. However, not all bipartite graphs are median, as the following examples show.

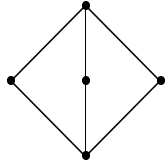 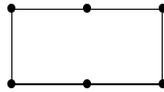 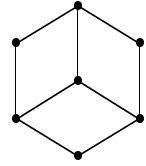

**Notes.**
There are many characterizations of median graphs in the literature, we just point the reader to BANDELT [1984] and BERRACHEDI [1994]. In §10 we prove a geometric characterization of median graphs, based on SAGEEV [1995], which is particularly handy to check finite graphs like the above examples.

# §3. The Structure of Poc Sets

The notion of a poc set may sound artificial, but it describes a situation that is typical for certain problems in low dimensional topology. Vice versa, we can use topological pictures and intuition to gain insight into poc sets. A **poc picture** consists of a collection $A$ of arcs in the unit disk $D$, with both end points on $\partial D$, such that any two arcs are either disjoint or intersect transversely in one point in the interior of $D$. The poc set $P$ obtained from this picture consists of a collection of subsets of $D$: We take $\varnothing$, $D$ as the improper elements, and for every arc $\alpha \in A$ we take the two components of $D \smallsetminus \alpha$, which are interchanged by $*$. The order is given by inclusion. Two proper elements of $P$ are transverse iff the corresponding arcs intersect; when the arcs don't intersect the corresponding elements of $P$ are nested (hence the symbols $\pitchfork$ and $\|$).

## 3.1. Examples.

Here are the poc pictures with up to three arcs, depicting the poc sets with up to eight elements.

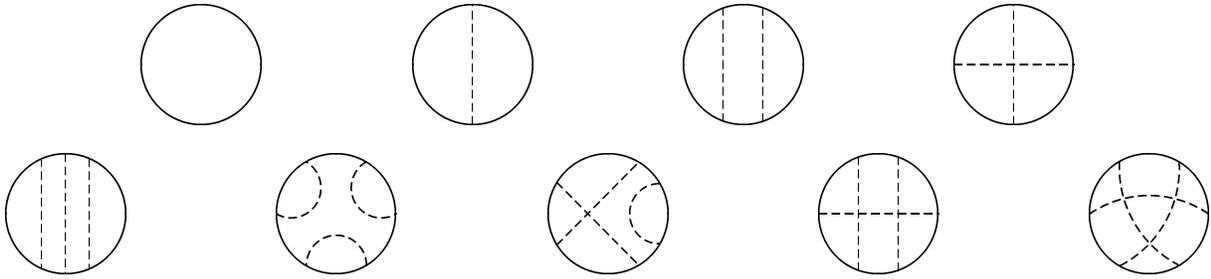

For example, the first picture in the first row represents the canonical poc set **2**. The first picture in the second row represents $\mathrm{Poc}(\{1,2,3\})$, where $\{1,2,3\}$ carries the usual ordering, whereas the last picture of the second row represents the orthogonal poc set $\mathrm{Poc}(\{\alpha,\beta,\gamma\})$. □

Now choose a point $x \in D$ that does not lie on any arc of $A$. Then for each $\alpha \in A$ we can distinguish the two components of $D \smallsetminus \alpha$, say $\alpha^+$ is the component which contains $x$, and $\alpha^-$ is the other component. In other words, we have chosen two orientations on $A$. The $+$-orientation has the characteristic property that for $\alpha, \beta \in A$ the components $\alpha^+$ and $\beta^+$ meet, because they both contain $x$. This can be expressed in terms of the ordering by saying that for $\alpha, \beta \in A$ we can never have $\alpha^+ \subseteq \beta^-$.

In the language of poc sets we take a slightly different point of view. If we specify an orientation by putting $+$ and $-$ signs on arcs, it is awkward to describe an orientation coming from a different point. Instead, we work with an operation $*$ which reverses the orientations. We call the positive orientations ultra filters, in formal analogy with the notion of an ultra filter in Boolean algebra.

Let $P$ be a partially ordered set. For any element $a \in P$ we define $\uparrow a := \{b \in A \mid b \geq a\}$, and for a subset $A \subseteq P$ we set $\uparrow A := \bigcup \{\uparrow a \mid a \in A\}$. We say that $A$ is an upper set if $A = \uparrow A$. Analogously, set $\downarrow a := \{b \in A \mid b \leq a\}$ and $\downarrow A := \bigcup \{\downarrow a \mid a \in A\}$. We say that $A$ is a lower set if $A = \downarrow A$.

Now let $P$ be a poc set. For a subset $A \subseteq P$ we define $A^* := \{a^* \mid a \in A\}$. The set $A \subset P$ is called orientable if $A \cap A^* = \varnothing$, and $A$ is an orientation, if $P$ is the disjoint union of $A$ and $A^*$.

A subset $U \subset P$ is called an ultra filter, if it is an orientation and an upper set. In other words, it satisfies the following two properties.

(**UF 1**) $\qquad \forall a \in P: \quad$ either $a \in U$ or $a^* \in U$, but not both.

(**UF 2**) $\qquad \forall a, b \in P: \quad a \not\leq b^*$.

We shall also use two weaker concepts. A non empty subset $F \subset P$ is called a filter, if it is an orientable upper set, i.e., it satisfies the following two properties.

(**Fil 1**) $\quad \forall a \in F \quad \forall b \in P: \quad b > a \quad \Rightarrow \quad b \in F$.

(**Fil 2**) $\qquad \forall a \in F: \quad a^* \notin F$.

Any subset $B$ of a filter $F$ must certainly satisfy (**UF 2**), and such a set is called a filter base.

Dually, a non empty subset $I \subset P$ is called an ideal, if $I^*$ is a filter, and $B \subset P$ is an ideal base if $B^*$ is a filter base. Observe that a subset $B \subset P$ may be both a filter base and an ideal base, e.g. if it is transverse or totally ordered.

### 3.2. Remark.

Though our notation emphasizes the analogy with filters and ultra filters in the Boolean sense, it is important to notice a fundamental difference. Any finite subset of a Boolean filter must have nonempty intersection. In the language of poc sets we can replace the condition $a \cap b \neq \varnothing$ by $a \not\leq b^*$, but the condition $a \cap b \cap c \neq \varnothing$ can not be expressed and is therefore not required for elements of a poc filter.

This accounts for the fact that our two dimensional poc pictures can sometimes be deceptive. For example, consider the picture below with the orientation specified by the center of the disc, which gives rise to the ultra filter $U_1 := \{\alpha^+, \beta^+, \gamma^+, 0^*\}$.

Observe that the components $\alpha^-$, $\beta^-$ and $\gamma^-$ also meet pairwise, thus we have another ultra filter $U_2 := \{\alpha^-, \beta^-, \gamma^-, 0^*\}$ which, however, does not correspond to any point in the left hand picture. But if we apply something like a Reidemeister move of type III to this picture we can "see" the ultra filter $U_2$ in the right hand picture.

This somewhat unsatisfactory explanation indicates that our pictures don't always tell us the whole truth about poc sets. We will come back to this problem in the next section. □

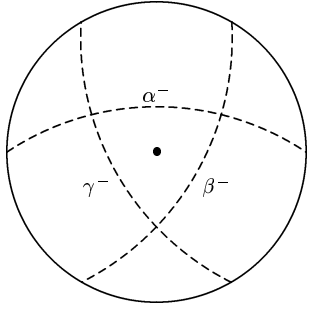 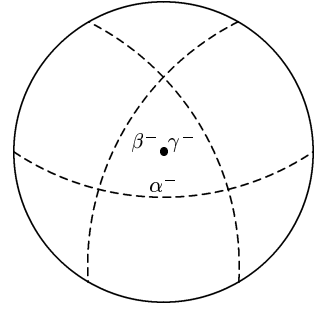

$$U_1 := \{\alpha^+, \beta^+, \gamma^+, 0^*\} \qquad U_2 := \{\alpha^-, \beta^-, \gamma^-, 0^*\}$$

Our main tool for constructing filters is the following algorithm.

### 3.3. Extension Algorithm.
*Consider a filter base $B \subset P$ and an element $a \in P$ such that $a, a^* \notin B$. Then there exists a filter base $B'$ that contains $B$ and either $a$ or $a^*$.*

CASE 1: *If there exists $b \in B$ with $b < a$, then put $B' := B \cup \{a\}$,*

CASE 2: *If there exists $b \in B$ with $b < a^*$, then put $B' := B \cup \{a^*\}$,*

CASE 3: *If neither CASE 1 nor CASE 2 holds, then both $B \cup \{a\}$ and $B \cup \{a^*\}$ are filter bases.*

*Proof.* The only way in which this algorithm could fail is that CASE 1 and CASE 2 happen simultaneously. This means that there exist $b_1, b_2 \in B$ with $b_1 < a$ and $b_2 < a^*$. But then $b_1 < a < b_2^*$, which is impossible if $B$ is a filter base. □

### 3.4. Applications.
(i) Given a filter base $B \subset P$, the set $\uparrow B$ is an extension of $B$ where CASE 1 applies throughout, thus it is a filter. As it is obviously the smallest filter containing $B$, we call it the filter generated by $B$.

(ii) A maximal filter base must also satisfy (**UF 1**), thus it is an ultra filter.

(iii) If $f \colon P \to P'$ is a poc morphism, then $f^{-1}(0^*)$ is a filter. If $f \colon P \to \mathbf{2}$, then its support $f^{-1}(1)$ is also a transversal to $^*$, hence an ultra filter. Conversely, any ultra filter in $P$ is the support of a poc morphism $P \to \mathbf{2}$.

(iv) By Zorn's Lemma every filter base is contained in an ultra filter.

(v) This last fact may also be expressed by saying that for any sub poc set $Q$ of $P$ the restriction map $\mathrm{Hom}_{\mathrm{Poc}}(P, \mathbf{2}) \to \mathrm{Hom}_{\mathrm{Poc}}(Q, \mathbf{2})$ is surjective.

(vi) More generally, if $F$ is a filter, $I$ an ideal and $F \cap I = \varnothing$, then there exists an ultra filter $U$ with $F \subset U$ and $I \cap U = \varnothing$. Simply start with the filter base $F \cup I^*$ and extend to an ultra filter as in (iv).

(vii) A filter is the intersection of all ultra filters containing it. □

### 3.5. Remark.
Let $P = \text{Poc}(O)$ be a binary poc set. Then we can write $P$ as a disjoint union $P = \{0, 0^*\} \sqcup O \sqcup O^*$, where $O$ is both a filter base and an ideal base. Conversely, any poc set which admits such a decomposition must be binary, i.e., no element of $O$ is comparable with any element of $O^*$.

In this case we can give a different description for the ultra filters in $P$, which, in the case of linear poc sets, is related to the notion of Dedekind cuts. Any ultra filter $U \subset P$ gives rise to a partition $O = U^+ \sqcup U^-$, where $U^+ := O \cap U$ is an upper set and $U^- := O \cap U^*$ is a lower set. Conversely, for any partition of $O$ into an upper set $U^+$ and a lower set $U^-$ we have an ultra filter $U := U^+ \cup U^{-*} \cup \{0^*\}$ on $\text{Poc}(O)$.

Here is a simple test for binary poc sets.

> $P$ is binary iff there does not exists a sequence $a_1, \ldots, a_n, a_{n+1} = a_1 \in P$ of proper elements, such that $a_i$ is comparable to $a_{i+1}^*$ for $i = 1, \ldots, n$, and $n$ is odd.

For example, the poc set pictured by ⊗ is not binary. □

## Transversality Graph and Decompositions

In this subsection we gather some general facts about the structure of poc sets. The involution $*$ in a poc set $P$ may be viewed as an equivalence relation whose classes have two elements. The quotient map will be denoted by $a \mapsto \overline{a} := \{a, a^*\}$. The quotient $\overline{P}$ is no longer a poc set, but the notion of nested and transverse pairs still makes sense in $\overline{P}$. We define the **transversality graph** $TP$ with vertex set $\{\overline{a} \in \overline{P} \mid a \text{ proper}\}$, where two vertices $\overline{a}$ and $\overline{b}$ are connected by an edge iff $a$ and $b$ are tranverse in $P$.

Now we recall the first non trivial result about general infinite graphs.

### 3.6. Ramsey's Theorem.
Let $X$ be a graph with infinite vertex set $VX$. Then there exists an infinite set $S \subseteq VX$ such that

(i) either any two vertices in $S$ are joined by an edge in $X$,

(ii) or no two vertices in $S$ are joined by an edge in $X$. □

Applied to the transversality graph of a poc set, we get an immediate corollary.

### 3.7. Corollary.
Every infinite poc set has an infinite subset that is either transverse or nested. □

This corollary will be used in conjunction with the following lemma, which is proved by inspecting the relevant definitions.

### 3.8. Lemma.
Suppose a subset $B$ of a poc set enjoys the following properties.

(i) $B$ is nested.
(ii) $B$ is a filter base.
(iii) $B$ is an ideal base.

Then $B$ is totally ordered. □

It is of course impossible, in general, to reconstruct a poc set from its transversality graph. For example, if the poc set is nested then the transversality graph has no edges. The next proposition shows, however, that one can always find some poc set that has a given graph as its transversality graph.

**3.9. Proposition.**
Let $\Gamma$ be a simple graph, i.e., $\Gamma$ has no loops and no multiple edges. Then there exists a poc set $P$, whose transversality graph is isomorphic to $\Gamma$.

*Proof.* Let $V$ and $E$ denote the vertex set and the unoriented edge set of $\Gamma$, respectively. Choose some object $z$, and set $X := \{z\} \sqcup V \sqcup E$. For any $v \in V$ define
$$A_v := \{v\} \cup \{e \in E \mid e \text{ incident with } v\}.$$
Let $P := \{\varnothing, X\} \cup \{A_v, X \smallsetminus A_v \mid v \in V\}$, then this is a poc set with inclusion as partial order and $A^* := X \smallsetminus A$. Since $v \in A_v$ and $z \in A_v^*$, the $A_v$'s are all proper. Consider two distinct vertices $v, w \in V$, then $A_v \not\subseteq A_w$, as $v \in A_v \smallsetminus A_w$ and $A_v^* \not\subseteq A_w$, as $z \in A_v^* \smallsetminus A_w$. The inclusion $A_v \subseteq A_w^*$ holds iff
$$\forall e \in E: \quad e \text{ incident with } v \quad \Rightarrow \quad e \text{ not incident with } w,$$
i.e., iff there is no edge between $v$ and $w$. Thus $A_v$ and $A_w$ are transverse iff $v$ and $w$ are connected by an edge. □

**3.10. Example.**
Let $\Gamma$ be a graph with no edges, and let $T$ be the cone on $\Gamma$, i.e., the tree with vertex set $V \sqcup \{z\}$, where every vertex in $V$ is connected to $z$. The poc set constructed in the previous proposition is isomorphic to $\text{Poc}(T)$, the poc set on the oriented edges of $T$ (see the Model Example 1.6). □

The length of a poc set $P$ is the maximal length of a chain of proper elements of $P$, and the dimension is the maximal cardinality of a transverse subset of proper elements. We say that $P$ has length $\omega$, if all chains in $P$ are finite, and dimension $\omega$, if all transverse subsets are finite. We say that $P$ is of type $\omega$, if it has both length $\omega$ and dimension $\omega$. A nested poc set, like the oriented edge set of a tree, always has dimension one, and the length corresponds to the length of a maximal geodesic. The dimension of $P$ is the maximal cardinality of a complete subgraph in $TP$.

The sum $P_1 \oplus P_2$ of two poc sets can be described in the following way: The set of proper elements of $P_1 \oplus P_2$ is the disjoint sum of the set of proper elements of $P_1$ and $P_2$, with the original orderings on $P_1$ and $P_2$, but any two elements $a_1 \in P_1$ and $a_2 \in P_2$ are transverse. A poc set $P$ is called prime, if $P = P_1 \oplus P_2$ implies that $P_1 = \mathbf{2}$ or $P_2 = \mathbf{2}$.

The summands of a poc set can be read off its transversality graph. The transversality graph $T(P_1 \oplus P_2)$ is just the graph theoretical join of $TP_1$ and $TP_2$. To say this in a different way consider the nesting graph $NP$, again with vertex set $\overline{P}$, where two distinct vertices $\overline{a}$ and $\overline{b}$ are joined by an edge iff $a$ and $b$ are nested elements of $P$; this is just the complementary graph to $TP$. Now $N(P_1 \oplus P_2)$ is the disjoint union of $NP_1$ and $NP_2$, hence the prime summands of $P$ correspond to connected components of $NP$.

**3.11. Proposition.**
$P$ is prime iff the nesting graph $NP$ is connected. $\hspace{1em}\square$

In this context we should mention an observation of Bell and VAN DE VEL [1986] Thm. 1.7, see also BANDELT-VAN DE VEL [1987], Cor. 3.6. Let us define the tree dimension of a poc set $P$ as the minimal number $n$ such that $P = P_1 \oplus P_2 \oplus \cdots \oplus P_n$, where the $P_i$ are nested poc sets (where $n = \infty$ is understood if no finite such decomposition exists). Obviously the dimension of a poc set is a lower bound of the tree dimension. Again the tree dimension is coded by the transversality graph. For if we have a decomposition of $P$ as above, then we have a partition of the vertex set of $TP$ into $n$ subsets such that no edge connects two vertices in the same subset, in other words a vertex colouring of $P$ with $n$ colours.

**3.12. Proposition.**
The tree dimension of $P$ is the chromatic number of the transversality graph $TP$. $\hspace{1em}\square$

**3.13. Example.**
A poc set of dimension 1 also has tree dimension 1. According to a famous theorem of Erdős there exist triangle free graphs of arbitrarily high chromatic number, thus by Proposition 3.9 and Proposition 3.12 there are poc sets of dimension 2 with arbitrarily high tree dimension. $\hspace{1em}\square$

## §4. Duality

Recall the poc pictures that we introduced in the previous section. The identification of regions of a picture with points suggests a duality similar to that of planar graphs. We further connect any two points with an edge, if they are separated by a single arc. In simple cases it turns out that we obtain a median graph, whose half spaces correspond to the elements of the poc set that we started out with.

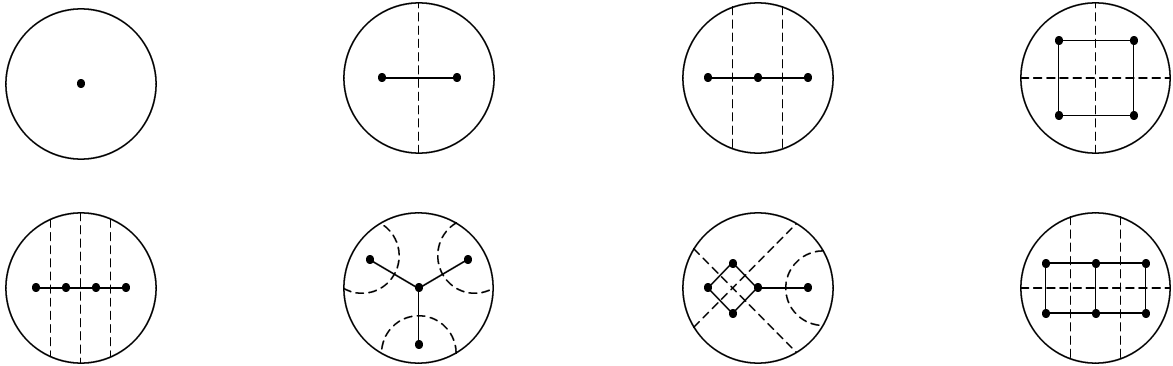

However, for the orthogonal poc set $\text{Poc}(\{\alpha, \beta, \gamma\})$, we found two different pictures representing the same poc set, which give rise to two different, but incomplete duals.

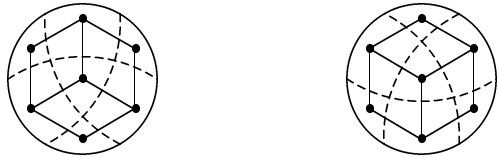

A more adequate representation of this poc set is given by a three-dimensional picture.

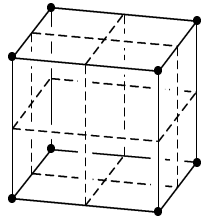

The main observation here is that the set of ultra filters on $P$ has a natural structure as a median algebra $M$, and the original poc set $P$ arises as the poc set of half spaces of $M$. Thus $M$ takes the place of the picture for $P$. For example, any finite set $\{a_1, \ldots, a_n\}$ of pairwise transverse elements of $P$ will manifest itself as an $n$-cube in $M$.

This duality works perfectly for finite poc sets and, more generally, for poc sets of type $\omega$. For infinite poc sets we can't expect to find a median algebra $M$ such that $P$ is isomorphic to the poc set of *all* half spaces in $M$. For a more detailed study we will have to introduce topologies on both poc sets and median algebras.

Before we formalize the duality between poc sets and median algebras, we recall the basic facts about its prototype, the duality theory of Stone spaces and Boolean algebras, due to STONE [1936].

A topological space $X$ is called a Stone space, if its topology has a subbasis $\mathscr{S}$ satisfying the following axioms.

(**Stone 1**) $\forall S \in \mathscr{S}: \quad X \smallsetminus S \in \mathscr{S},$ in particular $\mathscr{S}$ consists of clopen sets.

(**Stone 2**) $\mathscr{S}$ separates points, i.e., for distinct points $x, y \in X$ there exists an $S \in \mathscr{S}$ with $x \in S$ and $y \notin S$.

(**Stone 3**) $\mathscr{S}$ satisfies the finite intersection property (**FIP**): Every centered subfamily of $\mathscr{S}$ has nonempty intersection.

Equivalently, $X$ is a Stone space, iff it is Hausdorff, totally disconnected and quasi compact. Together with continuous maps the Stone spaces form a category.

On the other hand we have the category of Boolean algebras and Boolean homomorphisms, or, equivalently, the category of Boolean rings (i.e., commutative rings with 1, in which every element is idempotent) with ring homomorphisms.

The set $\mathbf{2} = \{0, 1\}$ is both a Stone space, when it is equipped with the discrete topology, and a Boolean ring, when it is viewed as $\mathbb{Z}_2$.

The dual of a Stone space $X$ is the ring $X^\circ$ of all continuous maps $X \to \mathbf{2}$, where the ring operations are defined point wise. The dual of a Boolean algebra $A$ is the space $A^\circ$ of all ring homomorphisms $A \to \mathbf{2}$, with the topology of pointwise convergence, i.e., the subspace topology of the Tychonov topology on $\mathscr{P}A$.

### 4.1. Stone Duality Theorem.

*The canonical evaluation map $ev(x)(f) := f(x)$ induces isomorphisms $ev\colon X \to X^{\circ\circ}$ and $ev\colon A \to A^{\circ\circ}$ of Stone spaces and Boolean algebras, respectively.*

*In other words: The categories of Stone spaces and Boolean algebras are anti equivalent.*
$\square$

Both poc sets and median algebras may be viewed as generalisations of Boolean algebras, and there is an analogue of Stone's Duality for either. In the remainder of this section we will consider a purely combinatorial version of this duality that does not give an anti equivalence in general: The double dual maps are not always surjective. In later sections we will study a canonical Stone topology on duals, which, amongst other things, allows us to determine exactly when the evaluation map is an isomorphism.

For a median algebra $M$ the dual $M^\circ := \mathrm{Hom}_{\mathrm{Med}}(M, \mathbf{2})$, as a subset of $\mathscr{P}M$, has a natural poc structure. More explicitly, for $h \in M^\circ$ we define $h^*(x) := h(x)^*$, and $h_1 \leq h_2$ iff $h_1(x) \leq h_2(x)$ for all $x \in M$. The constant map with value $0$ is the zero of $M^\circ$.

Similarly, for a poc set $P$ we have the **dual** $P^\circ := \text{Hom}_{\text{Poc}}(P, \mathbf{2})$. We will show that this is closed under the pointwise median operation, where for $f_1, f_2, f_3 \in P^\circ$ and $a \in P$ we define
$$m(f_1, f_2, f_3)(a) := m(f_1(a), f_2(a), f_3(a)).$$
Here $m(f_1, f_2, f_3)$ obviously preserves $0$ and $^*$, and for any pair $a, b \in P$ the values of $m(f_1, f_2, f_3)$ agree with those of one of the three functions, so there exists an $i \in \{1, 2, 3\}$ with $m(f_1, f_2, f_3)(a) = f_i(a)$ and $m(f_1, f_2, f_3)(b) = f_i(b)$, in particular, $m(f_1, f_2, f_3)$ is a monotone function and hence a poc morphism.

### 4.2. Convention.

We can view $M^\circ$ as the poc set of half spaces in $M$, ordered by inclusion, and $H^* := M \smallsetminus H$. Similarly we view $P^\circ$ as the median algebra of ultra filters in $P$ with the Boolean median operation.

More generally, we shall frequently identify functions $f\colon X \to \mathbf{2}$ with subsets of $X$. As usual, a function $f$ corresponds to its support $F := f^{-1}(1)$, and a subset $F \subseteq X$ to its characteristic function $f := \chi_F$. In this section we will switch between both notations, whereas in later sections we will prefer subsets over functions. □

The canonical evaluation function gives rise to **double dual** maps.
$$ev\colon M \to M^{\circ\circ} := \text{Hom}_{\text{Poc}}(\text{Hom}_{\text{Med}}(M, \mathbf{2}), \mathbf{2}), \quad ev(x)(h) := h(x),$$
$$ev\colon P \to P^{\circ\circ} := \text{Hom}_{\text{Med}}(\text{Hom}_{\text{Poc}}(P, \mathbf{2}), \mathbf{2}), \quad ev(a)(f) := f(a).$$
In the subset notation, $M^{\circ\circ}$ consists of all ultra filters of half spaces in $M^\circ$, and $ev(x)$ is the **principal** ultra filter defined by $ev(x) := \{H \in M^\circ \mid x \in H\}$. Similarly, $P^{\circ\circ}$ consists of all half spaces of ultra filters in $P^\circ$, and $ev(a)$ is a **principal** half space defined by $ev(a) := \{F \in P^\circ \mid a \in F\}$.

### 4.3. Examples.

(0) Although we vowed to disregard the trivial median algebra and poc set, it is amusing to see how they fit into our duality theory. For the trivial median algebra $M_{\text{triv}} = \varnothing$ there is a unique element in $\text{Hom}_{\text{Med}}(M_{\text{triv}}, \mathbf{2})$, thus $M^\circ_{\text{triv}}$ is the trivial poc set $P_{\text{triv}}$ with a single element $0 = 0^*$. There exists no poc morphism $P_{\text{triv}} \to \mathbf{2}$, thus $P^\circ_{\text{triv}} = M_{\text{triv}}$.

(i) Let $M$ be the median algebra with one element. Now there are two morphisms $M \to \mathbf{2}$, thus $M^\circ \cong \mathbf{2}$. On the other hand, every poc morphism $\mathbf{2} \to \mathbf{2}$ must be the identity, thus $\mathbf{2}^\circ$ is the one element median algebra.

(ii) The new phenomena that arise in infinite median algebras and poc sets can perhaps be seen most clearly when we study linear structures coming from well ordered sets.

Recall the definition of a von Neumann ordinal, which is a transitive set $\sigma$, i.e., if $\alpha \in \sigma$ then $\alpha \subset \sigma$, such that the relation $\in$ defines a well ordering on $\sigma$, which we denote by $<$. The **successor** of $\sigma$ is $\text{suc}(\sigma) := \{\alpha \mid \alpha \leq \sigma\}$.

Consider the linear median algebra $\mathrm{Med}(\sigma)$. For any half space $H \subseteq \mathrm{Med}(\sigma)$ either $H$ or $H^* = \mathrm{Med}(\sigma) \smallsetminus H$ is an initial segment of $\sigma$, thus the dual is

$$\mathrm{Med}(\sigma)^\circ = \{H_\alpha, H_\alpha^* \mid \alpha \leq \sigma\}, \qquad \text{where } H_\alpha := \{\beta \in \sigma \mid \beta < \alpha\},$$

in other words $\mathrm{Med}(\sigma)^\circ$ is isomorphic to $\mathrm{Poc}(\sigma \smallsetminus \{0\})$. If $\sigma$ is a finite ordinal greater than $0$ then $\sigma \smallsetminus \{0\}$ is isomorphic to $\sigma - 1$; if $\sigma$ is infinite then $\sigma \smallsetminus \{0\}$ and $\sigma$ are isomorphic as well ordered sets.

Similarly, for any ultra filter $F \subset \mathrm{Poc}(\sigma)$ either $\sigma \cap F = \varnothing$, or $\min(\sigma \cap F)$ exists and then $\sigma \smallsetminus F$ is an initial segment of $\sigma$, thus

$$\mathrm{Poc}(\sigma)^\circ = \{F_\alpha \mid \alpha \leq \sigma\}, \qquad \text{where } F_\alpha := \{0^*\} \cup \{\beta \in \sigma \mid \beta \geq \alpha\},$$

in other words $\mathrm{Poc}(\sigma)^\circ$ is isomorphic to $\mathrm{Med}(\mathrm{suc}(\sigma))$.

Now we can compute the double dual map $ev \colon \mathrm{Med}(\sigma) \to \mathrm{Med}(\sigma)^{\circ\circ}$. For $\alpha, \beta < \sigma$ we have

$$H_\beta \in ev(\alpha) \quad \Leftrightarrow \quad \alpha \in H_\beta \quad \Leftrightarrow \quad \alpha < \beta \quad \Leftrightarrow \quad \mathrm{suc}(\alpha) \leq \beta \quad \Leftrightarrow \quad \beta \in F_{\mathrm{suc}(\alpha)}.$$

If we identify $\mathrm{Med}(\sigma)^{\circ\circ}$ with $\mathrm{Med}(\mathrm{suc}(\sigma \smallsetminus \{0\}))$, then $ev$ agrees with the successor map on the underlying sets. If $\sigma$ is finite, then this map is a bijection, otherwise it misses precisely the limit ordinals.

For the corresponding poc sets the situation is analogous.

$$F_\beta \in ev(\alpha) \quad \Leftrightarrow \quad \alpha \in F_\beta \quad \Leftrightarrow \quad \alpha \geq \beta \quad \Leftrightarrow \quad \mathrm{suc}(\alpha) > \beta \quad \Leftrightarrow \quad \beta \in H_{\mathrm{suc}(\alpha)}.$$

(iii) On the other hand we have median algebras whose dual is orthogonal. In the finite case those are the cubes of our introductory example.

Let $X$ be any set, and let $\mathscr{F}X$ denote the set of finite subsets of $X$. This is a median subalgebra of the power set $\mathscr{P}X$ with its Boolean median operation. For any subsets $a, b, c \in \mathscr{P}X$ we have

$$a \in [b, c] \quad \Leftrightarrow \quad b \cap c \subseteq a \subseteq b \cup c.$$

Let $H \subseteq \mathscr{F}X$ be a proper half space and assume that $\varnothing \notin H$. Then any two elements of $H$ must meet, in particular $H$ contains at most one singleton. But if all singletons are in $H^*$, then an easy induction argument shows that $H^* = \mathscr{F}X$, contradicting the assumption that $H$ is proper. Thus $H$ contains precisely one singleton and is principal.

$$\mathscr{F}X^\circ = \{\varnothing, \mathscr{F}X\} \cup \{H_x, H_x^* \mid x \in X\}, \quad \text{where } H_x := \{a \in \mathscr{F}X \mid x \in a\}.$$

In other words, $\mathscr{F}X^\circ$ is isomorphic to $\mathrm{Poc}(X)$.

Any subset $a \subseteq X$ gives rise to an ultra filter

$$F_a := \{0^*\} \cup \{x \in \mathrm{Poc}(X) \mid x \in a\} \cup \{x^* \in \mathrm{Poc}(X) \mid x \notin a\},$$

and every ultra filter in $\mathrm{Poc}(X)^\circ$ arises in that way. Thus $\mathrm{Poc}(X)^\circ$ is isomorphic to $\mathscr{P}X$ with the Boolean median operation. With this identification, the double dual map $ev \colon \mathscr{F}X \to \mathscr{F}X^{\circ\circ}$ is just the inclusion $\mathscr{F}X \subseteq \mathscr{P}X$.

(iv) Now consider the Boolean median algebra $\mathscr{P}X$ and let $H \subseteq \mathscr{P}X$ be a half space. If some $a \in \mathscr{P}X$ and its complement $X \smallsetminus a$ both lie in $H$, then $H = \mathscr{P}X$. Thus a proper half space is a transversal to $*$. It follows immediately that $\mathscr{P}X^\circ$ is orthogonal. Now suppose that $H$ is proper and $X \in H$. We claim that $H$ must be an ultra filter in $\mathscr{P}X$ in the Boolean sense. For any $a \in H$ and $b \in \mathscr{P}X$ with $a \subseteq b$ we have $b \in [a, X]$, thus $b \in H$. For any $a, b \in H$ we further have $a \cap b \in [a, b] \subseteq H$. Conversely, every Boolean ultra filter is a proper half space of $\mathscr{P}X$. Thus $\mathscr{P}X^\circ$ is isomorphic to $\mathrm{Poc}(Y)$, where $Y$ is the set of Boolean ultra filters on $\mathscr{P}X$. According to a theorem of Hausdorff, $Y$ is as large as possible: If $X$ is finite, then $|Y| = |X|$ (because, as we showed above, every ultra filter is principal), and if $X$ is infinite, then $|Y| = 2^{2^{|X|}}$.

(v) The power set $\mathscr{P}X$ also carries a poc set structure, which is much more interesting: For finite $X$ its dual is the free median algebra on $X$. We will come back to this in §6.

(vi) Recall our model example of the first section. The vertex set $VT$ of a tree carries a natural median structure. The proper half spaces in this structure are precisely the components that we obtain after cutting through an edge. Thus $VT^\circ$ is isomorphic to the poc set $\mathrm{Poc}(T) = \tilde{E}T \cup \{0, 0^*\}$ on the oriented edge set of $T$. In fact, we used the embedding $ev \colon VT \to \mathrm{Poc}(T)^\circ$ to prove that $VT$ is a median algebra.

We state a few facts about the elements of $PT^\circ$, which are easy to prove using the nesting and the finite interval condition.

(a) If an ultra filter $F \in PT^\circ$ contains a minimal element $e$, then $F = ev(\tau(e))$.

(b) If $F$ contains no minimal element, then it contains an infinite ray, i.e., an infinite sequence $w := (e_1, e_2, \ldots)$ of edges in $\tilde{E}T$, such that $\tau(e_i) = \iota(e_{i+1})$ and $e_i \neq e_{i+1}^*$ for all $i$. In this case $F$ has the form $F_w := \{0^*\} \cup \{e \in \tilde{E}T \mid e > e_i \text{ for some } i\}$.

(c) If $w_1$ and $w_2$ are two different rays, then $F_{w_1} = F_{w_2}$ iff $w_1$ and $w_2$ agree from some point onwards.

This shows that $VT^{\circ\circ}$ consists of the image of $VT$ together with equivalence classes of rays, which are the ends of $T$. □

We call a morphism of median algebras or poc sets an **embedding**, if it is an isomorphism onto its image. Thus a median morphism is an embedding, iff it is injective, and a poc morphism $f \colon P_1 \to P_2$ is an embedding, if
$$\forall a, b \in P_1 : \quad a \not\leq b \quad \Rightarrow \quad f(a) \not\leq f(b).$$

### 4.4. Representation Theorem.
  (i) $ev \colon M \to M^{\circ\circ}$ is an embedding of median algebras.
  (ii) $ev \colon P \to P^{\circ\circ}$ is an embedding of poc sets.

*Proof.* (i) To show that $ev$ is a median morphism we only need to check the definitions.
$$ev(m(x,y,z))(h) = h(m(x,y,z)) = m(h(x), h(y), h(z)) = m(ev(x), ev(y), ev(z))(h).$$

If $x \neq y \in M$, then there exists a half space separating them, thus a morphism $h \in M^\circ$ with $h(x) \neq h(y)$, hence $ev(x) \neq ev(y)$.

(ii) It is clear that
$$ev(a^*)(f) = f(a^*) = f^*(a) = ev^*(a)(f).$$
If $a \leq b$ then for all $f \in P^\circ$ we have
$$ev(a)(f) = f(a) \leq f(b) = ev(b)(f),$$
thus $ev(a) \leq ev(b)$, therefore $ev$ is monotonic. Since $ev(0)$ is the constant zero map, i.e., the zero in $P^{\circ\circ}$, it follows that $ev$ is a poc morphism. If $a \nleq b$, then $\{a, b^*\}$ is a filter base, which can be extended to an ultra filter, thus there exists an $f \in P^\circ$ with $f(a) = 1$ and $f(b) = 0$, therefore $ev(a) \nleq ev(b)$. It follows that $ev$ is injective. □

In particular, the embedding $M \to M^{\circ\circ}$ represents any median algebra as a median subalgebra of a power set, and any poc set can be represented as a family of subsets of some set. Furthermore, dualizing preserves finiteness.

### 4.5. Corollary.
*A median algebra or a poc set is finite iff its dual is finite.*

*Proof.* The dual $M^\circ$ of a finite median algebra $M$ is certainly finite. If $M^\circ$ is finite, then $M^{\circ\circ}$ must be finite, but $M$ embeds into it, so it must be finite, too. For poc sets the argument is the same. □

### Notes.
The duality theory was discovered independently by ISBELL [1980] and WERNER [1981]. In both cases it arises as a special case of more general dualities, see also JOHNSTONE [1982].

# §5. Duality for Poc Sets

Let $P$ be a poc set. For an element $a \in P$ let $V(a) := ev(a) = \{U \in P^\circ \mid a \in U\}$. More generally, for $A \subseteq P$ we define $V(A) := \bigcap \{V(a) \mid a \in A\} = \{U \in P^\circ \mid A \subseteq U\}$. The set $\mathscr{V} := \{V(a) \mid a \in P\}$ is a subbase of a topology, which we call the Stone topology on $P^\circ$.

**5.1. Proposition.**
(i) $P^\circ$ is a Stone space.
(ii) The median operation on $P^\circ$ is continuous.

*Proof.* (i) Observe that $V(a^*) = P^\circ \smallsetminus V(a)$, thus (**Stone 1**) is satisfied.
For two distinct ultra filters $U_1, U_2 \in P^\circ$ there exists an element $a \in U_1$ with $a^* \in U_2$, thus $U_1 \in V(a)$ and $U_2 \notin V(a)$, thus (**Stone 2**) is true.
For any subset $\mathscr{T} \subseteq \mathscr{V}$ there exists a subset $B \subseteq P$ such that $\mathscr{T} = \{V(b) \mid b \in B\}$. Assume that the elements of $\mathscr{T}$ meet pairwise, so for any $a, b \in B$ we have $V(a) \cap V(b) \neq \varnothing$. This means that there exists an ultra filter $U \subseteq P$ such that $a, b \in U$. In particular, $a \not\leq b^*$ for any $a, b \in B$, thus $B$ is a filter base. By the Extension Algorithm 3.3 there exists an ultra filter $V$ containing $B$, and $V \in V(B) = \bigcap \mathscr{T}$. This shows that $\mathscr{V}$ satisfies (**Stone 3**).
(ii) Let $U_1, U_2, U_3 \in P^\circ$ be ultra filters. Then $m(U_1, U_2, U_3) \in V(a)$ iff $a$ lies in at least two of the three ultra filters, thus
$$m^{-1}(V(a)) = V(a) \times V(a) \times M \cup V(a) \times M \times V(a) \cup M \times V(a) \times V(a),$$
which is open in $P^\circ \times P^\circ \times P^\circ$. □

**5.2. Remarks.**
(i) A less geometric argument goes as follows: $P^\circ$ is topologized as a subspace of $\mathscr{P}P$ with the Tychonov topology, and it is closed because the point wise limit of poc morphisms is again a poc morphism.
(ii) As a consequence of Helly's Theorem 2.2 the subbasis $\mathscr{V}$, and therefore any family of closed convex sets in the Stone topology on $P^\circ$ satisfy the following Strong Finite Intersection Property.

(**SFIP**) If $\mathscr{T}$ is a subset of $\mathscr{V}$ whose elements meet pairwise, then $\bigcap \mathscr{T}$ is non empty.

A Stone median algebra is a median algebra $M$, together with a Stone topology on $M$, such that the median operation $m: M \times M \times M \to M$ is continuous. For example, the set **2** with the canonical median structure and the discrete topology is a Stone median algebra.
A morphism of Stone median algebras is both a continuous map and a median morphism. We write $M^\bullet$ for the poc set of Stone median morphisms $M \to \mathbf{2}$, and identify its elements with the clopen half spaces in $M$.

### 5.3. Theorem.
*The categories of poc sets and Stone median algebras are anti equivalent. The evaluation maps $ev\colon P \to P^{\circ\bullet}$ and $ev\colon M \to M^{\bullet\circ}$ induce isomorphisms of poc sets and Stone median algebras, respectively.*

Before we prove this theorem we need to gather more information about Stone median algebras. Our first observation is that $M^\bullet$ separates the points of $M$.

### 5.4. Proposition.
*Let $M$ be a Stone median algebra. For any two distinct elements $x, y \in M$ there exists a clopen half space $H \in M^\bullet$ with $x \in H$ and $y \notin H$.*

*Proof.* As $M$ is a Stone space, we can certainly find disjoint clopen sets $U_0$ and $U_1$ in $M$ with $x \in U_0$, $y \in U_1$ and $M = U_0 \cup U_1$. Now we define the following relation on $M$.
$$p \sim q \quad :\Leftrightarrow \quad \forall u, v \in M : \; m(p, u, v) \in U_0 \Leftrightarrow m(q, u, v) \in U_0.$$
Clearly, this is an equivalence relation. Note that $x$ and $y$ are not equivalent, because $x = m(x, x, y)$ and $y = m(y, x, y)$ lie in different sets.

If $p \sim q$, then for any $r, s \in M$ we also have $m(p, r, s) \sim m(q, r, s)$, because by the axiom (Med 3) we have
$$m(m(p, r, s), u, v) = m(p, m(r, u, v), m(s, u, v)) \quad \text{and}$$
$$m(m(q, r, s), u, v) = m(q, m(r, u, v), m(s, u, v)),$$
which always lie in the same $U_i$, because $p \sim q$. It follows that the quotient $M' := M/\sim$ is again a median algebra.

From the lemma below we deduce that the equivalence classes of $\sim$ are open. As $M$ is quasi compact there can only be finitely many classes, thus the classes are also closed, and $M'$ is finite. Now choose any half space in $M'$ separating the classes of $x$ and $y$, then the preimage of that half space is a clopen half space in $M$ separating $x$ and $y$. $\square$

### 5.5. Lemma.
*Let $X$, $Y$ and $Z$ be topological spaces and $f\colon X \times Y \to Z$ be a continuous function. Suppose that $X$ and $Y$ are quasi compact and that $Z$ is the disjoint union of two open sets $U_0$ and $U_1$. Then every $x \in X$ has an open neighbourhood $V_x$ such that for all $y \in Y$ we have $f(x, y) \in U_0 \Leftrightarrow f(V_x \times \{y\}) \subseteq U_0$.*

*Proof.* Suppose $f(x, y) \in U_i$, then by continuity there exist open neighbourhoods $x \in V_{xy} \subseteq X$ and $y \in W_{xy} \subseteq Y$ such that $f(V_{xy} \times W_{xy}) \subseteq U_i$. By quasi compactness of $X$, the open cover $\{V_{xy} \mid x \in X\}$ has a finite subcover $\{V_{x_1 y}, \ldots, V_{x_n y}\}$. Let $W_y := W_{x_1 y} \cap \ldots \cap W_{x_n y}$ and let $V'_{xy}$ be the $V_{x_i y}$ that contains $x$. By compactness of $Y$, the open cover $\{W_y \mid y \in Y\}$ has a finite subcover $\{W_{y_1}, \ldots, W_{y_m}\}$. Let $V_x := V'_{xy_1} \cap \ldots \cap V'_{xy_m}$. Then for all $x \in X$ and $y \in Y$ we have $f(x, y) \in U_i \Leftrightarrow f(V_x \times \{y\}) \subseteq U_i$. $\square$

### 5.6. Proposition.
Let $M$ be a Stone median algebra and $C \subseteq M$ a convex subset. Then $C$ is closed iff it is a retract.

*Proof.* Intervals in Stone median algebras are compact, because $[x, y] = m(x, y, M)$ is the continuous image of a compact set. Suppose that $C \subseteq M$ is closed and convex and pick an element $x \in M$. The family $\mathscr{C} := \{C \cap [x, c] \mid c \in C\}$ consists of closed sets which are convex and meet pairwise, because $m(x, c_1, c_2) \in C \cap [x, c_1] \cap [x, c_2]$. Thanks to (**SFIP**), $\mathscr{C}$ has nontrivial intersection, and $\bigcap \mathscr{C}$ contains at most one element, which is the element of $C$ closest to $x$. This shows that $C$ is a retract.

Conversely, assume that $C \subseteq M$ is a retract. For any $x \in M \smallsetminus C$ there exists an element $c \in C$ closest to $x$. By Proposition 5.4 we find a clopen half space $H$ containing $x$ but not $c$. But as $c$ is closest to $x$, the intersection $H \cap C$ must be empty, thus $M \smallsetminus C$ must be open. □

### 5.7. Remark.
The previous proposition has an important consequence: For any median algebra $M$ one can decide whether $M$ is the dual of some poc set $P$. We will demonstrate presently that $M$ can carry at most one topology that turns it into a Stone median algebra, and then we must have $P = M^\bullet$.

Let $M_{\text{biret}}$ denote the set $M$ endowed with the topology generated by the set of half spaces $H \in M^\circ$ such that both $H$ and $H^*$ are retracts. With an argument analogous to part (ii) of Theorem 5.3 one can show that the median is always continuous with respect to this topology. In general, $M_{\text{biret}}$ need neither be Hausdorff nor quasi compact, examples are easily found amongst linear median algebras.

For Stone median algebras we have shown that a half space $H \in M^\circ$ is clopen iff both $H$ and $H^*$ are retracts. Since the family of clopen half spaces separates points, it is also a subbasis of the Stone topology on $M$. So $M$ is a Stone median algebra iff $M_{\text{biret}}$ is Hausdorff and quasi compact, and then this is the Stone topology on $M$.

For poc sets the situation is entirely different. There exist non isomorphic median algebras $M_1$ and $M_2$ such that $M_1^\circ$ and $M_2^\circ$ are isomorphic as poc sets, though of course not as topological spaces. It is easiest to give examples where both $M_1^\circ$ and $M_2^\circ$ are orthogonal, see Example 4.3(iii). But one can even find examples where $M_1$ and $M_2$ are discrete median algebras, see Example 8.3(ii). □

Let $P$ be a poc set. We now show that the subsets of $P^\circ$ of the form $V(A)$ for some $A \subseteq P$ can be characterized by the median structure on $P^\circ$ as the retracts in $P^\circ$.

### 5.8. Proposition.
A convex set $C \subseteq P^\circ$ is a retract iff $C = V(B)$ for some filter base $B \subseteq P$.

*Proof.* Let $B \subset P$ be a filter base and $F := {\uparrow}B$ the filter generated by $B$. For an ultra filter $U \subseteq P$ we define $r_B(U) := (U \smallsetminus F^*) \cup F$. This is again an ultra filter, as a straight forward application of the Extension Algorithm 3.3 shows. It contains $F$ and differs from $U$ only by elements of $F \cup F^*$. Thus if $V$ is any ultra filter that contains $B$ and hence $F$, then
$$U \cap V \subseteq U \smallsetminus F^* \subseteq r_B(U) \subseteq U \cup F \subseteq U \cup V,$$
which means that $r_B(U) \in [U, V]$. This proves that $r_B(U) \in V(B)$ is closest to $U$, hence $V(B)$ is indeed a retract.

To prove the converse, observe first that for any subset $C \subseteq P^\circ$ the set $F := \bigcap C = \{a \in P \mid C \subseteq V(a)\}$ is a filter in $P$, and $C \subseteq \bigcap \{V(a) \mid a \in F\} = V(F)$.

To show the reverse inclusion we assume that $C$ is a retract, take an ultra filter $V \in P^\circ \smallsetminus C$ and let $r(V)$ be the element of $C$ closest to $V$. Since $V \neq r(V)$ there exists an element $a \in r(V) \smallsetminus V$. We claim that $a$ belongs to $F$. Otherwise there exists an element $U \in C \cap V(a^*)$, and the entire interval $[V, U]$ is contained in $V(a^*)$. As $r(V)$ is closest to $V$, this implies that $r(V) \in V(a^*)$, contrary to our choice of $a$. It follows that every element outside $C$ is also outside $V(F) = V(B)$. □

*Proof of Theorem 5.3.* Let $P$ be a poc set. The proof of the Representation Theorem 4.4 shows that $ev \colon P \to P^{\circ\circ}$ is an embedding of poc sets, and its image is contained in $P^{\circ\bullet}$. It remains to show that the map $ev \colon P \to P^{\circ\bullet}$ is surjective. Suppose that $H \subseteq P^\circ$ is clopen, so both $H$ and $H^*$ are retracts. This implies that $H = V(A)$ and $H^* = V(B)$, where $A = \bigcap H$ and $B = \bigcap H^* = A^*$ are nonempty filters. If $A$ contains two distinct elements $a$ and $b$, say $a \not\leq b$, then there exists an ultra filter $U$ containing $a$ and $b^*$, hence $U$ lies neither in $H$ nor in $H^*$, which is a contradiction. Thus $\bigcap H$ contains a unique element, say $a$, and then $H = V(a)$. This shows that $ev \colon P \to P^{\circ\bullet}$ is surjective.

Conversely, let $M$ be a Stone median algebra and $x \in M$, then the evaluation map is defined by $ev(x) = \{H \in M^\bullet \mid x \in H\}$. Using the fact that $M^\bullet$ separates points and the proof of the Representation Theorem 4.4, we can show that $ev$ is an injective median morphism.

Let $F \in M^{\bullet\circ}$ be a filter of clopen half spaces, then for any $H_1, H_2 \in F$ we have $H_1 \not\subseteq H_2^*$, which means that $H_1 \cap H_2 \neq \varnothing$. Since $M^\bullet$ satisfies (**SFIP**), the intersection $\bigcap F$ is nonempty. If $F$ is an orientation, then $\bigcap F$ contains at most one element, because the elements of $M^\bullet$ separate points of $M$. Thus if $F$ is an ultra filter in $M^\bullet$, then there exists an element $x \in M$ such that $\bigcap F = \{x\}$, and then $F = V(x) = ev(x)$. This shows that $ev \colon M \to M^{\bullet\circ}$ is surjective.

Now we know that this $ev$ is a bijective map between compact Hausdorff spaces, so it will be a homeomorphism if it is continuous. The median algebra $M^{\bullet\circ}$ is endowed with the Stone topology, given by the subbasis $\{V(H) \mid H \in M^\bullet\}$, and the preimage of $V(H)$ under $ev$ is $H$, since
$$ev(x) \in V(H) \quad \Leftrightarrow \quad H \in ev(x) \quad \Leftrightarrow \quad x \in H.$$
This concludes the proof of our theorem. □

In this subsection we consider the double dual $M^{\circ\circ}$ with the Stone topology arising from the poc structure on $M^\circ$.

**5.9. Proposition.**
*The image of the double dual map $M \to M^{\circ\circ}$ is dense in $M^{\circ\circ}$.*

*Proof.* Let $\xi \in M^{\circ\circ}$ and $O \subseteq M^{\circ\circ}$ be a basic open set containing $\xi$, i.e.,
$$O = V(H_1) \cap \ldots \cap V(H_n) \qquad \text{for } H_1, \ldots, H_n \in M^\circ.$$
The half spaces $H_i$ must meet pairwise, for if $H_i \cap H_j = \varnothing$, then $V(H_i) \cap V(H_j) = \varnothing$. In view of Helly's Theorem 2.2 this means that there exists an element $x \in \bigcap H_i$. But then $ev(x) \in O$, thus $ev(M)$ lies dense in $M^{\circ\circ}$. □

The identification topology induced on a median algebra $M$ from the Stone topology on $M^{\circ\circ}$ via the evaluation map $M \to M^{\circ\circ}$ has the set of half spaces in $M$ as a subbase. This subbase satisfies (**Stone 1**) and (**Stone 2**), but (**Stone 3**) is equivalent to (**SFIP**), which may or may not be satisfied. In particular, every convex set is closed. We call this the convex topology on $M$.

**5.10. Corollary.**
*The double dual map is bijective iff $M$ is compact in the convex topology. In particular, every finite median algebra is naturally isomorphic to its double dual.* □

**5.11. Examples.**
(i) Recall the definition of a linear median algebra $\mathrm{med}(O)$, where $O$ is a totally ordered set. The convex sets in $\mathrm{med}(O)$ are the intervals of $O$, so the convex topology is discrete: For any $x \in O$ both $\{y \in O \mid y > x\}$ and $\{y \in O \mid y < x\}$ are half spaces. The convex topology on $\mathrm{med}(O)$ is quasi compact iff $O$ is finite.

(ii) Let $X$ be any set and $\mathscr{P}X$ the power set of $X$ with the Boolean median structure (cf. Example 1.2(iii)). Here $a \in [b,c]$ iff $b \cap c \subseteq a \subseteq b \cup c$. For any $x \in X$ the set $H_x := \{a \subseteq X \mid x \in a\}$ is a half space, and the set $\mathscr{F}X$ of finite subsets of $X$ is a convex subset of $\mathscr{P}X$. Consider the family $\mathscr{H} := \{H_x \mid x \in X\}$. Obviously any subfamily of $\mathscr{H}$ has non empty intersection, and $\bigcap \mathscr{H} = \{X\}$. In the family $\{\mathscr{F}X\} \cup \mathscr{H}$ still every finite subfamily has non empty intersection, but $\mathscr{F}X \cap \bigcap \mathscr{H} = \varnothing$ for infinite $X$. Thus $\mathscr{P}X$ equals its double dual only if $X$ is finite. □

These examples are typical, as the next theorem shows. Recall that a poc set $P$ is of type $\omega$, if every chain in $P$ and every transverse subset of $P$ is finite.

**5.12. Theorem.**
*The convex topology on a median algebra $M$ is quasi compact if and only if $M^\circ$ is of type $\omega$.*

*Proof.* If $M$ is not compact then there exists a family $\mathscr{H}$ of half spaces meeting pairwise, such that $\bigcap \mathscr{H} = \varnothing$. Obviously $\mathscr{H}$ must be infinite. For $x \in M$ let $\mathscr{H}_x := \{H \in \mathscr{H} \mid x \notin H\}$, and for the purpose of this proof we say that $x$ is special if $\mathscr{H}_x$ is finite. We claim that not all $x \in M$ can be special.

Suppose that $x \in M$ is special, then by Helly's Theorem 2.2, $\bigcap \mathscr{H}_x$ is not empty. Choose any $y \in \bigcap \mathscr{H}_x$, then every $H \in \mathscr{H}$ contains either $x$ or $y$. Now suppose that $y$ is also special and let $C := \bigcap \{H \in \mathscr{H} \mid x, y \in H\}$, thus every $H \in \mathscr{H}_x \cup \mathscr{H}_y$ meets $C$ in either $x$ or $y$. Furthermore, the elements of $\mathscr{H}_x \cup \mathscr{H}_y$ meet pairwise. By Helly's property, again, $\bigcap \mathscr{H} = \bigcap \mathscr{H}_x \cap \bigcap \mathscr{H}_y \cap C \neq \varnothing$, contradicting our hypothesis on $\mathscr{H}$ and thus proving the claim.

Let $x$ be non-special, i.e., $\mathscr{H}_x$ is an infinite family of half spaces. By Ramsey's Theorem 3.6 there exists an infinite subset $\mathscr{H}' \subseteq \mathscr{H}_x$ that is either nested or transverse. If $\mathscr{H}'$ is transverse, we are done; otherwise recall that for any two $H_1, H_2 \in \mathscr{H}_x$ both $H_1 \cap H_2$ and $H_1^* \cap H_2^*$ are non empty, hence $\mathscr{H}'$ is totally ordered.

To prove the converse, suppose that $M$ is compact in the convex topology and assume first that there is an infinite chain $\mathscr{H}$ of half spaces in $M$. It is well known that an infinite chain contains either an infinite ascending or an infinite descending sequence; in a poc set the two possibilities are of course equivalent. Suppose $H_1 \subseteq H_2 \subseteq \ldots$ is an ascending sequence of proper half spaces in $M$, then $C := \bigcup H_i$ is a convex set which meets every $H_i^*$, but the family $\{C\} \cup \{H_i^* \mid i \in \mathbb{N}\}$ has empty intersection.

Secondly, if $\mathscr{H}$ is transverse then for any two half spaces $H_1, H_2 \in \mathscr{H}$ all four intersections
$$H_1 \cap H_2, \qquad H_1^* \cap H_2, \qquad H_1 \cap H_2{,}^*, \qquad H_1^* \cap H_2^*$$
are non empty. It follows by (**SFIP**) that the intersection $\bigcap \mathscr{H}$ is non empty. Let
$$C := \{x \in M \mid x \text{ is contained in all but finitely many } H \in \mathscr{H}\},$$
so $\bigcap \mathscr{H} \subseteq C$, and $C$ is obviously convex. By the same argument as above, for any half space $H \in \mathscr{H}$ the intersection
$$H^* \cap \bigcap \{H' \in \mathscr{H} \mid H' \neq H\}$$
is non empty, hence $C$ meets every $H^*$. But the family $\{C\} \cup \{H^* \mid H \in \mathscr{H}\}$ has empty intersection, since $\mathscr{H}$ is infinite. Thus in both cases we found a family of convex sets in $M$ that fails to satisfy (**SFIP**). $\square$

Note that infinite median algebras can still be compact in the convex topology.

**5.13. Examples.**
(i) Let $X$ be any set and define the starlet on $X$ to be a median algebra with underlying set $M := \{z\} \cup X$ (for some $z \notin X$), whose median operation is determined by the rule that $m(x_1, x_2, x_3) = z$ for any three distinct elements $x_1, x_2, x_3 \in M$. Thus a subset of $M$ that contains more than one element is convex iff it contains $z$. Clearly $M^\circ$ is of type $\omega$. The convex topology on $M$ is homeomorphic to the one point compactification of $X$ with the discrete topology.

(ii) A compact median algebra need not have finite diameter. For $n \in \mathbb{N}$ let $I_n := \{0, 1, 2, \ldots, n\}$ with the linear median structure, and let $M$ be obtained from the union of all $I_n$ by gluing together all points $0 \in I_n$. Thus $M$ is the vertex set of a tree with a single vertex of infinite valency, and every branch has finite length. Clearly $M^\circ$ contains no transverse pairs and only finite chains, so $M$ is compact. □

We conclude this section with a result whose proof is also based on the central idea of Theorem 5.12. Consider the following typical examples of convex sets that are not retracts.

**5.14. Examples**
(i) Let $\omega := \{0, 1, 2, \ldots\}$ denote the first infinite von Neumann ordinal. The linear median algebra $\mathrm{Med}(\omega + 1) = \{0, 1, 2, \ldots, \omega\}$ (see Example 1.2(iv)) contains $\omega$ both as a convex subset $C$ and as an element $x$ (in fact $\mathrm{Med}(\omega + 1) = C \cup \{x\}$). But there exists no element $y \in C$ closest to $x$.

(ii) Let $X$ be an infinite set. The median algebra $\mathscr{P}X$ contains $C := \mathscr{F}X$ as a convex subset and $X$ as an element. But there exists no element $y \in C$ closest to $X$. □

The point of the next proposition is that in a non discrete median algebra we can always find one of these two examples. In fact, discreteness of median algebras can be characterized by this property, as we will see in Theorem 8.2.

**5.15. Proposition.**
*A median algebra which is not discrete contains a convex set which is not a retract.*

*Proof.* Let $M$ be a median algebra which contains two elements $x$ and $y$ such that $[x, y]$ is infinite. Assume further that every convex subset of $M$ is a retract. By Ramsey's Theorem we can find an infinite set $\mathscr{H} \subseteq \Delta(x, y)$ which satisfies one of the following conditions.

(i) $\mathscr{H}$ is nested, in particular it is a chain.

(ii) $\mathscr{H}$ is transverse.

In case (i) we consider $C := \bigcup \mathscr{H}$. Let $r = \mathrm{retr}_C(y)$, then there exists an $H \in \mathscr{H}$ with $r \in H$. As $\mathscr{H}$ is an infinite chain, there also exists a $J \in \mathscr{H}$ with $H \subset J$. Choose any $z \in J \smallsetminus H$. For $m := m(r, z, y)$ we have

(a) $m \in C$, since $r, z \in C$;

(b) $m \in H^*$, as $z, y \in H^*$; and so $m \neq r$;

(c) $m \in [r, y]$.

This is a contradiction to the definition of $r$.

In case (ii) we consider $C := \{z \in M \mid \Delta(x, z) \cap \mathscr{H} \text{ is finite}\}$. The fact that $C$ is convex follows from Lemma 2.10. We next show that $C$ contains many points.

*Claim.* Let $\mathscr{F}$ be a finite subset of $\mathscr{H}$, let $U := \bigcap\{H^* \mid H \in \mathscr{F}\}$ and $r := \mathrm{retr}_U(x)$. Then $\Delta(x,r) \cap \mathscr{H} = \mathscr{F}$, which means that $r$ lies in $C$.

Consider a half space $J$ contained in $\mathscr{H}$ but not in $\mathscr{F}$. The elements of $\mathscr{H}$ are pairwise transverse, so by Helly's Theorem 2.2 there exists an element $z \in J \cap U$. For $m = m(x,z,r)$ we have the following.

(a) $m \in J$, since $x, z \in J$.

(b) $m \in U$, as $z, r \in U$.

(c) $m \in [r, x]$.

Recall that $r \in U$ is closest to $x$, therefore $m = r$. This means that $r$ lies in $J$ for all $J \in \mathscr{H} \setminus \mathscr{F}$, which proves our claim.

Now put $r := \mathrm{retr}_C(y)$, then $\mathscr{F} := \Delta(x,r) \cap \mathscr{H}$ is finite. As $\mathscr{H}$ is infinite there exists a half space $J \in \mathscr{H} \setminus \mathscr{F}$. The set $U' := J^* \cap \bigcap\{H^* \mid H \in \mathscr{F}\}$ is convex and non empty, thus $r' := \mathrm{retr}_{U'}(x)$ exists and lies in $C \cap J^*$. Again, we consider $m = m(r, r', y) \in [r, y]$; here $J$ separates $r$ from $m$, hence $m \in C$ is closer to $y$ than $r$, contrary to the definition of $r$.

Thus in both cases we have constructed a convex subset $C$ that is not a retract. $\square$

# §6. Duality for Median Algebras

Let $M$ be a median algebra. For $x \in M$ we define
$$V(x) := ev(x) = \{H \in M^\circ \mid x \in H\}, \qquad U(a) := ev(x)^* = \{H \in M^\circ \mid x \notin H\}.$$
More generally, we define $V(X) := \{H \in M^\circ \mid X \subseteq H\}$, and $U(X) := \{H \in M^\circ \mid X \subseteq H^*\}$ for any subset $X \subseteq M$.

We consider two topologies on the set $M^\circ$. The spectral topology has the set $\mathscr{U} := \{U(x) \mid x \in M\}$ as an open subbasis, or, what amounts to the same, the set $\mathscr{V} := \{V(x) \mid x \in M\}$ as a closed subbasis; we write $M^\circ_{\text{spec}}$ if we mean $M^\circ$ with the spectral topology. The family $\{U(X) \mid X \in M, X \text{ finite}\}$ is an open basis of the spectral topology. The Stone topology has $\mathscr{S} := \mathscr{U} \cup \mathscr{V}$ as subbasis, and if we write $M^\circ$ without further specification, the Stone topology should be understood.

The next proposition follows directly from the definitions.

**6.1. Proposition.**
Let $X, Y \subseteq M$.
  (i) $V(X) = V(\text{conv}(X))$.
  (ii) If $X$ and $Y$ are convex, then $V(X) \cap V(Y) = V([X, Y])$.
  (iii) $\bigcap V(X) = \text{conv}(X)$. □

The basic open sets in the Stone topology on $M^\circ$ are the separators
$$\Delta(X, Y) = V(X) \cap U(Y) = \{H \in M^\circ \mid X \subseteq H \text{ and } Y \subseteq H^*\},$$
for finite subsets $X, Y \subseteq M$.

**6.2. Proposition.**
  (i) $M^\circ$ is a Stone space.
  (ii) For $H, J \in M^\circ$ with $H \nsubseteq J$ there exists a clopen ultra filter $F \subset M^\circ$ with $H, J^* \in F$.

*Proof.* (i) Condition (**Stone 1**) is built into the definition of $\mathscr{S}$. If $H_1$ and $H_2$ are half spaces in $M$ with $H_1 \nleq H_2$, then there exists an $x \in M$ with $x \in H_1$ and $x \notin H_2$, which means that $V(x)$ separates $H_1$ from $H_2$, hence (**Stone 2**) holds.

To check (**Stone 3**) let $\mathscr{T} \subseteq \mathscr{S}$ such that any finite subset of $\mathscr{T}$ has non trivial intersection; we need to show that $\bigcap \mathscr{T} \neq \emptyset$. In particular there exist sets $X, Y \subseteq M$ such that $\mathscr{T} = \{V(x) \mid x \in X\} \cup \{U(y) \mid y \in Y\}$. If either $X$ or $Y$ is empty, then $\bigcap \mathscr{T}$ contains an improper half space; otherwise $\bigcap \mathscr{T} = \Delta(X, Y)$. We know that half spaces separate convex sets (Theorem 2.8), therefore it is enough to show that $\text{conv}(X)$ and $\text{conv}(Y)$ are disjoint. By Corollary 2.5 a convex hull is a direct limit of convex hulls of finite sets. The hypothesis on $\mathscr{T}$ says that for any finite subsets, $A \subseteq X$ and $B \subseteq Y$ the sets $\mathrm{conv}(A)$ and $\mathrm{conv}(B)$ are disjoint, hence $\mathrm{conv}(X) \cap \mathrm{conv}(Y) = \varnothing$, as required.

(ii) If $H \nsubseteq J$ then there exists an element $x \in H \cap J^*$, hence $H, J^* \in V(x)$, and $V(x)$ is a clopen ultra filter. □

### 6.3. Remark.
A less geometric argument goes like this: Consider $M^\circ$ as a subspace of the power set $\mathscr{P}M$ with the Tychonov topology. Since the point wise limit of a median morphism is again a median morphism, $M^\circ$ is a closed subspace, hence a Stone space. □

Now we describe the objects dual to median algebras. Let $P$ be a poc set that is endowed with a Stone topology. We say that $P$ is a **Stone poc set**, if the two structures are compatible in the following sense.

(**SPS**) For all $a, b \in P$ with $a \nleq b$ there exists a clopen ultra filter $U \subset P$ with $a, b^* \in U$.

The canonical poc structure on **2** together with the discrete topology make **2** into a Stone poc set.

A morphism of Stone poc sets is both a continous map and a poc morphism. We write $P^\bullet$ for the set of Stone poc morphisms $P \to \mathbf{2}$, and identify its elements with the clopen ultra filters in $P$.

### 6.4. Theorem.
*The categories of median algebras and Stone poc sets are anti equivalent. The evaluation maps $ev\colon M \to M^{\circ\bullet}$ and $ev\colon P \to P^{\bullet\circ}$ induce isomorphisms of median algebras and Stone poc sets, respectively.*

The proof of this theorem will occupy the rest of this subsection. We begin by investigating the relation between the spectral and the Stone topology.

Let $\mathrm{K}(M^\circ)$ denote the set of clopen subsets of $M^\circ$, let $\mathrm{Ko}(M^\circ)$ denote the set of quasi compact open subsets of $M^\circ_{\mathrm{spec}}$, and $\mathrm{Kc} := \{M^\circ \smallsetminus U \mid U \in \mathrm{Ko}(M^\circ)\}$. By $\mathrm{Lat}(\mathscr{S})$ we denote the sublattice of $\mathscr{P}M^\circ$ generated by $\mathscr{S}$, i.e., the family of all sets obtained from elements of $\mathscr{S}$ by taking finite unions and intersections. Similarly we define the sublattices $\mathrm{Lat}(\mathscr{V})$ and $\mathrm{Lat}(\mathscr{U})$.

### 6.5. Lemma.
(i) $\mathrm{K}(M^\circ) = \mathrm{Lat}(\mathscr{S})$.
(ii) $\mathrm{Ko}(M^\circ) = \mathrm{Lat}(\mathscr{U})$.
(iii) $\mathrm{Kc}(M^\circ) = \mathrm{Lat}(\mathscr{V})$.

*Proof.* (i) The set $\mathrm{Lat}(\mathscr{S})$ is an open basis for the Stone topology, whose elements are closed and open. Suppose $K \subseteq M^\circ$ is open, then it is a union of elements of $\mathrm{Lat}(\mathscr{S})$. If $K$ is also closed, then it is quasi compact, so the union is finite, hence $K$ lies itself in $\mathrm{Lat}(\mathscr{S})$.

The argument for (ii) is similar, and (iii) follows from (ii). □

By part (ii) of Proposition 6.2 the spectral topology on $M^\circ$ is a $T_0$ topology. Let $\mathrm{cl}_{\mathrm{spec}}$ denote the closure operator in the spectral topology. For elements $H, J \in M^\circ$ we say $H$ is a generalization of $J \in M^\circ$, and write $H \preceq J$ if $\mathrm{cl}_{\mathrm{spec}}(J) \subseteq \mathrm{cl}_{\mathrm{spec}}(H)$. This relation is clearly reflexive and transitive, and the anti symmetry condition is just the $T_0$ axiom. Thus generalization defines a partial ordering on $M^\circ$. In fact, this ordering agrees with the poc ordering on $M^\circ$.

### 6.6. Lemma.
$$H \preceq J \quad \Leftrightarrow \quad H \leq J.$$

*Proof.*
$$\begin{aligned}
\mathrm{cl}_{\mathrm{spec}}(J) \subseteq \mathrm{cl}_{\mathrm{spec}}(H) &\Leftrightarrow J \in \mathrm{cl}_{\mathrm{spec}}(H) \\
&\Leftrightarrow \forall x \in M : H \in V(x) \Rightarrow J \in V(x) \\
&\Leftrightarrow \forall x \in M : x \in H \Rightarrow x \in J \\
&\Leftrightarrow H \subseteq J \qquad \square
\end{aligned}$$

The next lemma recovers the spectral topology from the Stone topology and the generalization order.

### 6.7. Lemma.
For a subset $V \subseteq M^\circ$ the following are equivalent.
 (i) $V$ is closed in the spectral topology.
 (ii) $V$ is closed in the Stone topology and an upper set with respect to $\leq$.

*Proof.* If $V$ is closed in $M^\circ_{\mathrm{spec}}$, then it is an intersection of elements of $\mathrm{Lat}(\mathscr{V})$, in particular it is closed in $M^\circ$ and an upper set.
Conversely, let $V \subseteq M^\circ$ be closed and $a \in \mathrm{cl}_{\mathrm{spec}}(V)$. Then any basic open neighbourhood $U(X) \ni a$, with $X \subseteq M$ finite, also meets $V$. Thus the family
$$\mathscr{T} := \{V \cup U(X) \mid X \subseteq M \text{ finite}, \ a \in U(X)\}$$
is centered. Observe that this family consists of closed subsets of $M^\circ$, therefore there exists an element $b \in \bigcap \mathscr{T}$. In particular, $b \in V$ and $a \succeq b$, hence $a \geq b$.
If $V$ is also an upper set, then $a$ must also lie in $V$, so $V = \mathrm{cl}_{\mathrm{spec}}(V)$ is closed. $\square$

### 6.8. Proposition.
For an ultra filter $F \in M^{\circ\circ}$ the following are equivalent.
 (i) $F$ is clopen in the Stone topology.
 (ii) There exists an $x \in M$ with $F = V(x)$.

*Proof.* The implication (ii) $\Rightarrow$ (i) is clear, we only need to show the converse. If $F \subseteq M^\circ$ is a closed upper set, then $F$ is also closed in the spectral topology. $M^\circ \smallsetminus F$ is open in $M^\circ_{\mathrm{spec}}$, and it is quasi compact, so $M^\circ \smallsetminus F$ lies in $\mathrm{Ko}(M^\circ) = \mathrm{Lat}(\mathscr{V})$. This means that there exists an $n \in \mathbb{N}$ and finite sets $X_1, \ldots, X_n \in M$ such that
$$F = V(X_1) \cup \ldots \cup V(X_n).$$

By the definition of ultra filters, $M^\circ$ is the disjoint union of $F$ and $F^*$. Suppose first that for some $i$ and $j$ there exists an $H \in \Delta(X_i, X_j)$. But then $H \in F \cap F^*$, which is impossible. Thus the convex hulls of the $X_i$ meet pair wise. By Helly's Theorem 2.2 the intersection $C := \mathrm{conv}(X_1) \cap \ldots \cap \mathrm{conv}(X_n)$ is nonempty.

Assume next that $C$ contains two distinct elements, say $x$ and $y$. Choose an $H \in \Delta(x, y)$, then neither $H$ nor $H^*$ lies in any $V(X_i)$, thus $H \notin F \cup F^*$, which is again a contradiction. We conclude therefore that $C$ contains a single element, say $x$.

Now we claim that in fact $F = V(x)$. If $H \in F$, then $H$ contains some $X_i$, so $x \in C \subseteq \mathrm{conv}(X_i) \subseteq H$, thus $F \subseteq V(x)$. If $H \notin F$, then $H^*$ meets all $X_i$. By Helly's Theorem it follows that $H^* \cap \mathrm{conv}(X_1) \cap \ldots \cap \mathrm{conv}(X_n)$ is nonempty, but that implies that $x \in H^*$. This concludes the proof of our proposition. □

Now we have found $\mathscr{V} = ev(M)$, an isomorphic copy of the median algebra $M$, inside $M^{\circ\circ}$ as the set $M^{\circ\bullet}$ of clopen ultra filters in $M^\circ$. This proves one half of Theorem 6.4.

For the converse we start by recalling some facts about retractions, see Proposition 2.12 and Example 2.13(iii). Let $r_1, r_2 \colon M \to M$ be two retractions.

(i) Let $P$ be a poc set and $X \subseteq P^\circ$ a finite set of ultra filters. Then $\mathrm{conv}(X)$ is a retract with retraction $r(z) := \bigcap X \cup (z \cap \bigcup X)$.

(ii) The images $r_1(M)$ and $r_2(M)$ meet iff $r_1$ and $r_2$ commute, and then $r_1 \circ r_2$ is a retraction onto $r_1(M) \cap r_2(M)$.

### 6.9. Lemma.

Let $P$ be a poc set and $X, Y \subseteq P^\circ$ two finite sets of ultra filters. Then

$$\mathrm{conv}(X) \cap \mathrm{conv}(Y) \neq \varnothing \quad \Leftrightarrow \quad \bigcap X \cap \bigcap Y^* = \varnothing.$$

If $\mathrm{conv}(X) \cap \mathrm{conv}(Y) \neq \varnothing$, then this intersection meets the median algebra generated by $X \cup Y$.

*Proof.*

Let $r_X(x) = \bigcap X \cup (x \cap \bigcup X)$ be the retraction onto $\mathrm{conv}(X)$, and similarly $r_Y$. Then by (i) above

$$r_Y \circ r_X(x) = \bigcap Y \cup (\bigcup Y \cap \bigcap X) \cup (x \cap \bigcup X \cap \bigcup Y).$$

Our hypothesis $\bigcap X \cap \bigcap Y^* = \varnothing$ implies that $\bigcap X \subseteq P^\circ \smallsetminus \bigcap Y^* = \bigcup Y$, so

$$r_Y \circ r_X(x) = \bigcap Y \cup \bigcap X \cup (x \cap \bigcup X \cap \bigcup Y).$$

Applying $*$ to our hypothesis, we obtain $\bigcap X^* \cap \bigcap Y = \varnothing$, which in turn implies that $\bigcap Y \subseteq \bigcup X$. Thus the above formula for $r_Y \circ r_X$ is indeed symmetric in $X$ and $Y$, in other words, $r_X$ and $r_Y$ commute. According to (ii) above the two retracts $\mathrm{conv}(X)$ and $\mathrm{conv}(Y)$ must meet.

Moreover, for any $x \in X \cup Y$ the element $r_X \circ r_Y(x)$ lies in this intersection and is a median expression in terms of elements of $X$ and $Y$.

Conversely, if $x \in \bigcap X \cap \bigcap Y^*$, then the half space $V(x) \subseteq P^\circ$ separates $X$ and $Y$, so conv($X$) and conv($Y$) can not meet. □

For arbitrary sets $A, B \subseteq P^\circ$ let $I(A,B) := \bigcap A \cap \bigcap B^*$. To get information about $I(A,B)$ for infinite sets, we need to invoke compactness arguments.

**6.10. Proposition.**
Let $P$ be a Stone poc set. Let $A, B$ be convex subsets of $P^\bullet$.
   (i) If $A \cap B = \varnothing$, then $I(A,B) \neq \varnothing$.
   (ii) If $A \cup B = P^\bullet$, then $I(A,B)$ contains at most one element.

*Proof.* (i) Let $X \subseteq A$ and $Y \subseteq B$ be finite subsets. If $I(X,Y) = \varnothing$, then the intersection $\text{conv}_{P^\circ}(X) \cap \text{conv}_{P^\circ}(Y)$ contains an element $x$ of the median algebra generated by $X \cup Y$, in particular $x$ is a clopen ultra filter, so $x \in P^\bullet$. In view of Proposition 2.6 we have $x \in \text{conv}_{P^\bullet}(X) \cap \text{conv}_{P^\bullet}(Y) \subseteq A \cap B$. But that contradicts $A \cap B = \varnothing$.
Now the family $\{I(X,Y) \mid X \subseteq A \text{ finite}, Y \subseteq B \text{ finite}\}$ is centered and consists of closed subsets of $P$, so thanks to (FIP) its intersection, $I(A,B)$, is nonempty.
(ii) Suppose that $a, b \in P$ are distinct, say $a \nleq b$. By (**SPS**) there exists an element $x \in P^\bullet$ with $a \in x$ and $b \notin x$. If $x \in A$ then $b$ can't be in $I(A,B)$, and if $x \in B$ then $a$ can't be in $I(A,B)$. Thus $I(A,B)$ does not contain two distinct elements. □

Now we are ready to prove the second half of Theorem 6.4. Let $P$ be a Stone poc set. The evaluation map $ev \colon P \to P^{\bullet\circ}$ maps an element $a \in P$ to the half space $V(a) := \{F \in P^\bullet \mid a \in F\}$.

   (i) $ev$ is an embedding of poc sets.

This follows from the proof of the Representation Theorem 4.4, if we use (**SPS**) instead of the Extension Algorithm 3.3 to show that $a \nleq b$ implies $ev(a) \nleq ev(b)$.

   (ii) $ev$ is surjective.

Let $H \subseteq P^\bullet$ be a half space. By Proposition 6.10 there exists a unique element $a \in P$ with $I(H, P^\bullet \smallsetminus H) = \{a\}$, and obviously $H = V(a)$.

   (iii) $ev$ is a homeomorphism.

We already know that $ev$ is a bijective map between compact Hausdorff spaces, thus it is enough to show that it is continuous. The poc set $P^{\bullet\circ}$ is endowed with the Stone topology given by the subbasis $\{V(F), U(F) \mid F \in P^\bullet\}$. The preimage of these under $ev$ is easily seen to be $F$ and $F^*$, since

$$ev(a) \in V(F) \iff F \in ev(a) \iff a \in F,$$
$$ev(a) \in U(F) \iff F \notin ev(a) \iff a \notin F \iff a \in F^*,$$

and they are both clopen by definition of $P^\bullet$.

This concludes the proof of Theorem 6.4. □

In this subsection we consider the double dual $P^{\circ\circ}$ with the Stone topology arising from the median structure on $P^\circ$

### 6.11. Proposition.
*The image of the double dual map $P \to P^{\circ\circ}$ is dense.*

*Proof.* Let $\xi \in P^{\circ\circ}$ and let $O$ be a basic open neighbourhood, i.e.,

$$O = V(A) \cap U(B) \qquad \text{for some finite sets } A, B \subseteq P^\circ \text{ of ultra filters.}$$

We want to show that $O$ meets the image of $ev$. Let $F_A := \bigcap A$ and $F_B := \bigcap B$, both are filters in $P$. If there is an element $a \in F_A$ with $a^* \in F_B$, then $ev(a) \in O$ as required. If we can't find such an element, then $F_A \cup F_B$ is a filter in $P$ and can be extended to an ultra filter. However, this ultra filter corresponds to an element of $P^\circ$ that lies in the convex hull of both $A$ and $B$, contradicting the fact that these sets can be separated by $\xi$. □

The topology on $P$ induced via the evaluation map $ev\colon P \to P^{\circ\circ}$ from the Stone topology on $P^{\circ\circ}$ has the set $\{U, U^* \mid U \in P^\circ\}$ as a clopen subbasis. We call this the filter topology on $P$, because all filters and ideals are closed subsets; by the Extension Algorithm this topology is always Hausdorff.

### 6.12. Corollary.
*The double dual map $ev\colon P \to P^{\circ\circ}$ is bijective iff $P$ is quasi compact in the filter topology. In particular, every finite poc set is naturally isomorphic to its double dual.* □

The next theorem is the analogue of Theorem 5.12, only the situation here is somewhat simpler.

### 6.13. Theorem.
*The filter topology on a poc set $P$ is quasi compact if and only if $P$ is of type $\omega$.*

*Proof.* Assume that $P$ is not quasi compact in the filter topology, then there exists a set $\mathscr{V}$ of ultra filters and ultra ideals, such that any finite subset of $\mathscr{V}$ has non empty intersection, but $\bigcap \mathscr{V} = \varnothing$. Here $\mathscr{V}$ must contain both a filter $F$ and an ideal $I$, and the intersection $F \cap I$ must be infinite, otherwise the assumption $\bigcap \mathscr{V} = \varnothing$ can't be satisfied. Now the Ramsey Argument (Corollary 3.7) and the subsequent remark yield an infinite subset $S \subseteq F \cap I$ that is either transverse or totally ordered.

To prove the converse suppose first that $P$ contains an infinite chain $A$. Without loss of generality we may assume that $A$ has no maximal element, for otherwise every ascending and every descending sequence in $A$ would be finite, which would mean that $A$ is finite. Consider the filters $F_a := \{b \in P \mid b \geq a\}$ and the ideal $I := \{b \in P \mid b \leq a \text{ for some } a \in A\}$. If $S$ is a finite subset of $A$, then the intersection $I \cap \bigcap \{F_a \mid a \in S\}$ is not empty, because it contains the maximal element of $S$. On the other hand, $I \cap \bigcap \{F_a \mid a \in A\} = \varnothing$ because an element of $\bigcap \{F_a \mid a \in A\}$ can't be smaller than any element of $A$.

Next assume that $A \subseteq P$ is an infinite transverse set. Consider the filters $F_a := \{b \in P \mid b \geq c$ for some $c \in A$, but $b \neq a\}$ and the ideal $I := \{b \in P \mid b \leq a$ for some $a \in A\}$. If $S \subseteq A$ is a finite subset, then there still is an element $b \in A \smallsetminus S$, and $b \in I \cap \bigcap \{F_a \mid a \in S\}$. If $b \in I \cap \bigcap \{F_a \mid a \in A\}$, then there exist $a, a' \in A$ with $a \leq b \leq a'$. As $A$ is transverse, this can only happen when $a = a'$, but then $b \in A$, which contradicts the definition of the $F_a$.

Now we have shown that in both cases there exists a family of closed sets in the filter topology which fails to satisfy the finite intersection property, thus $P$ can't be quasi compact. □

Again, there exist examples of infinite compact poc sets.

**6.14. Example.** Given a set $X$, the starlet poc set on $X$ is defined as the set $\{0, 0^*\} \cup X \cup X^*$, with the usual bijection $^*\colon X \to X^*$, and we define an order where for distinct elements $x, y \in X$ we always have $x < y^*$. It is easy to see that starlet poc sets are precisely the duals of starlet median algebras, see Example 5.13(i). □

This example is typical, as the next theorem shows.

**6.15. Theorem.**
*A poc set that is quasi compact in the filter topology is either finite or contains an infinite starlet as a sub poc set.*

*Proof.* Assume that $P$ is an infinite poc set, which is compact in the filter topology and does not contain an infinite starlet. Choose an ultra filter $U \subseteq P$. By Theorem 6.13 every chain and every transverse subset of $U$ is finite. In particular, $U$ contains minimal elements. Let $A_0$ denote the set of all minimal elements of $U$. We claim that $A_0$ is finite. Two distinct elements $a, b \in A_0$ can't be comparable, and $a \not\leq b^*$, because $a$ and $b$ lie in an ultra filter. Only two possibilities remain,
  (i) $a < b^*$, or
  (ii) $a \pitchfork b$.
If $A_0$ is infinite, then Ramsey's Theorem yields an infinite poc set $B \leq A_0$, where all pairs of elements satisfy the same relation. This means that $B$ is either an infinite starlet or an infinite transverse subset. Both possibilities are excluded by our hypothesis, so $A_0$ must be finite.

Now we construct a sequence $A_0, A_1, A_2, \ldots$, where $A_i$ comprises the minimal elements of $U \smallsetminus (A_0 \cup \ldots \cup A_{i-1})$. By the argument just given, each $A_i$ is non empty and finite, and since $U$ is infinite, the sequence is infinite, too. For each $a \in A_i$, $i \in \mathbb{N}$, there exists an element $b \in A_{i-1}$ with $b < a$, thus for each $i$ we can choose a function $f_i \colon A_i \to A_{i-1}$ with $f(a) < a$ for all $a \in A_i$. Using König's Lemma, below, we find an infinite chain in $U$, which is a final contradiction. □

**6.16. Königs Lemma.** KÖNIG [1932]
*For every sequence of functions*

$$A_0 \xleftarrow{f_1} A_1 \xleftarrow{f_2} A_2 \xleftarrow{f_3} A_3 \longleftarrow \ldots,$$

*where the $A_i$ are finite, non empty sets, there exist an infinite sequence $(a_i \in A_i)$ with*

$$a_0 \xleftarrow{f_1} a_1 \xleftarrow{f_2} a_2 \xleftarrow{f_3} a_3 \longleftarrow \ldots \qquad \square$$

Let us consider some special cases. Recall that a binary poc set has the form $P = \{0, 0^*\} \cup O \cup O^*$, where $O$ is a partially ordered set and no element of $O$ is comparable to any element of $O^*$. If two proper elements $a, b \in P$ satisfy $a < b^*$, then one of them must lie in $O$ and the other in $O^*$, hence a binary poc set can't contain a starlet on three elements. Thus we have

**6.17. Corollary.**
*A binary poc set is compact in the filter topology if and only if it is finite.* $\qquad \square$

Observe that Theorem 6.15 for binary poc sets amounts to the fact that a partially ordered set in which all chains and all anti chains are finite must also be finite.

**6.18. Corollary.**
*An interval is compact in the convex topology iff it is finite.*

*Proof.*
The dual of an interval $[x, y]$ is a binary poc set, as $[x, y]^\circ = \{0, 0^*\} \cup \Delta(x, y) \cup \Delta(y, x)$. $\qquad \square$

CONGRUENCE RELATIONS

In this subsection we make a slight detour to study congruence relations on a median algebra. An equivalence relation $\sim$ on $M$ is said to be a **congruence relation**, if it satisfies

$$\forall x, y, u, v \in M : \quad x \sim y \quad \Rightarrow \quad m(x, u, v) \sim m(y, u, v).$$

In particular, the quotient $M/\sim$ carries a natural median structure.
We have already met a congruence relation briefly in the proof of Proposition 5.4. For general information about congruence relations see BANDELT-HEDLÍKOVÁ [1983].
Here we view the quotient map $M \to M/\sim$ as a contraction along a certain set of hyper planes. The question, which sets of hyper planes arise from congruence relations, will lead us back to the Stone topology on $M^\circ$.
Given a congruence relation $R \subseteq M \times M$ we define

$$\triangle(R) := \bigcup \{\Delta(x, y) \subset M^\circ \mid (x, y) \in R\}.$$

The set $U := \triangle(R) \subseteq M^\circ$ satisfies the following three conditions.

(**Con 1**)    $\varnothing, M \notin U$.

(**Con 2**)   $U = U^*$.

(**Con 3**)   $U$ is an open subset of $M^\circ$.

The set $\overline{U} := \{\overline{H} \mid H \in U\}$ is the set of hyperplanes mentioned above, but for technical reasons we prefer to work with the half spaces in $U$.

Now let $U \subseteq M^\circ$ satisfy conditions (**Con 1**) and (**Con 2**), and set

$$\nabla(U) := \{(x, y) \in M \times M \mid \Delta(x, y) \subseteq U\}.$$

Then $\nabla(U)$ clearly is a reflexive and symmetric relation. Lemma 2.16, the triangle equality, implies that it is also transitive. The fact that $\nabla(U)$ is a congruence relation follows immediately from the next lemma.

### 6.19. Lemma.
For all $x, y, u, v \in M$ we have $\overline{\Delta}(m(x, u, v), m(y, u, v)) = \overline{\Delta}(x, y) \cap \overline{\Delta}(u, v)$.

*Proof.* For a half space $H \in M^\circ$ the following are equivalent.

$\overline{H} \in \overline{\Delta}(m(x, u, v), m(y, u, v))$

$\Leftrightarrow \quad m(x, u, v) \in H \quad \text{xor} \quad m(y, u, v) \in H$

$\Leftrightarrow \quad (H \text{ contains at least two out of } x, u, v) \quad \text{xor}$
$\quad\quad (H \text{ contains at least two out of } y, u, v)$

$\Leftrightarrow \quad (x \in H \text{ xor } y \in H) \quad \text{and} \quad (u \in H \text{ xor } v \in H)$

$\Leftrightarrow \quad \overline{H} \in \overline{\Delta}(x, y) \cap \overline{\Delta}(u, v)$   □

### 6.20. Examples.
Let $C \subseteq M$ be a convex subset. Let $H \subseteq M$ be a half space. We say that $H$ cuts $C$, if both $C \cap H$ and $C \cap H^*$ are non empty. We say that $H$ touches $C$, if $C \cap H = \emptyset$, and $H$ is maximal with respect to this property.

Now consider the two congruence relations

$$\nabla\{H \in M^\circ \mid H \text{ cuts } C\}, \quad \nabla\{H, H^* \in M^\circ \mid H \text{ does not touch } C\}.$$

It is an easy exercise to show that these are the the smallest and largest congruence relation that have $C$ as an equivalence class.   □

### 6.21. Theorem.
*The operations $\triangle$ and $\nabla$ define a bijection between the set of congruence relations on $M$ and the family of subsets $U \subseteq M^\circ$ satisfying (**Con 1**), (**Con 2**) and (**Con 3**).*

The proof of this theorem comprises the next two propositions.

### 6.22. Proposition.
*For any congruence relation $R \subseteq M \times M$ we have $\nabla \circ \triangle(R) = R$.*

*Proof.* If $x \underset{R}{\sim} y$ then $\Delta(x, y) \subseteq \triangle(R)$, hence $x \underset{\nabla \triangle(R)}{\sim} y$.

Conversely, suppose that $x \underset{\nabla\triangle(R)}{\sim} y$ but $x \underset{R}{\not\sim} y$. Let $C$ denote the $R$-class of $x$. Observe that $C$ is convex, because if $u, v \in C$ and $w \in [u, v]$, then
$$w = m(u, v, w) \underset{R}{\sim} m(x, v, w) \underset{R}{\sim} m(x, x, w) = x.$$
By Theorem 2.7 there exists a half space $H \in M^\circ$ such that $y \in H$ and $H$ touches $C$. In particular $H \in \Delta(y, x)$ and hence by our assumption $H \in \triangle(R)$. This means that there exist $u, v \in M$ with $u \underset{R}{\sim} v$ and $H \in \Delta(v, u)$.

Let $D := \mathrm{conv}(H \cup \{u\})$. If $C \cap D = \varnothing$ then there exists a half space $J \in \Delta(D, C)$, contradicting the maximality of $H$. Thus $C$ and $D$ meet, i.e., there exist $h \in H$, $c \in C$ with $c \in [u, h]$. But then $H$ must meet $C$, because
$$x \underset{R}{\sim} c = m(h, c, u) \underset{R}{\sim} m(h, c, v) \in H,$$
as $h, v \in H$. This contradicts our choice of $H$, thus we must have $x \underset{R}{\sim} y$. □

**6.23. Proposition.**
Let $U \subseteq M^\circ$ satisfy (**Con 1**) and (**Con 2**). Then $\triangle \circ \nabla(U)$ is the interior of $U$ with respect to the Stone topology on $M^\circ$.

*Proof.* Obviously $\triangle \circ \nabla(U)$ is an open set contained in $U$. Conversely, to show that $H \in \triangle \circ \nabla(U)$ it is enough to find $x, y \in M$ such that $H \in \Delta(x, y) \subseteq U$. If $H$ lies in the interior of $U$ then there exists a basic open neighbourhood $\Delta(X, Y)$ with $X, Y \subseteq M$ finite and $H \in \Delta(X, Y) \subseteq U$. By (**Con 1**) neither $X$ nor $Y$ is empty. As $X$ and $Y$ are finite, $\mathrm{conv}(X)$ and $\mathrm{conv}(Y)$ are retracts. Choose a gate $(x, y)$ for these two sets, then by Lemma 2.14, $H \in \Delta(X, Y) = \Delta(\mathrm{conv}(X), \mathrm{conv}(Y)) = \Delta(x, y)$ as required. □

In the dual of a quotient we expect to see precisely the hyper planes that we didn't contract, and this is in fact the case. We will prove this in Proposition 7.13 in the next section.

## §7. Duality for Maps

Let $f\colon M_1 \to M_2$ be a median morphism. Then we have a natural induced morphism $f^\circ\colon M_2^\circ \to M_1^\circ$, called the **dual** of $f$. If we interprete $M^\circ$ as $\mathrm{Hom}_{\mathrm{Med}}(M, \mathbf{2})$, then $f^\circ\varphi := \varphi \circ f$, and if we view $M^\circ$ as a set of half spaces, then $f^\circ(H) := f^{-1}(H)$. Since $f^{-1}$ preserves inclusions and complements, $f^\circ$ is a morphism of poc sets. Of course this construction is functorial, i.e., $f_1^\circ \circ f_2^\circ = (f_2 \circ f_1)^\circ$, and $\mathrm{id}_M^\circ = \mathrm{id}_{M^\circ}$.

Similarly, for a poc morphism $g\colon P_1 \to P_2$ we have an induced **dual** dual map $g^\circ\colon P_2^\circ \to P_1^\circ$, defined by $g^\circ(U) := g^{-1}(U)$ for every ultra filter $U \in P_2^\circ$, and as $g^{-1}$ preserves unions and intersections, $g^\circ$ is also a median morphism.

### 7.1. Lemma.
*The following diagrams commute.*

$$\begin{array}{ccc} M_1^{\circ\circ} & \xrightarrow{f^{\circ\circ}} & M_2^{\circ\circ} \\ {\scriptstyle ev}\uparrow & & \uparrow{\scriptstyle ev} \\ M_1 & \xrightarrow{f} & M_2 \end{array} \qquad \begin{array}{ccc} P_1^{\circ\circ} & \xrightarrow{g^{\circ\circ}} & P_2^{\circ\circ} \\ {\scriptstyle ev}\uparrow & & \uparrow{\scriptstyle ev} \\ P_1 & \xrightarrow{g} & P_2 \end{array}$$

*Proof.* We only verify this for the first diagram, the second is exactly analogous. For all $x \in M_1$ we have

$$\begin{aligned} f^{\circ\circ}(ev(x)) &= (f^\circ)^{-1}(ev(x)) \\ &= \{\, H \in M_2^\circ \mid f^\circ(H) \in ev(x) \,\} \\ &= \{\, H \in M_2^\circ \mid x \in f^{-1}(H) \,\} \\ &= \{\, H \in M_2^\circ \mid f(x) \in H \,\} \\ &= ev(f(x)). \end{aligned}$$
$\square$

### 7.2. Corollary.
*The dual maps $f^\circ\colon M_2^\circ \to M_1^\circ$ and $g^\circ\colon P_2^\circ \to P_1^\circ$ are continuous in the Stone topology.*

*Proof.* The subbasis of the Stone topology on $M_2^\circ$ is the family $\{\, ev(x), ev(x)^* \subset M_2^\circ \mid x \in M_2 \,\}$, and the preimage of a set in this family under $f^\circ$ lies in the subbasis of $M_1^\circ$, because $(f^\circ)^{-1}(ev(x)) = f^{\circ\circ}(ev(x)) = ev(f(x))$, and $(f^\circ)^{-1}(ev(x)^*) = ev(f(x))^*$.
The proof for poc sets is analogous. $\square$

### 7.3. Remark.
Let $X$ be a set. In the following we consider the power set $\mathscr{P}X$ as a Boolean poc set (see Example 1.4(iii)). For any $x \in X$ let $\pi(x) := \{a \in \mathscr{P}X \mid x \in a\}$ denote the principal ultra filter generated by $x$. Thus we have an injection $\pi\colon X \to \mathscr{P}X^\circ$.

Let $M$ be a median algebra and $\varphi\colon X \to M$ be a function. The map
$$\varphi^{\#}\colon M^\circ \to \mathscr{P}X, \qquad H \mapsto \varphi^{-1}(H),$$
is in fact a poc morphism. By a computation analogous to the proof of Lemma 7.1 we can show that the following diagram commutes.

$$\begin{array}{ccc} \mathscr{P}X^\circ & \xrightarrow{(\varphi^\#)^\circ} & M^{\circ\circ} \\ \pi \uparrow & & \uparrow ev \\ X & \xrightarrow{\varphi} & M \end{array}$$

In particular we have a poc morphism $\pi^{\#}\colon \mathscr{P}X^{\circ\circ} \to \mathscr{P}X$. This is a left inverse to the double dual map $ev\colon \mathscr{P}X \to \mathscr{P}X^{\circ\circ}$; in other words, the composition of the following two maps yields the identity on $\mathscr{P}X$.
$$\mathscr{P}X \xrightarrow{ev} \mathscr{P}X^{\circ\circ} \xrightarrow{\pi^{\#}} \mathscr{P}X.$$
To check this observe that for all $x \in X$ and $a \in \mathscr{P}X$
$$x \in \pi^{\#} \circ ev(a) \quad \Leftrightarrow \quad \pi(x) \in ev(a) \quad \Leftrightarrow \quad a \in \pi(x) \quad \Leftrightarrow \quad x \in a. \qquad \square$$

As an application of the duality for maps we give a description of the free median algebra on a set $X$. Let $\mathrm{Med}(X)$ denote the median subalgebra of $\mathscr{P}X^\circ$ which is generated by the set $\pi(X)$. Let $\iota\colon X \to \mathrm{Med}(X)$ be defined by $\iota(x) := \pi(x)$, i.e., $\iota$ is the same as $\pi$, only with a restricted target. This careful distinction is necessary here, because the dual maps have different domains of definition.

**7.4. Theorem.**
(i) $\mathrm{Med}(X)$ *is the free median algebra on the set $X$. In other words, for every median algebra $M$ and every function $\varphi\colon X \to M$ there exists a unique median morphism $f\colon \mathrm{Med}(X) \to M$ with $\varphi = f \circ \iota$.*

$$\begin{array}{ccc} \mathrm{Med}(X) & \xrightarrow{f} & M \\ & \iota \nwarrow \quad \nearrow \varphi & \\ & X & \end{array}$$

(ii) $\iota^{\#}\colon \mathrm{Med}(X)^\circ \to \mathscr{P}X$ *is an isomorphism.*

(iii) $\mathrm{Med}(X) = \mathscr{P}X^\circ$ *iff $X$ is finite.*

*Proof.* (i) On the generating set $\iota(X)$ the map $f\colon \mathrm{Med}(X) \to M$ is determined by $\varphi$, thus there is at most one such function. It remains to show that $f$ is a well defined median morphism. We deduce from the previous remark that $(\varphi^\#)^\circ$ maps $\iota(X)$, and therefore $\mathrm{Med}(X)$, into the image of $M$ under $ev$. As the latter map is an isomorphism onto its image we have a median morphism $f := ev^{-1} \circ (\varphi^\#)^\circ$, which extends $\varphi$ as desired.

(ii) The function $\iota^{\#}\colon \mathrm{Med}(X)^\circ \to \mathscr{P}X$ is a morphism of poc sets by the previous remark. Applying the universal property (i) of the free median algebra $\mathrm{Med}(X)$ in the case of the target $M := \mathbf{2}$, we see that $\iota^{\#}$ is bijective. To show that $\iota^{\#}$ is in fact an embedding it is enough to show that the inverse map is a poc morphism.

Let $j\colon \operatorname{Med}(X) \hookrightarrow \mathscr{P}X^\circ$ be the inclusion. Dualizing the diagram

$$\begin{array}{ccc} \operatorname{Med}(X) & \xrightarrow{j} & \mathscr{P}X^\circ \\ & \searrow^{\iota} \quad \swarrow_{\pi} & \\ & X & \end{array}$$

we get a commutative diagram

$$\begin{array}{ccc} \operatorname{Med}(X)^\circ & \xleftarrow{j^\circ} & \mathscr{P}X^{\circ\circ} \\ & \searrow^{\iota^\#} \quad \swarrow_{\pi^\#} & \\ & \mathscr{P}X & \end{array}$$

i.e., $\pi^\# = \iota^\# \circ j^\circ$. Composing with $ev\colon \mathscr{P}X \to \mathscr{P}X^{\circ\circ}$ we obtain

$$\operatorname{id}_{\mathscr{P}X} = \iota^\# \circ j^\circ \circ ev,$$

as desired.

(iii) Consider the two maps

$$\mathscr{P}X^\circ \xrightarrow{(\iota^\#)^\circ} \operatorname{Med}(X)^{\circ\circ} \xleftarrow{ev} \operatorname{Med}(X)\,.$$

Since $\iota^\#$ is an isomorphism, so is $(\iota^\#)^\circ$. By Theorem 5.12 the map $ev$ is an isomorphism iff $\mathscr{P}X$, the dual of $\operatorname{Med}(X)$, is of type $\omega$, i.e., contains no infinite chains and no infinite transverse sets. Clearly, this holds iff $X$ is finite. □

In particular, the free median algebra generated by a finite set is itself finite. Thus we have the following corollary.

**7.5. Corollary.**
*A finitely generated median algebra is finite.* □

**7.6. Examples.**
It is instructive to work out the median graphs for the first few free median algebras.

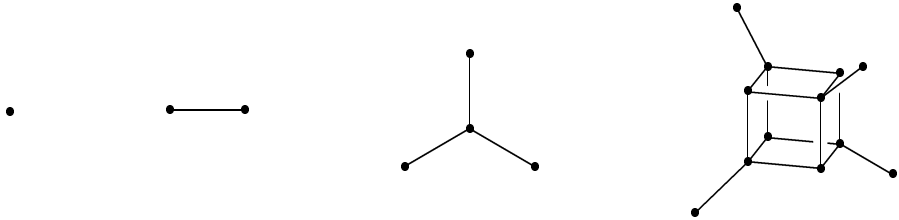

We will say more about the free median algebra on five generators in the next section.

**7.7. Remark.**
Perhaps not surprisingly, the free poc set on a set $X$ is simply $\operatorname{Poc}(X)$. The universal property is trivial to verify. □

Next we study how the properties of injectivity and surjectivity are exchanged by duality. The next result is as we expect it (and analogous to a similar result for Boolean algebra homomorphisms and continuous maps on Stone spaces), though the proofs are not always routine, because we don't have a full duality.

Recall that a poc morphism $g\colon P_1 \to P_2$ is an embedding, if it is an isomorphism onto its image; in particular, $g$ has to satisfy

$$\forall a, b \in P_1: \quad g(a) \le g(b) \quad \Rightarrow \quad a \le b.$$

**7.8. Proposition.**

Let $f\colon M_1 \to M_2$ be a median morphism and $g\colon P_1 \to P_2$ a poc morphism.
  (i) $g$ is an embedding iff $g^\circ$ is surjective.
  (ii) $g$ is surjective iff $g^\circ$ is injective.
  (iii) $f$ is injective iff $f^\circ$ is surjective.
  (iv) If $f$ is surjective then $f^\circ$ is an embedding.

*Proof.* (i) Suppose first that $g$ is not an embedding, i.e., there exist $a, b \in P_1$ such that $a \not\le b$ but $g(a) \le g(b)$. Then $\{a, b^*\}$ is a filter base, which can be extended to an ultra filter $U \in P_1^\circ$. If $U$ is in the image of $g^\circ$, then there exists an ultra filter $V \in P_2^\circ$ with $U = g^{-1}(V)$, therefore $\{g(a), g(b)^*\} \subset V$, but that is impossible. Hence $g^\circ$ is not surjective.

Now suppose that $g$ is an embedding, then for any ultra filter $U \in P_1^\circ$ the image $g(U)$ must be a filter base in $P_2$, hence there exists an ultra filter $V \in P_2^\circ$ containing it, and $U = g^{-1}(V)$.

(ii) Let $v_1, v_2\colon P_2 \to \mathbf{2}$ be two poc morphisms. If $g$ is surjective, then it can be cancelled on the right, i.e.,

$$v_1 \circ g = v_2 \circ g \quad \Rightarrow \quad v_1 = v_2,$$

which means that $g^\circ$ is injective.

Conversely, assume that $g$ is not surjective and choose an element $a \in P_2$ not in the image of $g$. The set $F := \{b \in P_2 \mid b > a \text{ or } b > a^*\}$ is a filter base in $P_2$. Extend $F$ to an ultra filter $U$ in $P_2 \smallsetminus \{a, a^*\}$. Then both $U \cup \{a\}$ and $U \cup \{a^*\}$ are ultra filters in $P_2$, whose preimages under $g$ are the same, thus $g^\circ$ is not injective.

(iii) Suppose $f$ is not injective, so there exist $x_1, x_2 \in M_1$ with $f(x_1) = f(x_2)$. Choose a half space $H \in M_1^\circ$ separating $x_1$ and $x_2$, then $H$ is not the preimage under $f$ of any half space in $M_2$, hence $f^\circ$ can not be surjective.

Conversely, suppose that $f\colon M_1 \to M_2$ is injective. A half space $H \in M_1^\circ$ is the preimage under $f$ of some half space $J \in M_2^\circ$ iff $f(H)$ and $f(H^*)$ can be separated by $J$, and by Theorem 2.8 such a $J$ exists iff

$$\operatorname{conv}(f(H)) \cap \operatorname{conv}(f(H^*)) = \varnothing.$$

Assume therefore that $\operatorname{conv}(f(H))$ and $\operatorname{conv}(f(H^*))$ meet. In view of Corollary 2.5 there must exist finite sets $A \subseteq H$ and $B \subseteq H^*$ such that $\operatorname{conv}(f(A))$ meets $\operatorname{conv}(f(B))$. Let $N$ denote the subalgebra of $M$ generated by $A \cup B$; from Corollary 7.5 we know that $N$ is finite. Here $ev\colon N \to N^{\circ\circ}$ is an isomorphism, hence by the commutative diagram of Lemma 7.1 the double dual $(f|_N)^{\circ\circ}\colon N^{\circ\circ} \to M_2^{\circ\circ}$ is also injective. From (ii) above we deduce that $(f|_N)^\circ\colon M_2^\circ \to N^\circ$ is surjective.

In particular there exists a half space $J \in M_2^\circ$ with $f^\circ(J) \cap N = H \cap N$, thus $J$ separates $f(A)$ from $f(B)$. But then

$$\mathrm{conv}(f(A)) \cap \mathrm{conv}(f(B)) = \varnothing,$$

contradicting the assumption that $\mathrm{conv}(f(H))$ and $\mathrm{conv}(f(H^*))$ meet. We conclude that $f^\circ$ is surjective.

(iv) We already know that $f^\circ$ is a poc morphism, so to prove that $f^\circ$ is an embedding it remains to show that

$$\forall H_1, H_2 \in M_2^\circ: \quad H_1 \nsubseteq H_2 \;\Rightarrow\; f^\circ(H_1) \nsubseteq f^\circ(H_2).$$

This is equivalent to

$$\forall H_1, H_2 \in M_2^\circ: \quad H_1 \cap H_2^* \neq \varnothing \;\Rightarrow\; f^{-1}(H_1) \cap f^{-1}(H_2)^* = f^{-1}(H_1 \cap H_2^*) \neq \varnothing,$$

which is true if $f$ is surjective. □

I don't know if the converse of part (iv) is true in general, but I can show it under additional hypothesis that $M_2$ is discrete. Recall that in a discrete median algebra all separators $\Delta(x,y)$ are finite. First we need two short lemmas.

### 7.9. Lemma.

Let $M$ be a median algebra and $x, y \in M$. Then $\bigcap \Delta(x,y) = \{z \in M \mid x \in [y,z]\}$.

*Proof.*
$$\begin{aligned}
x \in [y,z] &\Leftrightarrow \Delta(x,\{y,z\}) = \varnothing \\
&\Leftrightarrow \forall H \in \Delta(x,y): z \in H \\
&\Leftrightarrow z \in \bigcap \Delta(x,y)
\end{aligned}$$
□

### 7.10. Lemma.

Let $N$ be a median subalgebra of $M$ and let $C_1$, ..., $C_n$ be convex subsets of $M$. If $C_i \cap C_j \cap N \neq \varnothing$ for all $i$ and $j$, then $C_1 \cap \ldots \cap C_n \cap N \neq \varnothing$.

*Proof.* This lemma is a slight strengthening of Helly's Theorem 2.2, and indeed the same proof works here again. For $n = 2$ nothing is to prove. For $n = 3$ we choose elements $x_i \in C_j \cap C_k \cap N$, where $\{i,j,k\} = \{1,2,3\}$. Then $m(x_1,x_2,x_3)$ lies in $C_1 \cap C_2 \cap C_3$, but since $N$ is a median algebra, it also lies in $N$. The general case follows by induction. □

### 7.11. Proposition.

Let $f\colon M_1 \to M_2$ be a median morphism, where $M_2$ is discrete. If $f^\circ\colon M_2^\circ \to M_1^\circ$ is an embedding of poc sets, then $f$ is surjective.

*Proof.* Denote the image of $f$ by $N$ and assume that there exists an element $v \in M_2 \smallsetminus N$. Choose an arbitrary element $x \in N$. As $M_2$ is discrete, the set $\Delta(v,x)$ is finite. Furthermore, for any $H_1, H_2 \in \Delta(v,x)$, the intersection $H_1 \cap H_2$ is nonempty, because it contains $v$. As $f^\circ$ is an embedding, it follows that $H_1 \cap H_2 \cap N \neq \emptyset$. We conclude from Lemma 7.10 that $N \cap \bigcap \Delta(v,x)$ is nonempty.

Choose some element $y \in N \cap \bigcap \Delta(v,x)$. Again, the set $\mathscr{H} := \Delta(v,x) \cup \Delta(v,y)$ is finite and $N \cap \bigcap \mathscr{H} \neq \emptyset$ contains at least one element, say $z$. Now Lemma 7.9 says that

$$v \in [x,y] \cap [x,z] \cap [y,z].$$

As $M_2$ is a median algebra, this means that $v = m(x,y,z)$, and since $x$, $y$ and $z$ lie in the subalgebra $N$, we conclude that $v \in N$. This contradicts our assumption that $f$ is not surjective. $\square$

If $X$ is a subset of a median algebra $M$, and $j\colon X \hookrightarrow M$ is the inclusion, then the map $j^\#\colon M^\circ \to \mathscr{P}X$ is simply given by $H \mapsto H \cap X$. In this case we write $\operatorname{res}_X$ instead of $j^\#$. To say that $\operatorname{res}_X$ is an embedding means that

$$\forall H_1, H_2 \in M_2^\circ: \quad H_1 \cap H_2 \neq \emptyset \quad \Rightarrow \quad H_1 \cap H_2 \cap X \neq \emptyset.$$

**7.12. Proposition.**
Let $M$ be a discrete median algebra and $X \subseteq M$. The following are equivalent.
  (i) $M$ is generated by $X$.
  (ii) The map $\operatorname{res}_X\colon M^\circ \to \mathscr{P}X$ is an embedding.

*Proof.* Dualizing the diagram for the universal property of $\operatorname{Med}(X)$, we obtain another commutative diagram.

$$\begin{array}{ccc}
\operatorname{Med}(X)^\circ & \xleftarrow{f^\circ} & M^\circ \\
& \searrow_{\iota^\#} \quad \swarrow_{\operatorname{res}_X} & \\
& \mathscr{P}X &
\end{array}$$

Now condition (i) is equivalent to the surjectivity of $f\colon \operatorname{Med}(M) \to M$. By Proposition 7.8 and Proposition 7.11 this equivalent to the statement that $f^\circ\colon M^\circ \to \operatorname{Med}(X)^\circ$ is an embedding. But the arrow labeled $\iota^\#$ is an isomorphism, hence this equivalent to condition (ii). $\square$

Finally, let us state and prove the result that was left over at the end of the last section.

**7.13. Proposition.**
Let $M_1$ be a median algebra, $P_1 := M_1^\circ$ and let $U \subset P_1$ be a subset satisfying (**Con 1-3**). Let $P_2 := P_1 \smallsetminus U$ and $M_2 := M_1/\nabla(U)$, then $M_2^\circ$ is isomorphic to $P_2$.

*Proof.* Let $q\colon M_1 \to M_2$ denote the canonical quotient map. Proposition 7.8 says that the map $q^\circ\colon M_2^\circ \to M_1^\circ$, $H \mapsto q^{-1}(H)$ is an embedding of poc sets; here we will show that its image is $P_2$.

As $q$ is a surjective median morphism it maps convex sets onto convex sets. If $H \in P_2$ then $q(H)$ and $q(H^*)$ are disjoint, hence they are complementary half spaces of $M_2$. Therefore $H = q^{-1}(q(H))$, which means that $P_2 \subseteq q^\circ(M_2^\circ)$.

To prove the converse inclusion pick any proper $J \in M_2^\circ$, then $q^{-1}(J)$ is a half space of $M_1$, and for any $x, y \in M_1$ with $q^{-1}(J) \in \Delta(x, y)$ we certainly have $q(x) \neq q(y)$, which means that $q^{-1}(J) \notin \triangle \circ \triangledown(U) = U$, hence $q^{-1}(J) \in P_1 \smallsetminus U = P_2$. □

**Notes.**
A filter base in the poc set $\mathscr{P}X$ is also called a linked system, because it is a family of subsets of $X$ that meet pairwise. An ultra filter is therefore a maximal linked system. The set $\lambda(r)$ of maximal linked systems in $X := \{1, \ldots, r\}$ has been studied in various contexts, see VAN DE VEL [1983] and BANDELT-VAN DE VEL [1991]. In particular, the finite case of Theorem 7.4 seems to be folklore.

# §8. Discrete Median Algebras

We call two elements $x, y \in M$ almost equal, if $\Delta(x, y)$ is finite or, equivalently, the interval $[x, y]$ is finite. This relation is clearly reflexive and symmetric, it is transitive by the triangle equality, Lemma 2.16, and a congruence relation by Lemma 6.19. A median algebra is discrete in the sense of § 2 if it comprises a single almost equality class. In this section we provide some characterizations of discrete median algebras and discuss a few examples.

Our first lemma allows us to recognize the principal ultra filters in $M^{\circ\circ}$.

**8.1. Lemma.** *An element $\zeta \in M^{\circ\circ}$ lies in $ev(M)$ iff $\bigcap \zeta \neq \varnothing$.*

*Proof.* If $\zeta = ev(z)$, then $z \in \bigcap \zeta$. Conversely, suppose that $\bigcap \zeta$ is nonempty, then it must be a singleton, because for any $x \neq y$ in $M$ there exists a half space $H \in \Delta(x, y)$, and $\zeta$ contains either $H$ or $H^*$, hence $\bigcap \zeta$ contains either $x$ or $y$, but not both. Thus there exists a $z \in M$ with $\bigcap \zeta = \{z\}$, which means that $\zeta = ev(z)$. □

**8.2. Theorem.**
*For a median algebra $M$ the following are equivalent.*
  (i) *$M$ is discrete.*
  (ii) *Every non empty convex subset of $M$ is a retract.*
  (iii) *$ev(M)$ is an almost equality class in $M^{\circ\circ}$.*
  (iv) *$ev(M)$ is a convex subset of $M^{\circ\circ}$.*

*Proof.* (i) ⇒ (ii). Suppose that $M$ is discrete and $C \subseteq M$ is a non empty convex subset. For any $x \in M$ there exists a $y \in C$ such that $d(x, y)$ is minimal. We claim that $y \in C$ is in fact closest to $x$, i.e., for every $z \in C$ we have $y \in [z, x]$. To show this consider the element $m = m(x, y, z)$. Obviously $m$ lies in $[x, y]$, so $d(x, m) \leq d(x, y)$, where equality holds iff $m = y$. But $m$ lies in $C$ as $y, z \in C$, so by the choice of $y$ we have $y = m$, which means that $y \in [x, z]$ as required.

(ii) ⇒ (i). This is the content of Proposition 5.15.

(i) ⇒ (iii). Consider an element $\zeta \in M^{\circ\circ}$ as an ultrafilter of half spaces in $M$. For $x, y \in M$ we have

$$\zeta \in [ev(x), ev(y)] \quad \Leftrightarrow \quad ev(x) \cap ev(y) \subseteq \zeta \subseteq ev(x) \cup ev(y).$$

The difference between the two bounds is

$$(ev(x) \cup ev(y)) \smallsetminus (ev(x) \cap ev(y)) = ev(x) + ev(y) = \Delta(x, y) \cup \Delta(y, x).$$

Thus if $[x, y]$ is finite then also $[ev(x), ev(y)]$ is finite, hence if $M$ is discrete then $ev(M)$ is contained in an almost equality class of $M^{\circ\circ}$.

Next we show that if $\zeta \in M^{\circ\circ}$ is almost equal to some $ev(x)$, then $\zeta$ also lies in $ev(M)$. It is enough to consider the case when $d(ev(x), \zeta) = 1$; the general case follows from this special case by induction. Note that for any $H \in \zeta \smallsetminus ev(x)$ the half space $ev(H) \in M^{\circ\circ\circ}$ separates $\zeta$ from $ev(x)$. Therefore $\zeta$ and $ev(x)$ differ by a single half space in $M$, say $\zeta \smallsetminus ev(x) = \{H\}$.

Moreover, any half space $J \in \Delta(H, x)$ lies in $\zeta$, because $H \subseteq J$, but not in $ev(x)$, because $x \notin J$. It follows that $\Delta(H, x) = \zeta \smallsetminus ev(x) = \{H\}$. By the implication (i) $\Rightarrow$ (ii) we know that $H$ is a retract. In view of Lemma 2.11 this means that there exists a $z \in H$ with $\Delta(H, x) = \Delta(z, x)$, in other words $ev(z) = (ev(x) \smallsetminus \{H^*\}) \cup \{H\} = \zeta$. That completes the second half of this proof.

(iii) $\Rightarrow$ (iv). This is an immediate consequence of the fact that almost equality is a congruence relation.

(iv) $\Rightarrow$ (i). We will show that every interval $[x, y] \subseteq M$ is compact in the convex topology. Then Corollary 6.18 implies that all intervals are finite, hence $M$ must be discrete. In order to prove that an interval is compact in the convex topology we have to show that its convex subsets satisfy (**SFIP**). Let $\mathscr{C}$ be any non empty family of convex subsets of $[x, y]$, such that for any two $C_1, C_2 \in \mathscr{C}$ we have $C_1 \cap C_2 \neq \varnothing$. We have to show that $\mathscr{C}$ has nonempty intersection. Let $Z := \{H \in M^\circ \mid H \text{ contains some } C \in \mathscr{C}\}$. The set $Z$ is a filter in $M^\circ$, because the elements of $\mathscr{C}$ meet pairwise. By the extension algorithm there exists an ultra filter $\zeta \in M^{\circ\circ}$ containing $Z$.

We claim that $\zeta \in [ev(x), ev(y)]$. If any half space $H \in M^\circ$ contains both $x$ and $y$, then it contains every element of $\mathscr{C}$, so it must lie in $Z$; if $H$ contains neither $x$ nor $y$ then $H^*$ is contained in $Z$, but then $H$ can't lie in $\zeta$. Thus we have $ev(x) \cap ev(y) \subseteq \zeta \subseteq ev(x) \cup ev(y)$, which proves our claim.

Now condition (iv) says that there exists a $z \in M$ with $\zeta = ev(z)$. It remains to show that $z \in \bigcap \mathscr{C}$. Suppose that there exists a $C \in \mathscr{C}$ with $z \notin C$. Then there exists a half space $H \in \Delta(C, z)$, so $z \notin H \in Z$, contrary to the choice of $z$. $\square$

### 8.3. Examples.

(i) The median algebra $\mathscr{F}X$ of finite subsets of an infinite set $X$ is clearly discrete. Thus all convex subsets are retracts. This means that the biretract topology on $\mathscr{F}X$ agrees with the convex topology. But the convex topology is not quasi compact, as $\mathscr{F}X^\circ$ is an infinite orthogonal poc set. By Remark 5.7, $\mathscr{F}X$ is an example of a median algebra that is not the dual of any poc set.

(ii) In Remark 5.7 we announced the existence of two non isomorphic discrete median algebras with isomorphic duals. These we will construct now.

Let **3** denote the set $\{0, 1, 2\}$ with the obvious linear median structure, i.e., $m(0, 1, 2) = 1$. Let $M$ be the set of functions $f \colon \mathbb{Z} \to \mathbf{3}$ with the pointwise median operation. For each $n \in \mathbb{Z}$ we have a canonical projection $\pi_n \colon M \to \mathbf{3}$, $f \mapsto f(n)$. For $i \in \mathbf{3}$ we have the constant functions $c_i \colon \mathbb{Z} \to \mathbf{3}$ with value $i$. Let $M_i$ denote the almost equality class of $c_i$, i.e., the median algebra of functions $f\colon \mathbb{Z} \to \mathbf{3}$ such that $f^{-1}(i)$ is cofinite. Clearly, $M_i$ is a discrete median algebra.

Now we show that the duals $M_i^\circ$ are all the same. For any proper half space $H \subset M_i$ there exists an $n \in \mathbb{Z}$ and a half space $h \subset \mathbf{3}$ such that $H = M_i \cap \pi_n^{-1}(h)$. This follows from the fact that there exist functions $f_1, f_2 \in M_i$ such that $\{H\} = \Delta(f_1, f_2)$. In particular, $f_1$ and $f_2$ differ only in a single place $n \in \mathbb{Z}$, and $h$ is the separator of $f_1(n)$ and $f_2(n)$ in $\mathbf{3}$. It follows that
$$M_i^\circ \cong \bigoplus_{\mathbb{Z}} \mathbf{3}^\circ.$$

Next we show that $M_1$ and $M_0$ are not isomorphic as median algebras. Consider the ultra filter $U = ev(c_1) \in M_1^{\circ\circ}$. It can be written as
$$U = \{\pi_n^{-1}\{0,1\} \cap M_1,\ \pi_1^{-1}\{1,2\} \cap M_1 \mid n \in \mathbb{Z}\} \cup \{M_1\},$$
and it has the unique property that for every proper $H \in U$ there exists a $J \in U$ with $H^* \subsetneq J$. Every automorphism $\alpha$ of $M_1$ preserves the inclusion order of the half spaces, hence $U = \{\alpha H \mid H \in U\}$, and because $\bigcap U = \{c_1\}$ the constant function $c_1$ is fixed by every automorphism of $M_1$.

On the other hand we can construct a fixed point free automorphism $\alpha\colon M_0 \to M_0$. Let $\rho\colon \mathbf{3} \to \mathbf{3}$ denote the exchange of 0 and 2. Let
$$(\alpha f)(n) := \begin{cases} \rho(f(0)) & \text{if } n = 1, \\ f(n-1) & \text{if } n \neq 1. \end{cases}$$
Then
$$\alpha^k(c_0)(n) = \begin{cases} 0 & \text{for } n \leq 0 \text{ or } n \geq k, \\ 2 & \text{for } 1 \leq n \leq k. \end{cases}$$
It follows that the orbit of $c_0$ under the action of $\alpha$ is unbounded, hence $\alpha$ can't have a fixed point on $M_0$. (See Example 8.11 for another argument.) □

If we take the Stone topology on $M^\circ$ into account we get a particularly fitting characterization of discrete median algebras.

### 8.4. Theorem.
*A median algebra $M$ is discrete iff the set of proper half spaces, $M^\circ \smallsetminus \{\varnothing, M\}$, is a discrete topological subspace of $M^\circ$.*

*Proof.* If $M$ is discrete, then for every proper half space $H \in M^\circ$ there exist $x, y \in M$ with $\{H\} = \Delta(x, y)$, in particular $\{H\}$ is an open set, thus $M^\circ \smallsetminus \{\varnothing, M\}$ is discrete. Conversely, assume that $M^\circ \smallsetminus \{\varnothing, M\}$ is discrete. Then for each proper half space $H \in M^\circ$ the set $\{H\}$ is open. In particular, $\{H\}$ is a basic open set, so $\{H\} = \Delta(X, Y)$ for some finite sets $X, Y \subseteq M$. By Lemma 2.15 we can find a gate $(x, y)$ for the retracts $\operatorname{conv}(X)$ and $\operatorname{conv}(Y)$, and then by Lemma 2.14, $\{H\} = \Delta(x, y)$ for some $x, y \in M$.

Now let $R$ denote the congruence relation of almost equality. We compute the set $\triangle(R)$,

that we defined in §6.

$$M^\circ \smallsetminus \{\varnothing, M\} \supseteq \triangle(R) = \bigcup \{\triangle(x,y) \mid (x,y) \in R\}$$
$$\supseteq \{H \mid \{H\} = \triangle(x,y) \text{ for } x,y \in M\}$$
$$= M^\circ \smallsetminus \{\varnothing, M\}.$$

By Proposition 6.22 we can recover $R$ as $\triangledown \circ \triangle(R)$, but that is the trivial congruence relation with a single equivalence class, which means that $M$ is discrete. □

BOUNDARIES, STARS AND AN EXAMPLE

In a discrete median algebra each proper half space $H$ has a **boundary** which is defined as follows.

$$\partial H := \{x \in H \mid \text{ there exists } x' \in H^* \text{ with } d(x, x') = 1\}.$$

The next lemma gives a different expression for $\partial H$, which immediately implies that the sets $\partial H$, $\partial H \cup H^*$ and $\partial H \cup \partial H^*$ are all convex. Note that $H \smallsetminus \partial H$ is not convex in general.

**8.5. Lemma.**
$\partial H = H \cap \bigcap \{J \in M^\circ \mid H^* \subset J\}$.

*Proof.* Choose $x \in \partial H$ and a half space $J$ with $H^* \subset J$; we want to show that $x \in J$. But there exists an $x' \in H^* \subset J$ such that $J$ does not separate $x$ from $x'$, thus $x \in J$. Conversely, suppose $x \in H \cap \bigcap \{J \in M^\circ \mid H \subset J\}$ and choose an element $y \in H^*$ with $d(x,y)$ minimal. Assume there exists a half space $J \in \Delta(x,y)$ with $J \neq H$. We can't have $J \subset H$, because then $H^* \subset J^*$, so $x$ would have to lie in $J^*$. Thus we can choose an element $z \in H^* \cap J$, and we let $y' := m(x,y,z)$. Now we have

(i) $y' \in [x,y]$,

(ii) $y' \in J$, because $x, z \in J$, and

(iii) $y' \in H^*$, because $y, z \in H^*$.

Hence $d(x,y') < d(x,y)$, contradicting the choice of $y$. Thus no such $J$ exists, therefore $d(x,y) = 1$ and $x \in \partial H$. □

**8.6. Remark.**
The point $x'$ above is uniquely determined by $x$. For suppose $x', x'' \in H^*$ with $d(x,x') = d(x,x'') = 1$, then $\{m(x,x',x'')\} = [x,x'] \cap [x,x''] = \{x\}$, thus $x \in [x',x''] \subseteq H^*$, but that is impossible. Moreover, it is easy to show that the map $x \mapsto x'$ defines a median isomorphism $\partial H \to \partial H^*$. □

**8.7. Lemma.**
If $H_1, \ldots, H_n$ are pairwise transverse, then $\partial H_1 \cap \ldots \cap \partial H_n$ is nonempty and equals the set
$$\{x \in H_1 \cap \ldots \cap H_n \mid \text{there exists } y \in H_1^* \cap \ldots \cap H_n^* \text{ with } d(x,y) = n\}.$$

*Proof.* By transversality the sets $H_1 \cap \ldots \cap H_n$ and $H_1^* \cap \ldots \cap H_n^*$ are both nonempty. We can choose points $x$ and $y$ in either set with minimal distance and argue as in the previous lemma that $\Delta(x,y) = \{H_1, \ldots, H_n\}$. □

We have frequently encountered finite cubes, which are the median algebras $\mathbf{2}^n$. The infinite powers of $\mathbf{2}$ are no longer discrete, so we prefer the follwing definition, which agrees with the previous one for finite median algebras. A median cube is a discrete median algebra whose dual is orthogonal. A subcube of a median algebra $M$ is a convex subset $C \subseteq M$ which is a cube.
We collect some obvious facts about cubes.

**8.8. Lemma.**
  (i) If $C$ is a cube with $C^\circ \cong \mathrm{Poc}(X)$, then $C$ is isomorphic to the median algebra $\mathscr{F}X$.
  (ii) $C \subseteq M$ is a cube iff $\{\overline{H} \in \overline{M^\circ} \mid H \text{ cuts } C\}$ is transverse.
  (iii) A subset of a cube is convex iff it is a cube. □

A median algebra $M$ is called a star with centre $c$ if for all $x \in M$ the interval $[x,c]$ is a finite cube. Clearly, stars are always discrete. For example the starlets of Example 5.13(i) are precisely the 1-dimensional stars.

An interesting feature of stars is that their median structure is coded by the transverality graph of their dual. Let $M$ be a star and $\mathscr{H}$ the set of hyper planes. Let $\mathscr{T}$ be the family of finite, transverse subsets. This is a median subalgebra of $\mathscr{F}\mathscr{H}$ with the Boolean median structure.

**8.9. Proposition.**
Let $M$ be a star with center $c$. The map $\tau: M \to \mathscr{T}$, $x \to \overline{\Delta}(x,c)$ defines an isomorphism of median algebras.

*Proof.* If $\overline{\Delta}(x,c) \subseteq \overline{\Delta}(y,c)$ then $x \in [y,c]$; it follows that $\tau$ is injective. Suppose that $\{\overline{H}_1, \ldots, \overline{H}_n\}$ is a finite transverse set, and without loss of generality $c \in H_i^*$ for all $i$, then by Lemma 8.7 there exists a point $x \in H_1 \cap \ldots \cap H_n$ with distance $d(x,c) = n$, so $\Delta(x,c) = \{H_1, \ldots, H_n\}$, which means that $\tau$ is surjective. The fact that $\tau$ is a median morphism follows from the next lemma. □

**8.10. Lemma.**
If $c, x, y, z \in M$, then $\Delta(m(x,y,z), c) = m(\Delta(x,c), \Delta(y,c), \Delta(z,c))$.

*Proof.*  $H \in \Delta(m(x,y,z), c) \quad \Leftrightarrow \quad m(x,y,z) \in H \text{ and } c \notin H$

$\quad \Leftrightarrow \quad$ (at least two out of $x$, $y$, $z$ lie in $H$) and $c \notin H$

⇔   $H$ lies in at least two out of $\Delta(x,c), \Delta(y,c), \Delta(z,c)$

⇔   $H \in m(\Delta(x,c), \Delta(y,c), \Delta(z,c))$. □

**8.11. Example.**
Recall the two median algebras $M_0$ and $M_1$ from Example 8.3(ii). It is easy to see that $M_1$ is a star with center $c_1$, whereas for every $c \in M_0$ there exists an $x \in M_0$ with $[x,c] \cong \mathbf{3}$, hence $M_0$ is not a star. This gives another proof of the fact that the two median algebras are not isomorphic. □

Let $M$ be any discrete median algebra. We define the star at a vertex $v \in M$ as
$$\mathrm{star}(v, M) := \{x \in M \mid [x, v] \text{ is a cube}\}.$$
We will show that $\mathrm{star}(v, M)$ is a convex subset of $M$. This means in particular that it is a median subalgebra, so $\mathrm{star}(v, M)$ is indeed a star in the above sense.

**8.12. Proposition.**
If $M$ is a median algebra and $v \in M$, then $\mathrm{star}(v, M)$ is a convex subset of $M$.

*Proof.* Suppose $x, z \in \mathrm{star}(v, M)$ and $y \in [x, z]$. We need to show that any two distinct $H, J \in \Delta(y, v)$ are transverse. The previous lemma says that $\Delta(y, v)$ lies in the interval between $\Delta(x, v)$ and $\Delta(z, v)$, which means that
$$\Delta(x,v) \cap \Delta(z,v) \subseteq \Delta(y,v) \subseteq \Delta(x,v) \cup \Delta(z,v).$$
If $H$ and $J$ both lie in $\Delta(x,v)$ then they are transverse by the assumption $x \in \mathrm{star}(v, M)$, and similarly if they both lie in $\Delta(z,v)$. The only remaining possibility is that, say, $H \in \Delta(x, \{v, z\})$ and $J \in \Delta(z, \{v, x\})$. Now we have
$$x \in H \cap J^*, \qquad z \in H^* \cap J, \qquad y \in H \cap J, \qquad v \in H^* \cap J^*,$$
in other words, $H$ and $J$ are transverse. □

For a general median algebra $M$ we still can define bijection $\tau \colon M \to \mathscr{T}$, once we have chosen a base point, but it will no longer be a median morphism.

Let us first explain this for a simple example. Let $T$ be a finite tree, then we have $|VT| = |ET| + 1$. A constructive proof of this fact goes like this: Choose any vertex $v \in VT$ and orient all edges of $T$ towards $v$, then we get a bijection by mapping each edge of $T$ on its initial vertex. This bijection $ET \to VT \smallsetminus \{v\}$ even works for for infinite trees, and the following proposition is the appropriate generalization to discrete median algebras.

**8.13. Proposition.**
Let $M$ be a discrete median algebra, choose an element $v \in M$ and let $P := M^\circ$. Let $\mathscr{H} := \{\overline{H} \in \overline{P} \mid H \text{ proper}\}$ denote the set of hyper planes in $M$, and let $\mathscr{T}$ denote the *family of finite transverse subsets of $\mathcal{H}$ (including $\varnothing$ and all singletons). Then we have a bijection of sets,*

$$\tau \colon M \to \mathcal{T}, \quad x \mapsto \overline{\mathrm{Min}(\Delta(x,v))}.$$

*Sketch of Proof.* We omit the (routine) proof and just note that the inverse of $\tau$ can be described in the following way. Let $F \subseteq \mathcal{T}$ be a finite transverse subset, then the set $f := \{a \in P \mid \bar{a} \in F \text{ and } v \in a^*\}$ is a filter base with $\bar{f} = F$, thus

$$U := (ev(v) \smallsetminus (\uparrow\! f)^*) \cup \uparrow\! f$$

is an ultra filter. Moreover, since $U \dotplus ev(v)$ is finite, there exists an $x \in M$ with $U = ev(x)$, and $\tau(x) = F$. $\square$

**8.14. Example.**
The previous proposition is most useful if we know a median algebra by its dual, and if it has a canonical vertex. This is the case for the free median algebra $M_5$ with five generators. Recall from Theorem 7.4 that $M_5^\circ \cong \mathscr{P}\mathbf{5}$, where $\mathbf{5} := \{0,1,2,3,4\}$. Observe further, that $\mathscr{P}\mathbf{5}$ has a canonical ultra filter $U := \{a \subseteq \mathbf{5} \mid |a| \geq 3\}$.

In the following we give a brief description of the median structure of $M_5$. The set $\mathcal{H}$ of hyper planes will be identified with the set $\mathcal{H}' := \{a \subset \mathbf{5} \mid 1 \leq |a| \leq 2\}$. A finite subset $f \subseteq \mathcal{H}'$ corresponds to a transverse subset of $\mathcal{H}$ iff

$$\forall a, b \in f : \quad a \not\subseteq b \ \wedge \ b \not\subseteq a \ \wedge \ a \cap b \neq \varnothing,$$

the fourth condition $a^* \cap b^* \neq \varnothing$ is automatically satisfied. It is somewhat more economical to list the transverse sets in $\mathcal{H}'$ rather than the maximally linked systems in $\mathscr{P}\mathbf{5}$.

A simple case discussion yields six types of transverse subsets, which together make up the 81 elements of $M_5$. Suppose that $\mathbf{5} = \{x,y,z,s,t\}$.

| Shorthand | $f$ | Number |
|---|---|---|
| $\varnothing$ | $\varnothing$ | 1 |
| $(\ \vert x)$ | $\{x\}$ | 5 |
| $(xy)$ | $\{x,y\}$ | $\binom{5}{2} = 10$ |
| $(x\vert yz)$ | $\{x,y\}, \{x,z\}$ | $5\binom{4}{2} = 30$ |
| $(x\vert t)$ | $\{x,y\}, \{x,z\}, \{x,s\}$ | $5\binom{4}{1} = 20$ |
| $(xyz)$ | $\{x,y\}, \{y,z\}, \{z,x\}$ | $\binom{5}{3} = 10$ |
| $(x\vert\ )$ | $\{x,y\}, \{x,z\}, \{x,s\}, \{x,t\}$ | 5 |

The map $\tau$ sends $\varnothing$ to the canonical vertex $U$, which turns out to be the central vertex of $M_5$, and the principal ultra filter $\uparrow\! x$ to $(\ \vert x)$, which is furthest away from the centre.

With a little more work one can determine the intervals $M^x := [\varnothing, (\ |x)]$, which consists of a four-cube plus a hair, and $M^{xyz} := [\varnothing, (xyz)]$, which is a three-cube.

Now $M_5$ is the union of the $M^x$ and the $M^{xyz}$, where $M^x$ and $M^y$ meet in an edge and $M^x$ and $M^{xyz}$ meet in a square. In other words, $M_5$ consist of a star with center $\varnothing$, together with 5 hairs. □

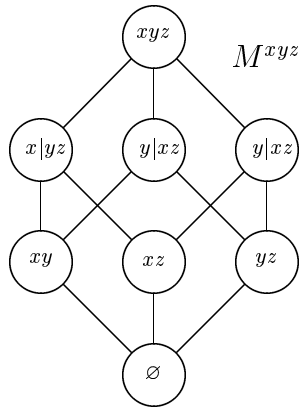

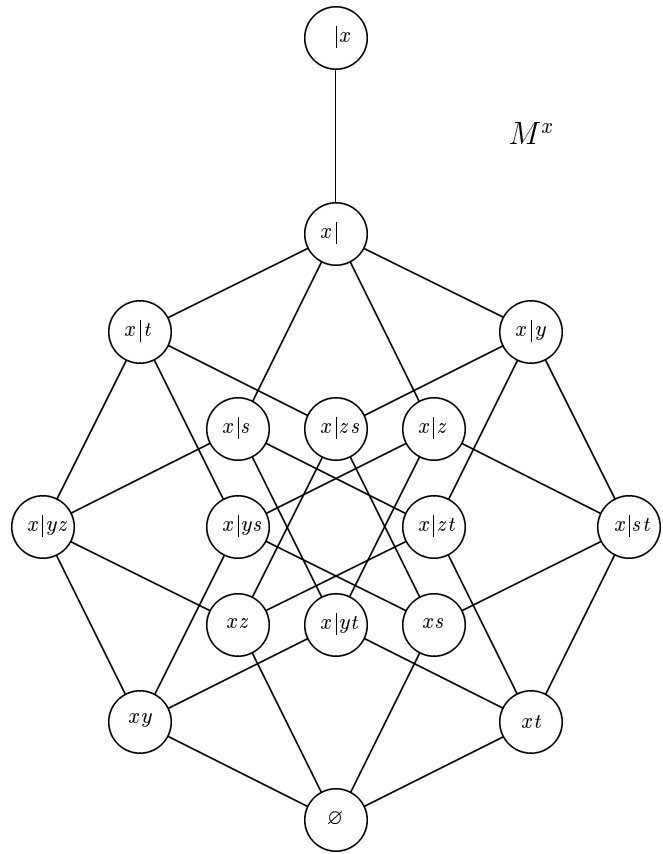

# §9. Discrete Poc Sets

So far we have studied median algebras for their own sake, but our ultimate motive is to use them as an interesting substratum for group actions. How does one construct median algebras with prescribed properties, e.g. discreteness and specified group actions? The difficulty usually lies in the verification of the median axiom (**Med 3**).

Our duality theory suggests the following strategy. We first construct the poc set $M^\circ$, whose structure is much easier to control. Then we try to find the image of the double dual map $M \to M^{\circ\circ}$; e.g. if $M$ is discrete then $ev(M)$ is an almost equality class of $M^{\circ\circ}$. But now a new difficulty arises:

> How can we recognize poc sets which are duals of discrete median algebras?

The prototypical answer to this question is the construction of DUNWOODY [1979]. He studied the partial ordering and orientation reversion on the oriented edge set of a tree, that we introduced in our Model Example 1.6. As the characteristic properties he identified nesting and interval finiteness.

We say that a poc set $P$ is discrete if for all proper elements $a, b \in P$ the interval $[a, b] := \{c \in P \mid a \leq c \leq b\}$ is finite. Dunwoody's result can now be expressed as follows.

**9.1. Theorem.** DUNWOODY [1979]
*If $P$ is a discrete nested poc set, then there exists a discrete median algebra $M$ with $P \cong M^\circ$.* □

In this section we push Dunwoody's construction forward in two directions. First we modify it to work for discrete poc sets of dimension $\omega$; then we give a construction for countable poc sets. While in the first case the median algebra is canonical, the second construction depends on an explicit enumeration, thus it may not be very useful for the study of group actions.

There is, however, a definite limit to these extensions of Dunwoody's construction. We will give an example of a discrete poc set that is not the dual of any discrete median algebra.

Finally we will study discrete representations of poc sets. Here we can use additional information that comes with poc set constructions typical for group theoretic applications.

A NECESSARY AND SUFFICIENT CONDITION

The first proposition shows that the condition of discreteness is certainly necessary for our purpose, but not sufficient. Nor can we hope to find many discrete median algebras as duals of discrete poc sets.

### 9.2. Proposition.
  (i) If $M$ is discrete, then $M°$ is a discrete poc set; the converse is false.
  (ii) $P°$ is discrete iff $P$ is of type $\omega$.

*Proof.* (i) If $H_1$ and $H_2$ are proper half spaces in $M$ with $H_1 \subseteq H_2$, then we can choose points $x \in H_1$, $y \in H_2^*$. Any half space $H$ with $H_1 \subseteq H \subseteq H_2$ also separates $x$ and $y$, so $[H_1, H_2] \subseteq \Delta(x, y)$ is finite.

The dual of the Boolean median algebra $\mathscr{P}X$ is orthogonal, hence always discrete, but $\mathscr{P}X$ is discrete only when $X$ is finite.

(ii) Suppose that $P$ is not of type $\omega$, then it has an infinite subset $X$ that is either a chain or transverse. By the Extension Algorithm 3.3 there exist ultra filters $U \supseteq X$ and $V \supseteq X^*$ in $P$, and the interval $[U, V]$ is infinite. Conversely, if $P$ is of type $\omega$, then intervals in $P°$ are finite by Corollary 6.17. □

In the previous section we have shown that $M$ is discrete iff $M° \smallsetminus \{0, 0^*\}$ is a discrete topological space. Thus all that we need to find is a system of neighbourhoods of $0^*$ that turns $M°$ into a Stone poc set. In fact, it is enough to specify a single clopen ultra filter. We call an ultra filter $U \subset P$ founded, if for every $a \in U$ there exist a minimal element $b \in U$ with $b \leq a$; in other words, $U$ is generated as a filter by the set $\mathrm{Min}(U)$. For example, the principal ultra filters in $M°$ are all founded; this is just a reformulation of Theorem 2.7. An ultra filter $U \subset P$ is called well founded, if for every proper element $a \in U$ the set $\{b \in U \mid b \leq a\}$ is finite.

### 9.3. Proposition.
*A poc set $P$ has a well founded ultra filter if and only if there exists a discrete median algebra $M$ with $P \cong M°$.*

*More precisely, if $U \in P°$ is well founded, then $M$ may be chosen as the almost equality class of $U$ in $P°$.*

*Proof.* If $M$ is discrete and $x \in M$, then $ev(x) \in M°°$ is well founded. We have $H \in ev(x)$ iff $x \in H$; as $H$ is proper there exists $y \in H^*$. Now $\{J \in ev(x) \mid J \subseteq H\} \subseteq \Delta(x, y)$, which is finite.

Conversely, let $U$ be a well founded ultra filter and put
$$\mathscr{S} := \{V, P \smallsetminus V \mid V \in P° \text{ and } U + V \text{ finite}\}.$$
We show that $\mathscr{S}$ is the subbasis of a topology that turns $P$ into a Stone poc set.

The condition (**Stone 1**) is clear. To prove (**Stone 2**) we use the well foundedness of $U$ in the following way. For every $a \in P$ the set $I_a := \{b \in P \mid b \leq a\}$ is finite, thus $V_a := (U \smallsetminus I_a) \cup I_a^*$ is an ultra filter which also lies in $\mathscr{S}$. Suppose that $a, b \in P$ are proper elements with $a \not\leq b^*$, we need to find an ultra filter $V \in \mathscr{S}$ with $a, b \in V$. We have to discuss four cases.

CASE 1: $a, b \in U$, then we are done.

CASE 2: $a, b^* \in U$, then $b \in V_{b^*}$ and $a \in V_{b^*}$ as $a \not\leq b^*$.

CASE 3: $a^*, b \in U$, then similarly $a, b \in V_{a^*}$.

CASE 4: $a^*, b^* \in U$. If $a < b$ or $a$ and $b$ are incomparable, then $a, b^* \in V_{a^*}$, continue as in CASE 2. If $a > b$ then $a^*, b \in V_{b^*}$, continue as in CASE 3.

Since we found a clopen ultra filter that separates $a$ from $b^*$ we have also verified condition (**SPS**).

To check condition (**Stone 3**) let $\mathscr{V} \subseteq \mathscr{S}$ be a centered subfamily. If all elements of $\mathscr{V}$ are filters, then $0^* \in \bigcap \mathscr{V}$, and if they are all ideals, then $0 \in \bigcap \mathscr{V}$. Otherwise $\mathscr{V}$ contains both a filter $V_1$ and an ideal $P \smallsetminus V_2$. Then $V_1 \smallsetminus V_2$ is finite, and since $\mathscr{V}$ is centered the intersection $\bigcap \mathscr{V}$ must be nonempty.

By the same argument every proper element of $P$ is contained in a finite open subset, and since the topology on $P$ is Hausdorff, the induced topology on $P \smallsetminus \{0, 0^*\}$ is discrete. Now we have furnished $P$ with a topology that turns it into a Stone poc set, such that $M := P^\bullet$ is a discrete median algebra with $M^\circ \cong P$. More precisely, the proof of Theorem 5.3 says that $M$ is the set $\{V \in P^\circ \mid V + U \text{ finite}\}$ with the Boolean median structure. □

POC SETS OF DIMENSION $\omega$

Recall that $P$ has dimension $\omega$ if all transverse subsets of $P$ are finite. If $P$ is also discrete, then there is no distinction between founded and well founded ultra filters.

**9.4. Proposition.**
*If $P$ is discrete and does not contain an infinite transverse subset, then every founded ultra filter in $P$ is also well founded.*

*Proof.* Let $U \in P^\circ$ be a founded ultra filter and $a \in U$ proper. Then the set $\{b \in \mathrm{Min}(U) \mid b \leq a\}$ is transverse and hence finite. Therefore
$$\{b \in U \mid b \leq a\} \subseteq \bigcup \{[c, a] \mid c \in \mathrm{Min}(U) \text{ and } c \leq a\}$$
is finite. □

**9.5. Theorem.**
*Suppose that $P$ is discrete and every transverse subset of $P$ is finite. Let $\mathscr{T}$ be the set of transverse subsets of proper elements of $P$. Let $\mathscr{C}$ consist of those elements of $\mathscr{T}$, which are maximal with respect to inclusion. For $A \in \mathscr{C}$ we define*
$$\tau(A) := \{b \in P \mid \exists a \in A: \quad b \geq a \quad \text{or} \quad b > a^*\},$$
*and let $M := \{\tau(A) \mid A \in \mathscr{C}\}$. Then $M$ is the median subalgebra of $P^\circ$ comprising all well founded ultra filters in $P$, and $P \cong M^\circ$.*

The condition that $P$ has dimension $\omega$ of course implies that $\mathscr{C}$ is non empty. If $P$ is nested then it has dimension 1, and we get Dunwoody's Theorem.

*Proof.* We have to check various facts.

(i) For every $A \in \mathscr{C}$, $\tau(A)$ is a well founded ultra filter.

Let $A \in \mathscr{C}$ and $b \in P$ be proper, then either $b$ or $b^*$ must be comparable with at least one element of $a$, otherwise $A \cup \{b\}$ would be transverse, and $A$ would not be maximal. Thus either $b$ or $b^*$ belongs to $\tau(A)$. Assume that $b_1, b_2 \in \tau(A)$ and $b_1 \leq b_2^*$. Then there exist elements $a_1, a_2 \in A$ and $\epsilon_1, \epsilon_2 \in \{1, *\}$ such that $b_i \geq a_i^{\epsilon_i}$ for $i = 1, 2$. It follows that

$$a_1^{\epsilon_1} \leq b_1 \leq b_2^* \leq (a_2^{\epsilon_2})^*.$$

Since $a_1$ and $a_2$ are transverse and proper, this can only hold when $a_1 = a_2$ and $\epsilon_1 \neq \epsilon_2$. But then one of the inequalities must have been strict, which is impossible.

Now we know that $\tau(A)$ is an ultra filter. It is founded by definition and well founded by the previous proposition.

(ii) Each founded ultra filter is a $\tau(A)$.

Let $U$ be a founded ultra filter; choose a maximal transverse subset $A \subseteq \mathrm{Min}(U)$. We claim that $A \in \mathscr{C}$, i.e., $A$ is maximal as a transverse subset of $P$.

Suppose that there exists an element $b \in P$ such that $A \cup \{b\}$ is transverse. Without loss of generality we may assume that $b \in U$. As $U$ is founded, we find an element $c \in \mathrm{Min}(U)$ with $c \leq b$. Now we can't have $c \in A$, because $A \cup \{b\}$ is transverse. But we must have $c \parallel a$ for some $a \in A$, otherwise $A$ is not maximal. Thus $a^* < c$, and then $a^* < b$, contradicting the assumption that $A \cup \{b\}$ is transverse, so our claim is true.

Now we know that for every $b \in U$ there exists some $a \in A$ with $a \parallel b$. If $a$ is comparable to $b$, then in fact $a \leq b$, because $a$ is minimal in $U$. If $a^*$ is comparable to $b$, then $a^* < b$, because $U$ is a filter. This proves that $U = \tau(A)$.

(iii) If $U$ and $V$ are founded ultra filters in $P$, then $U + V$ is finite.

As $U$ and $V$ are founded, we have

$$U \smallsetminus V = \bigcup \{[a, b] \mid a \in \mathrm{Min}(U \smallsetminus V),\ b^* \in \mathrm{Min}(V \smallsetminus U)\}.$$

But both $\mathrm{Min}(U \smallsetminus V)$ and $\mathrm{Min}(V \smallsetminus U)$ are transverse sets, hence finite by our hypothesis, and $U \smallsetminus V$ is a finite union of finite intervals.

(iv) If $U$ is well founded and $U + V$ is finite, then $V$ is also well founded.

If $a \in V \cap U$, then

$$\{b \in V \mid b \leq a\} \subseteq \{b \in U \mid b \leq a\} \cup (U + V),$$

and if $a \in V \smallsetminus U$, then $b \leq a$ implies $b^* \in U$, so

$$\{b \in V \mid b \leq a\} \subseteq V \smallsetminus U.$$

Altogether we have shown that $\tau(\mathscr{C})$ is an almost equality class in $P^\circ$, which consist entirely of founded ultra filters, in particular it is a median algebra. □

**Note.**

The map $\tau\colon \mathscr{C} \to M$ is usually not injective. In fact, the original construction of DUN-
WOODY [1979] starts by describing the equivalence relation $e \approx e' \Leftrightarrow \tau(e) = \tau(e')$ on the
oriented edge set $\widetilde{E}T$ of a tree in terms of the poc structure.

$$e \approx e' \qquad \Leftrightarrow \qquad e = e' \quad \text{or} \quad (e^* < e' \text{ and } [e^*, e'] = \{e^*, e'\}).$$

Dunwoody then identifies $VT$ with the set of equivalence classes $\widetilde{E}T/\approx$.

The idea to use the ultra filters $\tau(e)$ instead of $\approx$-classes originates with DICKS-DUN-
WOODY [1989], although it is somewhat disguised by their insistence on working with $ET$
rather than $\widetilde{E}T$.

COUNTABLE POC SETS

**9.6. Theorem.**
*If $P$ is a countable, discrete poc set, then $P$ has a well founded ultra filter.*

Before we prove this theorem, we explain the general idea. If we have two poc sets
$P_1 \subseteq P_2$, then there is a natural surjection $\pi\colon P_2^\circ \to P_1^\circ$. If we view $P_2^\circ$ as $\mathrm{Poc}(P_2, \mathbf{2})$, then
this map is simply the restriction to $P_1$, if we think of $P_2^\circ$ as the set of ultra filters, then
map $U \mapsto U \cap P_1$.

We now ask if we can find an embedding $\rho\colon P_1^\circ \to P_2^\circ$, such that $\pi \circ \rho = \mathrm{id}_{P_1}$. In other
words, can one extend all ultra filters on $P_1$ to $P_2$ in the same way? Consider three
elements $a < b < c$ with $a, c \in P_1$ and $b \in P_2 \smallsetminus P_1$. Then there exist ultra filters
$a \in U_1 \subset P_1$ and $c^* \in U_2 \subset P_1$, but when we extend $U_1$ we must put $b$ into it, and
when we extend $U_2$ we must put $b^*$ in. So the obstruction consists of elements of $P_2$ that
separate elements of $P_1$. If there is no such obstruction, then we call $P_1$ a convex subset
of $P_2$.

Now the idea is to build inductively, by a variation of the Extension Algorithm 3.3, an
ascending sequence $P_1 \subseteq P_2 \subseteq \ldots$, of finite, convex sub poc sets, whose union is $P$. Then
the corresponding duals also form a directed sequence $P_1^\circ \to P_2^\circ \to \ldots$, where as without
the convexity condition one would only get an inverse system $P_1^\circ \leftarrow P_2^\circ \leftarrow \ldots$. The sets
$P_i$ being finite makes all ultra filters of $P_i^\circ$ well founded. This property is shared by the
direct limit $\varinjlim P_i^\circ$. In the proof we will phrase this slightly differently. We build an ultra
filter stepwise, such that if an element $b$ is included at step $n$, then we will later only
include in elements above $b$ or incomparable to $b$, thus in the final ultra filter there will
only be finitely many elements below $b$.

One of the reasons why this works is that convex hulls in partially ordered sets are easy
to find. For any set $A \subseteq P$ just look at $C := \bigcup\{[a,b] \mid a, b \in A\}$. Clearly every convex
set containing $A$ must also contain $C$, but $C$ is also convex, so it is actually the convex
hull of $A$. If $A$ is finite and $P$ discrete, then $C$ is again finite.

*Proof of Theorem 9.6.* Recall the following definitions. For $a \in P$ let $\overline{a} := \{a, a^*\}$, and for $F \subseteq P$ let $\overline{F} := \{\overline{a} \mid a \in F\}$. A subset $\overline{F} \subseteq \overline{P}$ is called convex, if for any $\overline{a}, \overline{c} \in \overline{F}$ the interval

$$[\overline{a}, \overline{c}] := \{\overline{b} \mid a \leq b \leq c \quad \text{or} \quad a^* \leq b \leq c \quad \text{or} \quad a \leq b \leq c^* \quad \text{or} \quad a^* \leq b \leq c^*\}$$

also lies in $\overline{F}$.

Assume that $P$ is a countable poc set. We choose an enumeration for $\overline{P}$, say $\overline{P} = \{\overline{p}_0, \overline{p}_1, \ldots\}$, where $\overline{p}_0 := \overline{0}$. Assume further that $P$ is discrete, i.e., the interval between any two proper elements is finite.

We will inductively define a sequence of finite subsets $F_1 \subseteq F_2 \subseteq \ldots$, such that for all $n \in \mathbb{N}$ we have

   (i) $\overline{p}_n \in \overline{F}_{n+1}$.
   (ii) $F_n$ is a filter base.
   (iii) $\overline{F}_n$ is convex.
   (iv) If $a \in F_n$, $b \in F_{n+1}$ and $a > b$, then $b \in F_n$.

Now conditions (i) and (ii) guarantee that $F := \{0^*\} \cup \bigcup F_n$ is an ultra filter. An easy induction argument using condition (iv) says that for every $a \in F_n$ the set $\{b \in F \mid b < a\}$ is also contained in $F_n$, in particular it is finite. Thus $F$ is well founded.

Here is the construction. Start with $F_1 := \emptyset$. Assume that $F_n$ has been constructed for $n \in \mathbb{N}$. We consider the following cases:

CASE 1:   $p_n$ is transverse to all elements of $F_n$.
          Then put $F_{n+1} := F_n \cup \{p_n\}$.
CASE 2:   $\overline{p}_n \in \overline{F}_n$.
          Then do nothing, $F_{n+1} := F_n$.
CASE 3:   $\overline{p}_n \notin \overline{F}_n$, but some element $a \in F_n$ is comparable to $p_n$ or $p_n^*$.
          After replacing $p_n$ by $p_n^*$, if necessary, we may assume that either $a < p_n$ or $a^* < p_n$. Then let

$$F_{n+1} := F_n \cup \{b \in P \mid \overline{b} \notin \overline{F}_n, \text{ and } \exists a \in F_n \,:\, a \leq b \leq p_n \quad \text{or} \quad a^* \leq b \leq p_n\}.$$

All that remains is to check that everything works as we expect. It is clear that the three cases are disjoint. Could CASE 3 give conflicting advise? This means that there are $a_1, a_2 \in F_n$ such that $a_1^{\epsilon_1} < p_n$ and $a_2^{\epsilon_2} < p_n^*$. (Here $\epsilon \in \{1, *\}$ and we write $a^1 := a$.) Then $a_1^{\epsilon_1} < p_n < a_2^{\epsilon_2 *}$, but since $\overline{F}_n$ is convex, this means that $\overline{p}_n \in \overline{F}_n$, so we are in CASE 2. It is clear that the $F_n$ are all finite (by discreteness we only add a finite number of elements in CASE 3) and form an ascending sequence. Condition (i) is obviously true by our construction. Condition (ii) is also satisfied in CASES 1 and 2. To check CASE 3, we assume that there exist elements $b_1, b_2 \in F_{n+1}$ with $b_1 \leq b_2^*$. We consider three subcases.

   (a) Both $b_1, b_2 \in F_n$.
       This can't happen, because $F_n$ is already a filter base.

(b) $b_1 \in F_n$ and $b_2 \notin F_n$, so there are $a_2 \in F_n$ and $\epsilon_2 \in \{1, *\}$ with $a_2^{\epsilon_2} \leq b_2 \leq p_n$.
Then $b_1 \leq b_2^* \leq a_2^{\epsilon_2 *}$, so by convexity, $\bar{b}_2 \in \bar{F}_n$, contrary to our assumption.

(c) Both $b_1, b_2 \notin F_n$, so $a_i^{\epsilon_i} \leq b_i \leq p_n$ for $i = 1, 2$.
Then $a_1^{\epsilon_1} \leq b_1 \leq b_2^* \leq a_2^{\epsilon_2 *}$, and again $\bar{b}_1$ and $\bar{b}_2$ lie in $\bar{F}_n$.

Next we prove that $F_{n+1}$ is convex. This is trivial in CASE 2, and also clear in CASE 1, because transverse elements of $P$ can't be separated. In CASE 3 we have to show that

$$\forall b_1, b_3 \in F_{n+1} \; \forall b_2 \in P \; \forall \epsilon_1, \epsilon_2, \epsilon_3 \in \{1, *\} : \quad b_1^{\epsilon_1} < b_2^{\epsilon_2} < b_3^{\epsilon_3} \quad \Rightarrow \quad \bar{b}_2 \in \bar{F}_{n+1}.$$

Consider the following subcases.

(a) If both $b_1, b_3 \in F_n$, then this follows from the convexity of $F_n$.

(b) If $b_1 \in F_n$ and $b_3 \notin F_n$, then there exists $a_3 \in F_n$ and $\nu \in \{1, *\}$ with $a_3^\nu \leq b_3 \leq p_n$.
If $\epsilon_3 = *$, then $b_1^{\epsilon_1 *} \geq b_3 \geq a_3^{\nu *}$, thus $\bar{b}_3 \in \bar{F}_n$ and hence $\bar{b}_2 \in \bar{F}_n$.
If $\epsilon_3 = 1$, then $b_1^{\epsilon_1} \leq b_2^{\epsilon_2} \leq b_3 \leq p_n$, which implies that $\bar{b}_2 \in \bar{F}_{n+1}$.

(c) Both $b_1, b_3 \notin F_n$, so $a_i^{\nu_i} \leq b_i \leq p_n$ for $i = 1, 3$, $a_i \in F_n$ and $\nu_i \in \{1, *\}$.
Here we may assume $\epsilon_1 = 1$ without loss of generality.
If $\epsilon_3 = *$, then $a_1^{\nu_1} \leq b_1 \leq b_3^* \leq a_3^{\nu_3 *}$, thus $\bar{b}_3 \in \bar{F}_n$ and hence $\bar{b}_2 \in \bar{F}_n$.
If $\epsilon_3 = 1$, then $a_1^{\nu_1} \leq b_1 \leq b_2^{\epsilon_2} \leq b_3 \leq p_n$, so $\bar{b}_2 \in \bar{F}_{n+1}$.

Finally, condition (iv) is also trivially satisfied in CASES 1 and 2. To check it in CASE 3, assume that $b \in F_{n+1} \setminus F_n$, so there exists $a_1 \in F_n$ and $\epsilon \in \{1, *\}$ with $a_1^\epsilon \leq b \leq p_n$. If there also exists an $a \in F_n$ with $b < a$, then we have $a_1^\epsilon \leq b < a$, so in fact $b \in F_n$. This completes the proof of our construction. □

A COUNTER EXAMPLE

The condition that $P$ is countable can not be omitted in Theorem 9.6. To understand this we first study the property of being well founded in the special case of a binary poc set $P = \{0, 0^*\} \cup O \cup O^*$, see Example 1.4(vi). An ultra filter $U \subset P$ gives rise to a pair of sets $U^+ := O \cap U$ and $U^- := O \setminus U$, which partition $O$, such that $U^+$ is an upper set, i.e., whenever $u \in U^+$ and $x \in O$ with $x > u$ then $x \in U^+$, and similarly $U^-$ is a lower set. Conversely, every partition of $O$ into an upper and a lower set implies an ultra filter $\{0^*\} \cup U^+ \cup (U^-)^*$ in $P$.

The condition that $U$ is well founded translates into the condition that $U^+$ is down finite, i.e., for every $u \in U^+$ the set $\{x \in U \mid x \leq u\}$ is finite, and similarly $U^-$ is up finite. Now we give an example of a discrete partially ordered set $O$, that does not permit a partition $O = U^+ \cup U^-$ with $U^+$ a down finite upper set and $U^-$ an up finite lower set. Consequently, the corresponding binary poc set has no well founded ultra filters, so it is not the dual of a discrete median algebra.

### 9.7. Example.
Let $I$ be a countable infinite set and $J$ an uncountable set. On the set $O := I \sqcup (I \times J) \sqcup J$ we define an ordering by declaring
$$i < (i,j) < j \qquad \text{for all } i \in I \text{ and } j \in J.$$
In particular, $O$ is discrete, as any interval contains at most 3 elements.

Suppose $U^+$ and $U^-$ partition $O$ as desired. Assume first, that some element $i \in I$ lies in $U^+$. Since $U^+$ is an upper set, $(i,j) \in U^+$ for all $j \in J$, and then $J \subset U^+$. Since $U^+$ is down finite, only finitely many elements of $I$ can lie in $U^+$. We may swap these elements from $U^+$ to $U^-$, thus without loss of generality we may assume that $I \subseteq U^-$ and $J \subseteq U^+$. For each $i \in I$ the set $\{(i,j) \in I \times J \mid (i,j) \in U^-\}$ is finite, so $X := (I \times J) \cap U^-$ is a countable union of finite sets, hence countable. The projection of $X$ onto $J$ can't be surjective, since $J$ is uncountable. Thus there exists an element $j_0 \in J$ such that $I \times \{j_0\} \subset U^+$, but every pair $(i, j_0)$ is smaller than $j_0 \in U^+$, hence $U^+$ can't be down finite. Contradiction! $\square$

### Discrete Representations

A **representation** of a poc set $P$ on a set $X$ is a poc embedding $\rho \colon P \to \mathscr{P} X$. For $x, y \in X$ we can define $\Delta(x, y) := \{a \in P \mid x \in \rho(a) \text{ and } y \in \rho(a^*)\}$. We say that $\rho$ is a **discrete representation** if $\Delta(x, y)$ is finite for all $x, y \in X$.

### 9.8. Proposition.
*A poc set $P$ admits a discrete representation if and only if there exists a discrete median $M$ with $M^\circ \cong P$.*

*Proof.* If $M$ is a discrete median algebra, then the inclusion $M^\circ \subset \mathscr{P} M$ is of course a discrete representation. Conversely suppose that $\rho \colon P \to \mathscr{P} X$ is a discrete representation. We have a canonical map $\iota \colon X \to P^\circ$, with $\iota(x) := \{a \in P \mid x \in \rho(a)\}$. The image of $\iota$ lies in an almost equality class, because $\iota(x) \smallsetminus \iota(y) = \Delta(x, y)$ is finite. Further more, every $\iota(x)$ is a well founded ultra filter in $P^\circ$. By Proposition 9.3 the whole almost equality class of $\iota(X)$ in $P^\circ$ is the desired discrete median algebra $M$. $\square$

We will meet these discrete representations again in the context of groups. Here we add just a few remarks on metric aspects.

### 9.9. Remarks.
(i) On the set $X$ we have a pseudo metric $d(x,y) := |\Delta(x,y)|$, i.e., $d$ satisfies the triangle inequality (this is a consequence of the Lemma 2.16, which holds with the same proof), but we may have $d(x,y) = 0$ for distinct elements $x, y \in X$. Write $X_d$ for the quotient of $X$ modulo the equivalence relation $x \sim y \Leftrightarrow d(x,y) = 0$ and let $M_\rho$ denote the median algebra provided by the previous proposition. We have an induced map $\iota\colon X_d \to M_\rho$ which is injective. The metric $d$ on $X$ is extended by the metric $d$ on $M_\rho$.

(ii) Proposition 7.12 says that $M_\rho$ is generated as a median algebra by $\iota(X)$.

(iii) Apart from the triangle inequality, the pseudo metric $d\colon X \times X \to \mathbb{N}_0$ also satisfies Dunwoody's pattern condition (DICKS-DUNWOODY [1989]).

$$\forall x, y, z \in X: \quad d(x,y) + d(y,z) + d(z,x) \text{ is even.}$$

It is conceivable that this condition actually characterizes the pseudo metrics arising from discrete representations.

(iv) If one starts with a (pseudo) metric space that satisfies the pattern conditions, then the implied poc set or median algebra are not uniquely determined. For example, the simplex $X := \{x_1, x_2, x_3, x_4\}$ with $d(x_i, x_j) = 2$ can generate two different median algebras.

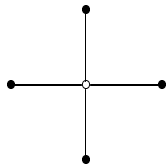
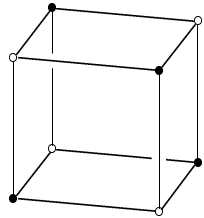

# §10. Cubings

The main object of SAGEEV [1995] is a cubing, i.e., a cube complex satisfying a generalized CAT(0) condition. Sageev constructs a cubing out of a discrete poc set by taking an almost equality class of its dual. Although in this work we take a more algebraic point of view, we will recast Sageev's construction here, because it gives a geometric meaning to discrete median algebras. First we recall some definitions.

By the $n$-dimensional unit cube we mean the set $[0,1]^n \in \mathbb{R}^n$, with its usual $CW$-structure as a finite product of intervals. A cube complex is a CW-complex, all whose characteristic functions have unit cubes as their domains and are isomorphisms of cell complexes onto their image.

Following BRIDSON-HAEFLIGER [1995], we define a Euclidean metric on $K$, such that every cell is in fact an isometric copy of a unit cube, and where two cubes meet the metrics agree. There is a notion of a piecewise linear path, whose length can be measured in the obvious way inside each cube. The distance of two points is then the infimum of the lengths of piecewise linear paths connecting them. The arguments of BRIDSON-HAEFLIGER [1995], § 1.7, can be used to show that this defines a metric, even when $K$ is infinite dimensional. According to a theorem of Bridson, this metric space is complete if $K$ is finite dimensional, even if it is not locally finite. For infinite dimensional cube complexes the Euclidean metric need not be complete.

The topology induced by the Euclidean metric on the underlying set of $K$ does in general not agree with the CW-topology. However, we are only interested in the homotopy type of $K$, and by Smale's version of the Vietoris Mapping Theorem the identity map $K_{Eucl} \to K$ is a homotopy equivalence.

Our next condition on cube complexes is a local condition, and to express it we need to recall some terminology for simplicial complexes. Let $X$ be a graph with vertex set $VX$. A flag in $X$ is the vertex set of a complete subgraph, i.e., a set $F \subseteq VX$ where any pair $x_1, x_2 \in F$ is adjacent in $X$. The flag complex of $X$ is the simplicial complex flag$(X)$ on the vertices of $X$, whose simplexes are the finite flags of $X$. In particular, the 1-skeleton of flag$(X)$ is $X$ itself.

Let $K$ be a cube complex. The link of a cube $C \in K$ is a simplicial complex link$(C, K)$, whose vertices are the cubes $C'$ which strictly contain $C$, and $\{C_0, C_1, \ldots, C_n\}$ are the vertices of a simplex in link$(C, K)$ if there exists a cube $C'$ in $K$ with $C_i \subseteq C'$ for all $i = 0, 1, \ldots, n$.

BRIDSON-HAEFLIGER [1995] show that a finite dimensional cube complex $K$ satisfies the local CAT(0) condition iff the link of every vertex in $K$ is a flag complex. (See also GROMOV [1987], 4.2.C. Gromov uses the no-bigons and no-triangle condition, but that is equivalent to the flag condition.)

A metric space satisfies the global CAT(0) condition if it enjoys the local CAT(0) condition and is simply connected. The uniqueness of geodesics in CAT(0) spaces then shows that the space is in fact contractible.

A cubing, in the sense of Sageev, is a simply connected cube complex $K$ so that for every vertex $v \in K$ the link $\text{link}(v, K)$ is a flag complex.

We now show how to associate a cubing to each discrete median algebra. As in §8, a finite cube is the median algebra of the power set $\mathscr{P}X$ on some finite set $X$. If two cubes meet, then their intersection is a subcube of both. Thus the finite subcubes of a discrete median algebra $M$ form a combinatorial cube complex, and we denote the realisation of this complex as a CW-complex by $\mathscr{K}M$, the cubical nerve of $M$. The 0-skeleton of $\mathscr{K}M$ is of course $M$, and 1-skeleton of $\mathscr{K}M$ is the median graph associated to $M$.

Let $x \in M$ and consider the star $S := \text{star}(x, M)$ as defined in §8. A half spaces $H \in M^\circ$ contains $x$ and cuts $S$ precisely when $H$ is a minimal element of $ev(x)$. In Proposition 8.9 we showed that the cubes in $M$ which contain $x$ correspond bijectively to the finite transverse subsets of $\text{Min } ev(x)$, the set of minimal elements of $ev(x)$. This can be rephrased in the following way.

**10.1. Proposition.**
$\text{link}(x, \mathscr{K}M)$ is the flag complex of the transversality graph of $\text{Min } ev(x)$. □

**10.2. Proposition.**
The cube complex $\mathscr{K}M$ is contractible.

*Proof.* As $\mathscr{K}M$ is a CW-complex, it is enough to show that every finite subcomplex is contained in a contractible complex. But a finite subcomplex of $\mathscr{K}M$ involves only finitely many elements of $M$, thus it lies in the cube complex obtained from a finite subalgebra of $M$. Thus it is enough to show that $\mathscr{K}M$ is contractible whenever $M$ is finite.

Choose a proper half space $H \in M^\circ$. The cubes in $M$ fall into three classes:

  (i) Those contained in $H$,
 (ii) those contained in $H^*$, and
(iii) those cut by $H$.

The cubes in the first two classes make up the subcomplexes $\mathscr{K}H$ and $\mathscr{K}H^*$ of $\mathscr{K}M$. If $C$ is a cube in the third class, then for any point $x \in C \cap H$ there is a unique point $x' \in C \cap H^*$ with $\Delta(x, x') = \{H\}$. The map $C \cap H \to C \cap H^*$ sending $x$ to $x'$ is an isomorphism of median algebras. Thus $C$ is a product $\mathbf{2} \times (C \cap H)$, and the cubes in the third class make up a subcomplex homeomorphic to $I \times \mathscr{K}\partial H$, where $I = [0, 1]$ is the cube complex of the median algebra $\mathbf{2}$.

We form a new median algebra $M'$ by identifying every point $x \in \partial H$ with its opposite point $x' \in H^*$. $M'$ can also be described as the quotient $M/\triangledown\{H, H^*\}$, see §6. The natural quotient map $q\colon M \to M'$ induces a continuous map $q\colon \mathscr{K}M \to \mathscr{K}M'$, which collapses the factor $I$ in the "handle" $I \times \mathscr{K}\partial H$. Clearly, this is a homotopy equivalence.

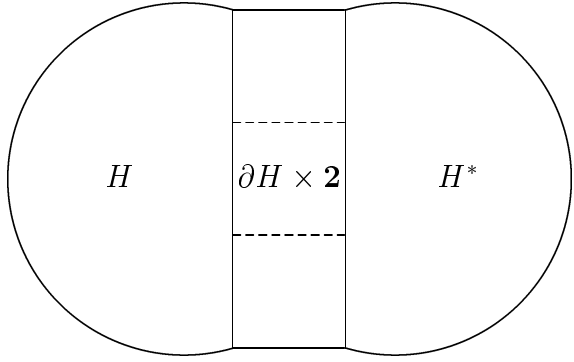

If $M$ is a finite median algebra, we can perform this step inductively, reducing the number of proper half spaces in $M$ each time, until we finish with a median algebra that consists of a single element, whose cube complex is a single point. Thus we have shown that $\mathscr{K}M$ is contractible. □

Putting the two propositions together we have the first half of the following theorem.

**10.3. Theorem.**
*If $M$ is a median algebra, then $\mathscr{K}M$ is a cubing. Conversely, for every cubing $K$ the 0-skeleton $M := K^0$ carries a natural structure of a median algebra such that $K = \mathscr{K}M$.*

*Proof.* The converse relies heavily on §4 of SAGEEV [1995]. He introduces a notion of geometric hyper planes in a cubing $K$ and shows that hyper planes embed as totally geodesic codimension-1 subspaces, which dissect $K$ into half spaces. Our definition here follows NIBLO-REEVES [1996]. Inside a unit cube $I^n$ one has $n$ midplanes $M_i := \{(x_1, \ldots, x_n) \in I^n \mid x_i = \frac{1}{2}\}$ (Sageev calls them dual blocks), and we also call their images under the characteristic maps of $K$ midplanes. For two midplanes $M_1$ and $M_2$ of two cubes $C_1$ and $C_2$ in $K$ we write $M_1 \sim M_2$ if $M_1 \cap M_2$ is again a midplane (and then it is a midplane of $C_1 \cap C_2$). The transitive closure of this symmetric relation is an equivalence relation, and the union of all midplanes in an equivalence class is called a geometric hyper plane.

Using surgery and minimal discs arguments, Sageev proves several important properties of geometric hyper planes. First, Theorem 4.10 of SAGEEV [1995],

(i) The complement of a geometric hyper plane in $K$ has two connected components.

If the two components are $K_1$ and $K_2$, then we call the sets $M \cap K_1$ and $M \cap K_2$ geometric half spaces, where $M$ is the 0-skeleton of $K$. These half spaces form the proper elements of a poc set $P$ with the inclusion order. We will use the symbols $H$ and $H^*$ to denote a geometric half space and its complement, and $\overline{H}$ to denote the corresponding geometric hyperplane.

As $K$ is a connected cube complex, any two vertices $x$ and $y$ in $M$ can be connected by a path in the 1-skeleton $K^1$. A path with the minimal number, $d(x,y)$ of edges is called a geodesic. Let $\overline{\Delta}(x,y)$ denote the set of hyper planes that separate $x$ from $y$. It is clear that any geodesic from $x$ to $y$ must meet all hyper planes in $\overline{\Delta}(x,y)$, and from Theorem 4.13 of SAGEEV [1995] it follows that a geodesic crosses a hyper plane at most once, so

  (ii) $d(x,y) = |\overline{\Delta}(x,y)|$.

For $x, y \in M$ define the interval $[x,y] := \{z \in M \mid d(x,z) + d(z,y) = d(x,y)\}$. The next statement expresses the fact that intervals are convex.

  (iii) $z \in [x,y]$ iff no hyper plane separates $z$ from $\{x,y\}$.

*Proof.* With a proof identical to that of Lemma 2.16, we have
$$\overline{\Delta}(x,y) = \overline{\Delta}(x,z) + \overline{\Delta}(z,y).$$
Thus $d(x,y) = d(x,z) + d(z,y)$ is equivalent to
$$|\overline{\Delta}(x,z) + \overline{\Delta}(z,y)| = |\overline{\Delta}(x,z)| + |\overline{\Delta}(z,y)|,$$
and that means $\overline{\Delta}(x,z) \cap \overline{\Delta}(z,y) = \varnothing$, which proves (iii).

The next statement captures the idea that geometric hyper planes are totally geodesic; Theorem 4.13 of SAGEEV [1995].

  (iv) Let $\{x,y\}$ and $\{x',y'\}$ be edges in $K^1$, whose mid points lie in the same hyper plane, such that $x$ and $x'$ lie in the same component of the complement. Then for every geodesic $x_0 := x$, $x_1$, ..., $x_n := x'$ there exists a geodesic $y_0 := y$, $y_1$, ..., $y_n := y'$ such that $\{x_{i-1}, x_i, y_{i-1}, y_i\}$ are the vertices of a square in $K^2$ for all $i = 1, \ldots, n$.

We use this to prove

  (v) Let $e_1 := \{x,y\}$ and $e_2 := \{x,z\}$ be two distinct edges in $K^1$ with a common vertex $x$. Let $\overline{H}_1$ and $\overline{H}_2$ be the hyper planes containing the midpoints of $e_1$ and $e_2$, respectively. Then $\overline{H}_1$ and $\overline{H}_2$ meet if and only if there is a square in $K^2$ containing $e_1$ and $e_2$ as edges.

*Proof.* If $e_1$ and $e_2$ are edges of a square, then $\overline{H}_1$ and $\overline{H}_2$ both contain the midpoint of this square.
Conversely, assume that $\overline{H}_1$ and $\overline{H}_2$ meet in a square $\{x', y', z', w'\}$. We can choose components $H_1$ and $H_2$ and labels for the vertices of the square, in such a way that
$$x, x' \in H_1 \cap H_2, \quad y, y' \in H_1^* \cap H_2, \quad z, z' \in H_1 \cap H_2^*, \quad w' \in H_1^* \cap H_2^*.$$
Choose a geodesic $x_0 := x$, $x_1$, ..., $x_n := x'$. As in (iv) we find further geodesics $y_0 := y$, $y_1$, ..., $y_n := y'$ and $z_0 := z$, $z_1$, ..., $z_n := z'$
Now the edges $\{x_{n-1}, x_n\}$, $\{x_n, y_n\}$ and $\{x_n, z_n\}$ lie pairwise in a common square. Here we apeal to the local CAT(0) condition, which says that they are edges of a cube in $K^3$.

In particular, there exists a vertex $w_{n-1}$ such that $\{x_{n-1}, y_{n-1}, z_{n-1}, w_{n-1}\}$ form a square. Working our way backwards inductively, we find a square $\{x, y, z, w\}$, which proves (v).

(vi) For all $x, y, z \in M$ the intersection $[x,y] \cap [y,z] \cap [z,x]$ is nonempty.

*Proof.* Choose an $m \in [x,y] \cap [x,z]$ such that $d(x,m)$ is as large as possible. Suppose that $m \notin [y,z]$, then by (iii) there exists a hyper plane $\overline{H}$ separating $m$ from $\{y,z\}$. Choose such an $\overline{H}$ and a point $w$ which is separated from $m$ by $\overline{H}$, such that $d(m,w)$ is as small as possible.

We claim that in fact $d(m,w) = 1$. If not, then the geodesic $w_0 := w$, $w_1$, ..., $w_n := m$ has length at least 2. By minimality of $n$, $\overline{H}$ separates $w_0$ from $w_1$. Let $\bar{J}$ denote the hyper plane that separates $w_1$ from $w_2$. If $\overline{H}$ and $\bar{J}$ are disjoint, then $\bar{J}$ also separates $m$ from $\{y,z\}$, thus $n$ was not minimal. If $\overline{H}$ and $\bar{J}$ meet, then according to (v) they meet in a square containing $w_0$, $w_1$ and $w_2$. The fourth vertex, $w_3$, is separated from $m$ by $\overline{H}$, but closer to $m$ than $w_0$, again contradicting the minimality of $n$. Thus we have proved our claim that $\overline{\Delta}(m,w) = \{\overline{H}\}$.

Next we claim that $w \in [x,y]$. If not, then by (iii) there exists a hyper plane separating $w$ from $\{x,y\}$, say $w \in J$ and $x,y \in J^*$. We can't have $m \in J$, because $m \in [x,y]$. Thus $\bar{J}$ separates $m$ and $w$, which means, as we have shown in the previous paragraph, that $\bar{J} = \overline{H}$. But $\overline{H}$ also separates $m$ from $\{y,z\}$, so $y \in J$, which is a contradiction.

By the same argument we can show that $w \in [x,z]$. But $d(x,w) = d(x,m) + 1$, which contradicts our choice of $m$. Thus we must have $m \in [y,z]$, so $m$ lies in the intersection of all three intervals. Hence (vi) is proved.

It remains to show that $m =: m(x,y,z)$ is unique and satisfies the median axioms. Here we employ the same trick as in our model example of §1. Let $\mathscr{H}$ denote the set of all half spaces that contain at least two out of $x$, $y$ and $z$. Then any $H \in \mathscr{H}$ contains a whole interval, and in particular the intersection of all three intervals. Thus

$$\bigcap \mathscr{H} \supseteq [x,y] \cap [y,z] \cap [z,x].$$

On the other hand, statement (ii) says that any two distinct elements of $M$ are separated by a hyper plane $\overline{H}$, and exactly one of $H$ and $H^*$ belongs to $\mathscr{H}$, so $\bigcap \mathscr{H}$ contains at most one element. Hence $m$ is unique.

But now it is clear, that the map $x \mapsto \{H \in P \mid x \in H\}$ embeds $M$ into the (Boolean) median algebra $P^\circ$, and the median half spaces coincide with the geometric half spaces. This proves our theorem. □

**Notes.**
MULDER [1978] characterizes finite median algebras by a blowing up process that is the inverse to our contraction in Proposition 10.2. NIBLO-REEVES [1996] give an account of Sageev's work for finite dimensional cubings. Here one can use the Cartan-Hadamard Theorem to prove that hyper planes embed and are totally geodesic. GERASIMOV [1997] contains a proof of Theorem 10.3 using different methods.

# §11. Groups Acting on Discrete Median Algebras

In this section we study the dynamics of the action of a group $G$ on a discrete median algebra $M$. Viewing the theory of group actions on simplicial trees as a special case of our theory, we try to follow the survey in ROLLER [1993] as closely as possible — keeping in mind that new difficulties arise from the possibility that infinitely many half spaces may meet transversally.

Earlier we introduced the median graph $\Gamma M$, whose vertices are the elements of $M$, and two distinct vertices $x$ and $y$ are joined by an edge when $[x,y] = \{x,y\}$. $\Gamma M$ is connected iff $M$ is discrete; in this case the edge metric on $\Gamma M$ agrees with the metric $d(x,y) = |\Delta(x,y)|$ defined on $M$. A median automorphism $\varphi\colon M \to M$ preserves adjacency in $\Gamma M$, hence it induces a graph automorphism of $\Gamma M$. Conversely, if $\Gamma M$ is connected, any graph automorphism of $\Gamma M$ maps intervals to intervals and thus preserves medians in $M$. Thus, a group action on a discrete median algebra is the same as a group action on the corresponding median graph.

We introduce some terminology to describe the qualitative behaviour of the group action with respect to half spaces. We call a subset $X \subseteq M$ bounded, if it has bounded diameter, i.e., if there exists an $N \in \mathbb{N}$ such that $d(x,y) \leq N$ for all $x,y \in X$.

(i) We say that a half space $H \in M^\circ$ cuts $M$ properly if for some $x \in M$ (and then for all $x \in M$) both $H \cap Gx$ and $H^* \cap Gx$ are unbounded.

(ii) We say that an element $g \in G$ shifts a half space $H$ if $H$ and $gH$ are distinct and comparable, i.e., either $H \subset gH$ or $H \supset gH$.

(iii) A half space $H$ is called essential if $H$ cuts properly and the orbit $GH$ is infinite.

It is well known that for group actions on trees the three conditions are equivalent. For median algebras one can still show that shifted half spaces cut properly and have infinite orbit, but the converse need not be true. There are two reasons, why this could happen: $H$ could have a finite orbit under the action of $G$, or the orbit of $H$ consists of pairwise transverse elements, in which case the median algebra contains an infinite cube.

**11.1. Proposition.**
If $g$ shifts $H$ then $H$ is essential.

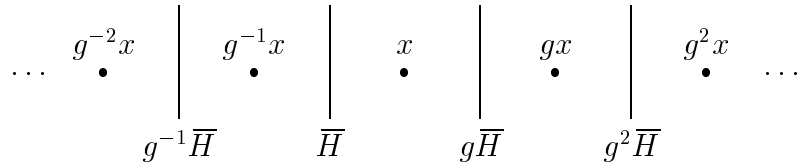

*Proof.* Suppose we can choose an element $x \in H \smallsetminus gH$. Then we have, schematically, the above picture. Thus $d(x, g^n x) \geq |n|$, and both $H \cap Gx = \{g^n x \mid n \in \mathbb{N}_0\}$ and $H^* \cap Gx = \{g^{-n} \mid n \in \mathbb{N}\}$ are unbounded, so $H$ cuts properly and $GH$ contains the infinite set $\{g^n H \mid n \in \mathbb{Z}\}$. □

## 11.2. Examples.
(i) Let $X$ be an infinite set. Let $M := \mathscr{F}X$, the family of finite subsets of $X$ with the Boolean median structure and $G := \mathscr{F}X$ with the Abelian group structure defined by the symmetric difference. It is easy to check that the action of $G$ on $M$ preserves the median. Every half space of $M$ is either of the form $H_x = \{a \in \mathscr{F}X \mid x \in a\}$ or $H_x^*$, and $GH_x = \{H_x, H_x^*\}$, moreover every $H_x$ cuts properly because $G$ acts transitively on $M$.

(ii) Consider the median algebra $\mathscr{F}\mathbb{Z}$, whose dual is orthogonal, so no half space can be shifted. The map $\sigma \colon \mathbb{Z} \to \mathbb{Z}$, $n \mapsto n+1$ induces a morphism $\sigma \colon \mathscr{F}\mathbb{Z} \to \mathscr{F}\mathbb{Z}$. Let $\rho \colon \mathscr{F}\mathbb{Z} \to \mathscr{F}\mathbb{Z}$, $a \mapsto a + \{0\}$. Then $g := \rho \circ \sigma$ is an automorphism of $\mathscr{F}\mathbb{Z}$. Consider $\mathscr{F}\mathbb{Z}$ under the action of the cyclic group $G$ generated by $g$.
A simple computation gives
$$g^n(\varnothing) = \begin{cases} \varnothing & \text{for } n = 0, \\ \{0, \ldots, n-1\} & \text{for } n > 0, \\ \{n, \ldots, -1\} & \text{for } n < 0. \end{cases}$$
The median metric on $\mathscr{F}\mathbb{Z}$ is given by $d(x, y) = |x + y|$, so in particular the orbit of $\varnothing$ under $G$ is unbounded. For the half spaces $H_n := \{x \in \mathscr{F}\mathbb{Z} \mid n \in x\}$ we have
$$g^n H_0 = \begin{cases} H_n & \text{for } n \geq 0, \\ H_n^* & \text{for } n < 0. \end{cases}$$
In particular $H_0 \cap G\varnothing = \{g^n \varnothing \mid n \in \mathbb{N}\}$ and $H_0^* \cap G\varnothing = \{g^{-n}\varnothing \mid n \in \mathbb{N}_0\}$, which are both unbounded, so $H_0$ is essential. □

In the absence of the above mentioned obstructions we can indeed find shifted half spaces.

## 11.3. Proposition.
*Assume that $M$ has dimension $\omega$. Suppose there exists a half space $H \subset M$ whose orbit $GH$ is infinite. Then some $g \in G$ shifts $H$.*

*Proof.* By Ramsey's argument there exists an infinite subset $\mathscr{S} \subseteq GH$ which is either pairwise transverse or pairwise nested. The first possibility is excluded by the assumption that $M$ has dimension $\omega$, hence the second possibility must hold. In particular there exist three distinct and pairwise nested half spaces $g_1 H$, $g_2 H$ and $g_3 H$. At least two out of these must be comparable, say $g_1 H \subset g_2 H$, and then $g_1^{-1} g_2$ shifts $H$. □

Next we classify group actions on discrete median algebras.
  (i) We say that the action of $G$ on $M$ is elliptic if for some $x \in M$ (and then for all $x \in M$) the orbit $Gx$ is bounded.

(ii) The action of $G$ on $M$ is called parabolic if $M$ is not elliptic and no half space cuts $M$ properly.

(iii) We say that the group action on $M$ is essential if there exists an essential half space.

These names are justified by the following results. We will show that $G$ acts elliptically on $M$ iff there exists a finite $G$-fixed cube in $M$, in other words, if there exists a a $G$-fixed point in the cubing $\mathscr{K}M$.

A group that acts parabolically on $M$ has a canonical fixed point in $M^{\circ\circ}$ which does not lie in $ev(M)$. In this case both $M$ and $G$ must be large, in the sense that $G$ can't be finitely generated and $M$ has an infinite number of $G$-orbits of half spaces.

The most interesting case is when $G$ acts essentially on $M$. We will characterize an essential half space $H$ by the way in which its stabilizer $G_H$ is embedded in $G$. To formulate this algebraically we need the theory of ends of pairs of groups.

The above Example 11.2(i) shows that there are group actions on discrete median algebras which don't fit into any of these three categories, i.e., they have unbounded orbits and half spaces which cut properly, but none of these half spaces has an infinite orbit. We will have nothing further to say about these actions.

Ends of Pairs of Groups

Let $K$ be a subgroup of $G$ and $K\backslash G := \{Kg \mid g \in G\}$ the set of right cosets of $K$. We identify the set $\mathscr{P}K\backslash G$ with the set $\{A \subseteq G \mid A = KA\}$ of $K$-invariant subsets of $G$ (i.e., a set of cosets is identified with its union). As usual, $\mathscr{P}K\backslash G$ carries the structure of a right $\mathbb{F}_2 G$-module, where $\mathbb{F}_2$ is the field of two elements, $G$ acts by multiplication on the right and addition is the symmetric difference.

An element $A \in \mathscr{P}K\backslash G$ is called $K$-finite, if $A = KF$ for some finite subset $F \subseteq G$; those sets form a submodule $\mathscr{F}K\backslash G$. The number of ends of the pair $(G, K)$ is defined as the dimension of the $G$-fixed submodule of the quotient.

$$e(G, K) := \dim_{\mathbb{F}_2} \left(\frac{\mathscr{P}K\backslash G}{\mathscr{F}K\backslash G}\right)^G.$$

**11.4. Remarks.**
(i) Observe that $e(G, K) \geq 1$ if and only if $K$ has infinite index in $G$.

(ii) If $K$ has infinite index in $G$ then we have the following formula for $e(G, K)$ involving the first cohomology of $G$ and $K$.

$$e(G, K) = 1 + \dim_{\mathbb{F}_2} \ker\left(H^1(G, \mathscr{F}K\backslash G) \to H^1(K, \mathbb{F}_2)\right).$$

This follows from the exact sequence

$$0 \longrightarrow \mathbb{F}_2 \longrightarrow \left(\frac{\mathscr{P}K\backslash G}{\mathscr{F}K\backslash G}\right)^G \longrightarrow H^1(G, \mathscr{F}K\backslash G) \longrightarrow H^1(G, \mathscr{P}K\backslash G)$$

and the fact that $\mathscr{P}K\backslash G$ is isomorphic to the coinduced module $\mathrm{Coind}_K^G \mathbb{F}_2$. See Krop-holler-Roller [1989, 1996] for more information about the cohomological nature of ends.

(iii) For a finitely generated group the invariant $e(G, K)$ also has a geometric interpretation. Choose some finite set $X$ of generators of $G$ and let $\Gamma(G, X)$ be the Cayley graph of $G$, then $e(G, K)$ measures the number of ends of the quotient space $K\backslash\Gamma$. In particular, $e(G, \{1\})$ is the classical number of ends of $G$. For more information about this geometric interpretation we refer to Epstein [1962] and Scott [1977]. □

A subset $A \in \mathscr{P}K\backslash G$ for which the symmetric difference $A + Ag$ is $K$-finite for all $g \in G$ is called *K-almost invariant*, those subsets represent $G$-invariant elements in the quotient $(\mathscr{P}K\backslash G)/(\mathscr{F}K\backslash G)$. A $K$-almost invariant set is *trivial* if either $A$ or $G \smallsetminus A$ is $K$-finite, otherwise $A$ is a *K-proper* almost invariant set. Thus to say that $e(G, K) > 1$ means that there exists a subset $A \subset G$, with the following properties.

(**End 1**)   $A$ is $K$-invariant.

(**End 2**)   $A$ is $K$-almost invariant.

(**End 3**)   $A$ is $K$-proper.

Now we can formulate the main results of this section. As we will be interested in a single orbit of half spaces only, it will be useful to consider median algebras where this is in fact the only orbit of proper half spaces. We say that $M$ is a *simple $G$-median algebra*, if for any two proper half spaces $H, J \in M^\circ$ there exists an element $g \in G$ such that $J = gH$ or $J^* = gH$, in other words, if $G$ acts transitively on the hyper planes of $M$.

### 11.5. Theorem.
*Let $G$ be a group. The following are equivalent.*

(a) *There exists a subgroup $K < G$ with $e(G, K) > 1$.*

(b) *$G$ admits a non elliptic action on a simple discrete median algebra.*

(c) *$G$ admits an essential action on a discrete median algebra.*

*In case (b) let $H \in M^\circ$ be any proper half space, and in case (c) let $H$ be an essential half space; then $e(G, G_H) > 1$. Conversely, if $e(G, K) > 1$, then the stabilizer $G_H$ of the half space $H$ provided in (b) and (c) contains $K$ as a subgroup of finite index.*

### 11.6. Theorem.
*Let $G$ be a finitely generated group. The following are equivalent.*

(a) *There exists a subgroup $K < G$ with $e(G, K) > 1$.*

(b) *$G$ admits a non elliptic action on a discrete median algebra.*

Before we prove these theorems we have to justify our claims about elliptic and parabolic groups in the next two subsections.

Clearly, if there exists a $G$-fixed point in $M$, or, more generally, if there exists a finite orbit, then the action must be elliptic. A well known theorem of Tits says that for trees the converse is true: An elliptic action leaves either a vertex or an edge fixed. The proof is very simple: A bounded orbit is contained in a subtree of finite diameter; the action of $G$ permutes the leaves of this tree, i.e., the vertices of valency one, amongst themselves, and the remaining vertices span a $G$-tree of smaller diameter, thus a $G$-tree of minimal diameter is either a single vertex or a single edge.

We will extend this result to discrete median algebras. However, we can't use valencies in our situation and have to use topological arguments instead.

### 11.7. Proposition.
Let $M$ be a discrete median algebra, $x \in M$ and $r \in \mathbb{N}_0$. The ball $B = B_r(x) := \{y \in M \mid d(x, y) \leq r\}$ of radius $r$ is compact in the convex topology on $M$.

*Proof.* We will show that the set of intersections $H \cap B$, where $H$ is a half space in $M$, satisfies the finite intersection property (FIP). Let $\mathscr{H}$ be a set of half spaces in $M$, such that

$$(*) \qquad \forall n \in \mathbb{N} \quad \forall H_1, \ldots, H_n \in \mathscr{H} \;:\; B \cap H_1 \cap \ldots \cap H_n \neq \varnothing.$$

We have to show that $B$ meets the intersection of $\mathscr{H}$.

We claim, that for all $y \in B$ the set $\mathscr{H}_y := \{H \in \mathscr{H} \mid y \notin H\}$ has at most $2r$ elements. For if $H_1, \ldots, H_n \in \mathscr{H}_y$ are distinct, then by $(*)$ there exists a $z \in B \cap H_1 \cap \ldots \cap H_n$. This means that $H_1, \ldots, H_n \in \Delta(z, y)$, so we have

$$n \leq d(y, z) \leq d(y, x) + d(x, z) \leq 2r,$$

which proves our claim.

By assumption $(*)$, $B$ meets $\bigcap \mathscr{H}_x$. Choose any $y \in B \cap \bigcap \mathscr{H}_x$, then every $H \in \mathscr{H}$ contains either $x$ or $y$. Let $\mathscr{H}_0 := \mathscr{H} \smallsetminus (\mathscr{H}_x \cup \mathscr{H}_y)$, then $\bigcap \mathscr{H}_0$ contains both $x$ and $y$, $x \in \bigcap \mathscr{H}_y$, $y \in \bigcap \mathscr{H}_x$, and $\bigcap \mathscr{H}_x$ meets $\bigcap \mathscr{H}_y$, according to $(*)$. By Helly's Theorem 2.2 there exists an element $z \in \bigcap \mathscr{H} = \bigcap \mathscr{H}_0 \cap \bigcap \mathscr{H}_x \cap \bigcap \mathscr{H}_y$.

Now the median $m := m(x, y, z)$ lies in $[z, x] \cap [z, y]$, hence $m$ lies in every $H \in \mathscr{H}$. Finally $d(x, m) \leq d(x, y) \leq r$, so $m \in B \cap \bigcap \mathscr{H}$. □

### 11.8. Lemma.
The set of accumulation points $B_r(x)'$ is contained in $B_{r-1}(x)$.

*Proof.* As $B_r(x)$ is compact, we have $B_r(x)' \subseteq B_r(x)$. Let $y \in M$ with $d(x, y) = r$. In view of Lemma 7.9 it is clear that $\{y\} = B_r(x) \cap \bigcap \Delta(y, x)$. But $\Delta(y, x)$ is a finite set of half spaces, which are open in the convex topology, hence $\{y\}$ is an open subset of $B_r(x)$ and $y$ can't be an accumulation point. □

**11.9. Theorem.**
*Let $G$ be a group that acts on a discrete median algebra $M$. If there exists a bounded $G$-subset of $M$, then there also exists a finite $G$-cube in $M$.*

*Proof.* We first claim that there exists a finite $G$-subset of $M$. Let $B \subseteq M$ be a bounded, infinite subset that is permuted by the action of $G$. Then there exists an element $x \in M$ and an $r > 0$ such that $B \subseteq B_r(x)$. Consider the sequence of derived sets $B^{(1)} := B'$, and $B^{(i+1)} := B^{(i)\prime}$, they are all $G$-subsets of $M$. As $B_r(x)$ is compact, we have $B_1 \subseteq B_r(x)' \subseteq B_{r-1}(x)$, and inductively $B^{(i)} \subseteq B_{r-i}(x)$. Thus after at most $r$ steps we have $B^{(i+1)} = \varnothing$. Thus in the penultimate step we have a set $B^{(i)}$ whose derived set is empty, so $B^{(i)}$ must be finite, which proves our claim.

Now let $C$ be the convex closure of a finite $G$-subset of $M$, so $C$ is again finite a $G$-set. Let
$$\mathscr{H} := \{H \in M^\circ \,\big|\, |C \cap H| > \tfrac{1}{2}|C|\}.$$
Suppose that $H_1, H_2 \in \mathscr{H}$, but $H_1 \cap H_2 \cap C = \varnothing$. Then $H_1 \cap C \subseteq H_2^* \cap C$, thus
$$|C| = |C \cap H_2| + |C \cap H_2^*| \geq |C \cap H_2| + |C \cap H_1| > \tfrac{1}{2}|C| + \tfrac{1}{2}|C| = |C|,$$
which is a contradiction. The family $\{C \cap H \mid H \in \mathscr{H}\}$ is a $G$-set consisting of finitely many convex subsets of $C$ which meet pairwise, so by Helly's Theorem 2.2 the intersection $W := C \cap \bigcap \mathscr{H}$ is a nonempty $G$-set.

We have to show that $W^\circ$ is orthogonal. Let $H_1, H_2$ be half spaces in $M$ which meet $W$ properly, i.e., both $W \cap H_i$ and $W \cap H_i^*$ are non empty. In particular neither $H_i$ nor $H_i^*$ is contained in $\mathscr{H}$, thus $|W \cap H_i| = \tfrac{1}{2}|W|$. Now we have
$$W \cap H_1 \subseteq W \cap H_2 \quad \Rightarrow \quad W \cap H_1 = W \cap H_2,$$
so $W \cap H_1$ and $W \cap H_2$ can only be nested half spaces of $W$ if $H_1 = H_2$ or $H_1 = H_2^*$. This proves that $W$ is a cube. □

In the case of simple median algebras we can be more precise.

**11.10. Proposition.**
*If $M$ is a simple, elliptic, discrete $G$-median algebra, then either $M$ is a finite cube or there exists a unique $G$-fixed point on $M$.*

*Proof.* By Theorem 11.9 we know that there exists a finite $G$-subcube $W \subseteq M$. If $W$ consists of more than one point, then some half space of $M$ cuts $W$. The hypothesis that $M$ is simple implies that every proper half space of $M$ cuts $W$, which means that in fact $M = W$.

If $W$ is a singleton then of course it comprises a fixed point for the $G$-action. Now suppose that we have two distinct fixed points $x, y \in M$, then $\overline{\Delta}(x, y)$ contains a $G$-orbit of hyper planes so, again, the fact that $M$ is simple implies that $M = [x, y]$. Next suppose that we have two distinct half spaces $H, J \in \Delta(x, y)$ with $H \parallel J$. We can't have $J^* = gH$, because $x \in gH \cap J$, so $J = gH$ for some $g \in G$, which implies that $g$ shifts $H$. But then we have an infinite set $\{\ldots, H, gH, g^2H, \ldots\} \subseteq \Delta(x,y)$ which contradicts the hypothesis that $M$ is discrete. Hence $M^\circ$ is orthogonal and $M$ is a finite cube.

Taking both cases together we see that $M$ is a finite cube unless there is a unique $G$-fixed point in $M$. □

### 11.11. Remark.

Suppose that $M$ is simple and infinite with a fixed point $x$. The arguments of the previous proof show that for every $y \in M$ the interval $[x,y]$ is a finite cube. In other words, $M$ is a star in the sense of § 8. □

Our fixed point theorem may also be seen as an analogue of the Tits Fixed Point Theorem for complete CAT(0) spaces (see DE LA HARPE-VALETTE[1989]). Following Serre we say that a group $G$ has $G$ has property (FH) if any isometric action of $G$ on a Hilbert space has bounded orbits, which by Tits' Theorem is the same as saying there exists a $G$-fixed point. It is well known that if $G$ has (FH), then any action of $G$ on a tree has bounded orbits — this is Serre's property (FA).

We now generalize this implication to discrete median algebras.

### 11.12. Proposition.

*If $G$ has property (FH) then every action of $G$ on a discrete median algebra is elliptic.*

*Proof.* Let $M$ be a discrete $G$-median algebra, choose a vertex $v \in M$ and let $\mathscr{H}$ denote the set of hyper planes of $M$. We will embed $M$ as a subset of the unit cube in $\mathbb{H} := \ell^2(\mathscr{H})$, the Hilbert space of square summable functions $\mathscr{H} \to \mathbb{R}$, such that the $G$-action on $M$ extends to an isometric action on $\mathbb{H}$. Furthermore, we will show that if $G$ acts on $M$ with unbounded orbits, then it also acts on $\mathbb{H}$ with unbounded orbits, which is the contrapositive of our proposition.

For $x \in M$ and $\overline{H} \in \mathscr{H}$ define

$$s_x(\overline{H}) := \begin{cases} 0 & \text{if } \overline{H} \notin \overline{\Delta}(v,x), \\ 1 & \text{if } \overline{H} \in \overline{\Delta}(v,x). \end{cases}$$

This function has finite support, hence it lies in $\mathbb{H}$. Observe that $s_v$ is the zero function.

For $g \in G$, $s \in \mathbb{H}$ and $\overline{H} \in \mathscr{H}$ we define

$$gs(\overline{H}) := \begin{cases} s(g^{-1}\overline{H}) & \text{if } \overline{H} \notin \overline{\Delta}(v,gv), \\ 1 - s(g^{-1}\overline{H}) & \text{if } \overline{H} \in \overline{\Delta}(v,gv). \end{cases}$$

The function $gs$ is again in $\mathbb{H}$, because the sum $\sum_{\overline{H} \in \mathscr{H}} s(\overline{H})^2$ differs from $\sum_{\overline{H} \in \mathscr{H}} gs(\overline{H})^2$ only in finitely many summands.

We have to check that our formula defines a $G$-action.

$$g_1(g_2 s(\overline{H})) = \begin{cases} g_2 s(g_1^{-1}\overline{H}) & \text{if } \overline{H} \notin \overline{\Delta}(v, g_1 v), \\ 1 - g_2 s(g_1^{-1}\overline{H}) & \text{if } \overline{H} \in \overline{\Delta}(v, g_1 v), \end{cases}$$

$$= \begin{cases} s(g_2^{-1} g_1^{-1}\overline{H}) & \text{if } \left(\overline{H} \in \overline{\Delta}(v, g_1 v) \Leftrightarrow g_1^{-1}\overline{H} \in \overline{\Delta}(v, g_2 v)\right), \\ 1 - s(g_2^{-1} g_1^{-1}\overline{H}) & \text{if } \left(\overline{H} \in \overline{\Delta}(v, g_1 v) \text{ xor } g_1^{-1}\overline{H} \in \overline{\Delta}(v, g_2 v)\right), \end{cases}$$

$$= \begin{cases} s((g_1 g_2)^{-1}\overline{H}) & \text{if } \overline{H} \notin \overline{\Delta}(v, g_1 g_2 v), \\ 1 - s((g_1 g_2)^{-1}\overline{H}) & \text{if } \overline{H} \in \overline{\Delta}(v, g_1 g_2 v), \end{cases}$$

$$= (g_1 g_2) s(\overline{H}).$$

It is easy to see that $gs_x = s_{gx}$, and for $s_1, s_2 \in \mathbb{H}$ we have $gs_1(\overline{H}) - gs_2(\overline{H}) = \pm(s_1(g^{-1}\overline{H}) - s_2(g^{-1}\overline{H}))$, hence

$$\|gs_1 - gs_2\| = \|s_1 - s_2\|,$$

i.e., $G$ acts by isometries on $\mathbb{H}$. This completes our proof. $\square$

We just emphasize one consequence of this proposition. For countable groups Serre's property (FH) coincides with Kazhdan's property (T), so we have the following.

**11.13. Corollary.**
*If $G$ is countable and has property (T), then every action of $G$ on a discrete median algebra is elliptic.*

Parabolic Actions

**11.14. Remark.**
If $M$ is a parabolic $G$-median algebra, then there exists a canonical $G$-fixed point for the action of $G$ on $M^{\circ\circ}$. Choose $x \in M$ and let

$$\xi := \{H \in M^\circ \mid H \cap Gx \text{ is unbounded}\}.$$

It is clear that exactly one of $H \cap Gx$ and $H^* \cap Gx$ are unbounded for every $H \in M^\circ$, and if we have $H_1, H_2 \in \xi$ with $H_1^* \geq H_2$ then both $H_1 \cap Gx$ and $H_2 \cap G_x \subseteq H^* \cap Gx$ would be unbounded, but then $H_1$ cuts properly. So if $M$ is parabolic then $\xi$ is an ultra filter which is fixed under the action of $G$. Also, $\xi$ can't lie in $ev(M)$, because $G$ has no fixed point in $M$. $\square$

As we have just observed, for any hyper plane $\overline{H}$ in a parabolic median algebra $M$ exactly one of $H \cap Gx$ and $H^* \cap Gx$ is unbounded; let $d_x \overline{H}$ denote the diameter of the bounded set. As $G$ acts by isometries we have

$$d_x \overline{H} = d_x g \overline{H} = d_{gx} \overline{H}.$$

**11.15. Lemma.**
If $M$ is parabolic, then $\{d_x\overline{H} \mid H \in M°\}$ is an unbounded subset of $\mathbb{N}$.

*Proof.* Let $N \in \mathbb{N}$. As $Gx \subseteq M$ is unbounded, we can find $g_1, g_2 \in G$ with $d(x, g_1x) \geq N$ as well as $d(x, g_2x) > 2d(x, g_1x)$. Suppose there exists an element $y \in [x, g_1x] \cap [g_2x, g_2g_1x]$, then
$$d(x, g_2x) \leq d(x, y) + d(y, g_2x) \leq d(x, g_1x) + d(x, g_1x),$$
contrary to our second assumption. Hence the intervals $[x, g_1x]$ and $[g_2x, g_2g_1x]$ are disjoint, and since they are convex there exists a halfspace $H \in M°$ separating them. By our first assumption we have $d_x\overline{H} \geq d(x, g_1x) \geq N$. □

**11.16. Proposition.**
If $M$ is discrete and parabolic, then the following hold.
  (i) There are infinitely many $G$-orbits of half spaces in $M$.
  (ii) $G$ can't be finitely generated.

*Proof.* (i) follows from the $G$-invariance of $d_x\overline{H}$.
(ii) Suppose that $G$ is generated by $g_1, \ldots, g_n$. Let $\mathscr{H} := \bigcup_{i=1}^n \overline{\Delta}(x, g_ix)$; this is a finite set of hyper planes. Now any $g \in G$ can be written as a product $g = g_{i_1} \ldots g_{i_k}$ of generators. Then
$$\overline{\Delta}(x, gx) = \overline{\Delta}(x, g_{i_1}x) + \overline{\Delta}(g_{i_1}x, g_{i_1}g_{i_2}x) + \cdots + \overline{\Delta}(g_{i_1} \ldots g_{i_{k-1}}x, gx) \subseteq G\mathscr{H}.$$
This implies that $d_x\overline{H} = 0$ unless $\overline{H} \in G\mathscr{H}$, so $d_x$ takes only finitely many values, which contradicts the previous lemma. □

**11.17. Example.**
Let $G$ be a group with an ascending sequence $G_1 < G_2 < \ldots$ of subgroups which exhaust $G$, i.e., $G = \bigcup_i G_i$. Then there exists a well known parabolic $G$-tree $T$. We give a direct description of it here.

The vertex set is $VT := \bigcup_i G/G_i$. We say that a vertex $gG_i$ belongs to the level $i$. For every $g \in G$ and $i \in \mathbb{N}$ there is an edge connecting $gG_i$ to $gG_{i+1}$. In particular, every vertex on level $i$ is connected to exactly one vertex on level $i+1$, hence $T$ is a forest. For any two vertices $g_1G_{i_1}$ and $g_2G_{i_2}$ there exists a number $j$ such that $g_1, g_2 \in G_j$. Then $g_1G_{i_1}$ and $g_2G_{i_2}$ are both connected by a path to the vertex $G_j$, thus $T$ is in fact a tree.

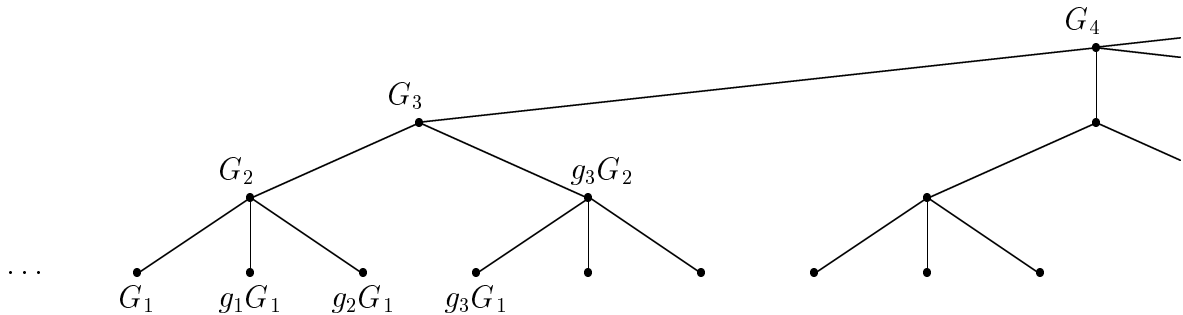

Corresponding to the edge from $gG_i$ to $gG_{i+1}$ we have two proper half spaces, viz. $H_{g,i} := \{hG_j \in VT \mid j \leq i \text{ and } hG_j \subseteq gG_i\}$ and its complement $H_{g,i}^*$. If we choose any vertex $x$ on level 1, then the diameter of $Gx \cap H_{g,i}$ equals $2i$, whereas $Gx \cap H_{g,i}^*$ is unbounded, hence $T$ is a parabolical $G$-tree. The non principal ultra filter $\xi = \{0^*, H_{g,i}^* \mid g \in G, i \in \mathbb{N}\}$ is fixed by $G$. □

Essential Actions

Suppose that $M$ is a discrete median algebra and $H \in M^\circ$ is a proper half space. To study the behaviour of the orbit $GH$ it is useful to consider a median algebra in which $G\overline{H}$ is the only orbit of hyper planes. This can be obtained as a quotient of $M$ by compressing all hyper planes not in the orbit of $\overline{H}$ in the following way.

Let $P := M^\circ$ and let $P_H$ denote the sub poc set generated by the orbit $GH$, i.e., $P_H = \{\varnothing, M, gH, gH^* \mid g \in G\}$. As $P \smallsetminus \{\varnothing, M\}$ is a discrete topological space (see Theorem 8.4), $P \smallsetminus P_H$ is an open subset, hence by Proposition 7.13 there exists a median algebra $M_H := M/\nabla(P \smallsetminus P_H)$ whose dual is isomorphic as a $G$-poc set to $P_H$. Let $q: M \to M_H$ denote the canonical quotient map; by construction this is a $G$-map. For $x, y \in M$ let $d_H(x, y) := |\overline{\Delta}(x, y) \cap G\overline{H}|$, then the median metric in $M_H$ is given by $d(q(x), q(y)) = d_H(x, y)$.

**11.18. Proposition.**
Let $M$ be a discrete median algebra and $H \in M^\circ$ a proper half space. Then one of the following holds.

  (i) $GH$ is finite.
  (ii) Either $\bigcap GH \neq \varnothing$ or $\bigcap GH^* \neq \varnothing$, in which case $H$ does not cut properly.
  (iii) $H$ cuts properly.

Moreover, $M_H$ is not elliptic iff (iii) holds.

*Proof.* Suppose that (i) does not hold, i.e., $M_H$ is infinite. If $x \in M_H$ is a $G$-fixed point, say $x \in q(H)$, then $q^{-1}(x) \subseteq gH$ for all $g \in G$, hence $q^{-1}(x) \subseteq \bigcap GH$. In particular for any $y \in \bigcap GH$ we have $Gy \cap H^* = \varnothing$, so $H$ does not cut properly.

If $M_H$ has no $G$-fixed point, then by Proposition 11.10, $G$ acts with unbounded orbits. By Proposition 11.16, $M_H$ can't be parabolic, hence $q(H)$ cuts $M_H$ properly. As the map $q$ shortens distances, it follows that also $H$ cuts $M$ properly. □

**11.19. Remark.**
This has an interesting consequence. Suppose that $H$ has an infinite $G$-orbit and does not cut properly. On the face of it, this means that for any $x \in M$ either $H \cap Gx$ or $H^* \cap Gx$ is bounded. But the proof of (ii) says that we can in fact find an $x \in M$ such that one of the two sets is empty.

We arrived at this fact at a somewhat circuitous route, thus it is instructive to give a direct proof as well. In case $G$ acts elliptically it follows easily from Theorem 11.9, hence we may assume that $G$ acts with unbounded orbits and the diameter of $H \cap Gx$ is bounded by $N$. We can find three points $x_1$, $x_2$ and $x_3$ in the orbit of $x$ with pairwise distance greater than $N$, e.g., by requiring that $d(x_1, x_2) > N$ and $d(x_1, x_3) > 2d(x_1, x_2)$. Now no translate $gH$ can contain more than one out of $x_1$, $x_2$ and $x_3$, therefore $m(x_1, x_2, x_3)$ does not lie in $gH$ for any $g \in G$. □

The next proposition explains the importance of essential half spaces, in particular it proves the equivalence of conditions (b) and (c) of Theorem 11.5. Observe that conditions (ii) and (iii) do not depend on the choice of $x \in M$.

**11.20. Proposition.**
Let $M$ be a discrete median algebra and $H \in M^\circ$ a proper half space. The following are equivalent.

(i) $H$ is essential, i.e., $GH$ is infinite and $H$ cuts properly.

(ii) For some $x \in M$ both sets $\{gH \mid x \in gH\}$ and $\{gH \mid x \notin gH\}$ are infinite.

(iii) For some $x \in M$ the set $\{d_H(x, gx) \mid g \in G\} \subseteq \mathbb{N}$ is unbounded.

(iv) $M_H$ is non elliptic.

*Proof.* (i) $\Rightarrow$ (iv). Suppose that $M_H$ is elliptic, then either condition (i) or (ii) of Proposition 11.18 holds, which means that $H$ is not essential.

(iv) $\Rightarrow$ (iii). If $M_H$ is non elliptic then the orbit $Gx$ is unbounded in the metric $d_H$.

(iii) $\Rightarrow$ (ii). Suppose that $\{gH \mid x \in gH\}$ contains $n$ elements, then so does $\{gH \mid kx \in gH\}$ for any $k \in G$. Now if $g\overline{H} \in \overline{\Delta}(x, kx)$ then either $x \in gH$ or $kx \in gH$, hence $\overline{\Delta}(x, kx) \cap G\overline{H}$ contains at most $2n$ elements.

(ii) $\Rightarrow$ (i). Clearly (ii) implies that $GH$ is infinite. Suppose that $H$ does not cut properly, then by Proposition 11.18 there exists a $y \in \bigcap GH$, say. It follows that
$$\{gH \mid x \notin gH\} = \Delta(y, x)$$
is finite. □

The connection between group actions on a discrete median algebra $M$ and almost invariant subsets of $G$ is based on the following pairing. For $x \in M$ and $H \in M^\circ$ we define
$$G[-,-] \colon M^\circ \times M \to \mathscr{P}G, \qquad G[H, x] := \{g \in G \mid gx \in H\}.$$
Recall that we have endowed the power set $\mathscr{P}G$ with various structures:

- A $\mathbb{F}_2 G$-bimodule, with $G$-actions by multiplication from the left and the right.
- A poc set, ordered by inclusion, with $*$ acting as complement.
- A median algebra with $m(x, y, z) = (x \cap y) \cup (y \cap z) \cup (z \cap x)$.

We have called $G[-,-]$ a pairing, because it is compatible with all these structures in the following sense.

**11.21. Proposition.**
Let $M$ be a $G$-median algebra. Let $g \in G$, $x \in M$, $H \in M^\circ$.

(i) $gG[H,x] = G[gH,x]$.

(ii) $G[H,x]g^{-1} = G[H,gx]$.

(iii) $G[-,x]\colon M^\circ \to \mathscr{P}G$ is a poc morphism.

(iv) $G[H,-]\colon M \to \mathscr{P}G$ is a median morphism.

*Proof.*
(i)  $h \in gG[H,x] \iff g^{-1}h \in G[H,x] \iff g^{-1}hx \in H \iff hx \in gH \iff h \in G[gH,x]$.

(ii)  $h \in G[H,x]g^{-1} \iff hg \in G[H,x] \iff hgx \in H \iff h \in G[H,gx]$.

(iii)  $g \in G[H^*,x] \iff gx \in H^* \iff gx \notin H \iff g \notin G[H,x] \iff g \in G[H,x]^*$.

Suppose that $H_1 \subseteq H_2$, then for all $g \in G$ we have $gx \in H_1 \Rightarrow gx \in H_2$, thus $G[H_1,x] \subseteq G[H_2,x]$.

(iv)  $g \in G[H, m(x,y,z)] \iff gm(x,y,z) \in H \iff m(gx,gy,gz) \in H$
$\iff (gx \in H \wedge gy \in H) \vee (gy \in H \wedge gz \in H) \vee (gy \in H \wedge gx \in H)$
$\iff g \in m(G[H,x], G[H,y], G[H,z])$. □

**11.22. Proposition.**
Let $M$ be a median $G$-algebra. Let $x \in M$, $H \in M^\circ$ and $g_1, g_2 \in G$. Then
$$G[H,x]g_1 + G[H,x]g_2 = \{g \in G \mid g^{-1}\overline{H} \in \overline{\Delta}(g_1^{-1}x, g_2^{-1}x)\}.$$

*Proof.*  $g \in G[H,x]g_1 + G[H,x]g_2 \iff g \in G[H,g_1^{-1}x] + G[H,g_2^{-1}x]$
$\iff gg_1^{-1}x \in H \;\; \text{xor} \;\; gg_2^{-1}x \in H \iff g_1^{-1}x \in g^{-1}H \;\; \text{xor} \;\; g_2^{-1}x \in g^{-1}H$
$\iff g^{-1}\overline{H} \in \overline{\Delta}(g_1^{-1}x, g_2^{-1}x)$. □

*Proof of the implication* (b) $\Rightarrow$ (a) *of Theorem* 11.6.
Let $M$ be a simple, non elliptic, discrete $G$-median algebra, choose a point $x \in M$ and a proper half space $H \in M^\circ$. Let $K := G_H = \{g \in G \mid gH = H\}$ be the stabilizer of $H$ and $\overline{K} := G_{\overline{H}} = \{g \in G \mid gH = H \text{ or } gH = H^*\}$ the stabilizer of the hyper plane defined by $H$, then $\overline{K}$ contains $K$ as a subgroup of index at most 2. Finally we set $A := G[H,x]$.
In view of Proposition 11.21(i) we have $KA = A$, so $A$ satisfies (**End 1**).
Since $M$ is simple and discrete, for any $g \in G$ the set $\overline{\Delta}(x, g^{-1}x) = \{f_1\overline{H}, \ldots, f_k\overline{H}\}$ is finite, thus by Proposition 11.22 $G[H,x] + G[H,x]g = \overline{K}\{f_1^{-1}, \ldots, f_k^{-1}\}$ is $K$-finite, in other words $A$ satisfies (**End 2**).
Notice that $H$ is essential, so by Proposition 11.20(iii) neither $A$ nor $G \smallsetminus A$ is $K$-finite, which means that $A$ satisfies (**End 3**). □

For the converse implication we start with a subgroup $K < G$ and a subset $A \subset G$ satisfying (**End 1-3**). Our goal is to construct a simple, discrete $G$-median algebra $M$. In fact we already know the dual of $M$; it is the poc set $P := \{\varnothing, G, gA, gA^* \mid g \in G\}$ with the Boolean poc structure. Observe that the canonical representation $P \subset \mathscr{P}G$ is discrete. The reason for this is the following lemma, which is the analogue of Proposition 11.22.

**11.23. Lemma.**
For all $g, g_1, g_2 \in G$ we have
$$\overline{gA} \in \overline{\Delta}(g_1, g_2) \quad \Leftrightarrow \quad g^{-1} \in Ag_1^{-1} + Ag_2^{-1}. \qquad \square$$

Let $K' := G_A = \{g \in G \mid gA = A\}$. In view of (**End 1**) we have that $K \leq K'$, and (**End 2**) together with (**End 3**) implies that $K'$ can't be too large.

**11.24. Lemma.**
$K$ has finite index in $K'$.

*Proof.* Since $A$ is $K$-proper, there exist elements $x \in A$ and $y \in G \smallsetminus A$. So for $g := x^{-1}y$ we have $A \neq Ag$. The set $A + Ag$ is $K$ finite and non empty, so it contains a coset $K'h$. It follows that $K'$ is $K$-finite, which is another way of saying that $K$ has finite index in $K'$. $\qquad \square$

*Proof of the implication* (a) $\Rightarrow$ (b) *of Theorem 11.6.*
Now we can argue as in Proposition 9.8. We have a canonical map $\iota \colon G \to P^\circ$ with $\iota(g) := \{a \in P \mid g \in A\}$. This is obviously a $G$ map, so $g_1\iota(g_2) = \iota(g_1g_2)$. The previous lemma shows that the entire orbit $G\chi_1$ is contained in an almost equivalence class $M$ of $P^\circ$, moreover $M^\circ$ is isomorphic to $P$. This $M$ is discrete and simple by construction, and it contains a half space $H$ and a point $x = \iota(1)$ such that $A = G[H, x]$ and $G_H = K'$. The fact that $A$ satisfies condition (**End 3**) implies that $H$ satisfies Proposition 11.20(ii), but then $M$ must be non elliptic. This concludes the proof of Theorem 11.5. $\qquad \square$

*Proof of Theorem 11.6.*
The implication (a) $\Rightarrow$ (b) is a weak form of Theorem 11.5.
Suppose that $G$ admits a non elliptic action on a discrete median algebra $M$. Choose a finite set $\{g_1, \ldots, g_n\}$ of generators of $G$ and an element $x \in M$. The set
$$\mathscr{H} := \{g\overline{H} \mid g \in G \text{ and } \overline{H} \in \overline{\Delta}(x, g_ix) \text{ for some } i \in \{1, \ldots, n\}\}$$
is a union of finitely many $G$-orbits of hyper planes; choose a transversal $\overline{H}_1, \ldots, \overline{H}_k$ to the action of $G$ on $\mathscr{H}$.
Recall that the distance $d(x, gx)$ is computed by counting the hyper planes in $\overline{\Delta}(x, gx)$. The proof of Proposition 11.16 shows that $\overline{\Delta}(x, gx) \subseteq \mathscr{H}$, hence we can write
$$d(x, gx) = d_{H_1}(x, gx) + \ldots + d_{H_k}(x, gx).$$

By the hypothesis of Theorem 11.6 the function $g \mapsto d(x, gx)$ is unbounded. Therefore also $g \mapsto d_{H_i}(x, gx)$ must be unbounded for some $i \in \{1, \ldots, k\}$, which means that $H_i$ is essential and by Theorem 11.5 we have $e(G, G_{H_i}) > 1$. □

**Remark.**

We record two consequences of the above proofs.

1) If $G$ acts non elliptically on a discrete median algebra $M$ such that $M°$ consists of finitely many $G$-orbits of half spaces, then the action is essential.

2) Suppose $G$ is finitely generated and acts on a discrete median algebra $M$. Then for any $x \in G$ the $G$-subalgebra $M' = \operatorname{conv}(Gx)$ has only a finite number of $G$-orbits of hyper planes. □

For a median algebra $M$ that is simple with respect to a group $G$ one can characterize the subgroups $L$ of $G$ that act elliptically or essentially in the following way.

**11.25. Proposition.**

*Let $M$ be a simple, discrete $G$-median algebra. Let $H \subset M$ be a proper half space with stabilizer $K = G_H$, choose a point $x \in M$ and let $A := G[H, x]$.*

  (i) *A subgroup $L \leq G$ acts elliptically on $M$ iff there exists a subset $A' \subset G$ such that $A + A'$ is $K$-finite and $A' = A'L$.*

  (ii) *The half space $g^{-1}H$ is essential with respect to the action of $L$ iff neither $A \cap KgL$ nor $A^* \cap KgL$ is $K$-finite.*

*Proof.* (i) Suppose that $L$ acts elliptically, then there exists a finite $L$-cube $W \subseteq M$. Choose an element $w \in W$, then $L_w$ has finite index in $L$. By Proposition 11.21 the set $A_w := G[H, w]$ satisfies $A_w = A_w L_w$ and the symmetric difference $A_w + A$ is $K$-finite. Let $\{l_1, \ldots l_n\}$ be a transversal to $L_w/L$, then the set

$$A' = A_w L = A_w l_1 \cup \ldots \cup A_w l_n$$

lies in the $K$-almost equality class of $A$. Note that $A' = \{g \in G \mid gW \cap H \neq \varnothing\}$. Conversely, suppose we are given a set $A'$ with the above properties, in particular suppose that $A + A'$ is contained in a union of $N$ cosets of $\overline{K}$, the stabilizer of $\overline{H}$. Consider the right action of $G$ on hyper planes defined by $\overline{H}g := g^{-1}\overline{H}$. For any $l \in L$ we have

$$\overline{\Delta}(x, l^{-1}x) = \overline{H}(A + Al) = \overline{H}(A + A') + \overline{H}(A + A')l.$$

It follows that $d(x, lx) \leq 2N$, which means that $M$ is elliptic as an $L$-median algebra.

(ii) follows from Proposition 11.20. □

**Notes.**

The notation $G[H, x]$ for the pairing is borrowed from DICKS 1980. Theorem 11.5 is an extended version of ROLLER-NIBLO [1998], which in turn is based on SAGEEV [1995]. Sageev's original definition of essentiality is analogous to condition (ii) of Proposition 11.20. The embedding of a median algebra into a Hilbert Space used in the proof of Proposition 11.12 also originates in ROLLER-NIBLO [1998]. GERASIMOV [1997] contains a proof of Theorem 11.6 using different methods. The second half of the proof of Theorem 11.9 uses an idea of BANDELT-VAN DE VEL [1987].

# Open Problems

### CUBIC SPACES

Recall the cube $\mathscr{F}X$, where the distance between two points $x, y \in \mathscr{F}X$ is given by $|x + y|$. A (pseudo) metric space $(M, d)$ is called cubic if for some set $X$ there exists a a map $\iota \colon M \to \mathscr{F}X$ which preserves the metric. Many of our arguments carry over to cubic spaces, for example one can study the median algebra they generate inside $\mathscr{F}X$.

We observed that $d$ has to satisfy the pattern condition

$$\forall x, y, z \in M \colon \quad d(x, y) + d(y, z) + d(z, x) \text{ is even.}$$

Is this condition sufficient to show that M is cubic? This may be known to people who study embeddings of metric spaces into $\ell^1$-spaces.

### DISCRETE POC SETS

Characterize discrete poc sets which are not duals of discrete median algebras. It is conceivable that there is a (finite?) list of minimal examples.

### DISCRETE MEDIAN ALGEBRAS ARISING FROM SUBGROUPS

Given a group $G$ and a subgroup $K$ with $e(G, K) > 1$ we construct a simple, essential, discrete median algebra $M$ with an action of $G$. If $G$ is finitely generated then $M$ has finitely many orbits of hyper planes. It would be interesting to know, under which conditions on $G$ and $K$ the dimension of $M$ is finite or at least $\omega$, furthermore, when the $G$-action is cocompact. SAGEEV[1996] has shown that both is the case when $G$ is hyperbolic and $K$ is a quasi convex subgroup. RUBINSTEIN-WANG [1996] have given a geometric example, where $G$ and $K$ are finitely presented, but $M$ contains infinite dimensional cubes. The main point of that example is that $G$ is a finite union of double cosets $KgK$.

### COXETER GROUPS

A Coxeter group $G$ acts on a vector space such that every involution in $G$ acts as a reflection in a hyper plane. This gives rise to a poc set with a discrete representation, and hence to a discrete median algebra on which $G$ acts. One can show that the median algebra has finite dimension when $G$ is finitely generated. It is locally finite at the points which generate the median algebra; that is a consequence of the Finiteness Theorem of BRINK-HOWLETT [1993]. A global bound on the size of the links of vertices would imply a strong Parallel Walls Theorem for Coxeter groups. For which Coxeter groups is the median algebra cocompact?

## Property (FH)

We have shown that a group $G$ with property (FH) has no subgroup $K$ with $e(G,K) > 1$. What about the converse?

## Median Groups

A group $G$ is called median if it acts freely and transitively on a median algebra. This is equivalent to saying that the Cayley graph with respect to a certain set of generators is a cubing. For example free groups and free abelian groups are median. BASARAB [1997] studies a class of median groups that is a mixture of the two classes. It is easy to see that a median group has a presentation containing relations of length 2 and 4 only (corresponding to edges and squares in the star of 1 in $G$). Characterize all presentations of median groups!

## Splitting Conjecture

KROPHOLLER-ROLLER [1989] contains the following conjecture: Let $G$ be a group with a subgroup $K$ and a subset $A$ which is proper $K$-almost invariant and satisfies $A = KAK$, then $G$ splits algebraically over some subgroup $K'$ which is commensurable to a subgroup of $K$.

DUNWOODY-ROLLER [1993] and SAGEEV [1996] contains some progress with that conjecture.

## Almost Stability Theorem

Is there a version of the Almost Stability Theorem of DICKS-DUNWOODY [1989] for group actions on trees with arbitrary (i.e., not necessarily finite) edge stabilizers?

# Index